\documentclass[11pt]{article}
\usepackage{amssymb}
\usepackage{amsmath}
\usepackage{amsthm}
\newcommand{\cA}{{\cal A}}
\newcommand{\cB}{{\cal B}}
\newcommand{\cC}{{\cal C}}

\newcommand{\cH}{{\cal H}}
\newcommand{\cE}{{\cal E}}

\newcommand{\cI}{{\cal I}}
\newcommand{\cO}{{\cal O}}
\newcommand{\cL}{{\cal L}}

\newcommand{\cN}{{\cal N}}
\newcommand{\cF}{{\cal F}}

\newcommand{\cP}{{\cal P}}
\newcommand{\cQ}{{\cal Q}}
\newcommand{\cS}{{\cal S}}

\newcommand{\cV}{{\cal V}}
\newcommand{\cW}{{\cal W}}

\newcommand{\cY}{{\cal Y}}
\newcommand{\cZ}{{\cal Z}}

\renewcommand{\AA}{{\mathbb A}}

\newcommand{\NN}{{\mathbb N}}
\newcommand{\ZZ}{{\mathbb Z}}

\newcommand{\CC}{{\mathbb C}}
\newcommand{\PP}{{\mathbb P}}
\newcommand{\OO}{{\mathbb O}}

\newcommand{\gp}{\mathfrak{p}}
\newcommand{\gq}{\mathfrak{q}}

\newcommand{\on}{\operatorname}

\newcommand{\Rep}{{\on{Rep}}}

\newcommand{\Qlb}{\mathbb{\bar Q}_\ell}
\newcommand{\Gm}{\mathbb{G}_m}
\newcommand{\Ga}{\mathbb{G}_a}
\newcommand{\A}{\mathbb{A}}

\newcommand{\toup}[1]{\stackrel{#1}{\to}}
\newcommand{\hook}[1]{\stackrel{#1}{\hookrightarrow}}

\newcommand{\getsup}[1]{\stackrel{#1}{\gets}}
\newcommand{\Sp}{\on{\mathbb{S}p}}

\newcommand{\GSp}{\on{G\mathbb{S}p}}

\newcommand{\CT}{\on{CT}}
\newcommand{\Hom}{\on{Hom}}

\newcommand{\Ext}{\on{Ext}}

\newcommand{\Sym}{\on{Sym}}

\newcommand{\Ker}{\on{Ker}}

\newcommand{\Pic}{\on{Pic}}

\newcommand{\Bun}{\on{Bun}}

\newcommand{\Bunt}{\on{\widetilde\Bun}}

\newcommand{\rk}{\on{rk}}
\newcommand{\Spec}{\on{Spec}}

\newcommand{\supp}{\on{supp}}

\newcommand{\Gr}{\on{Gr}}

\newcommand{\GL}{\on{GL}}

\newcommand{\Eis}{{\on{Eis}}}
\newcommand{\pr}{\on{pr}}
\newcommand{\id}{\on{id}}

\newcommand{\QED}{$\square$} 
\newcommand{\Fq}{\mathbb{F}_q}  
\newcommand{\Fp}{\mathbb{F}_p}  
\newcommand{\iso}{{\widetilde\to}}

\newcommand{\comp}{\circ}
\newcommand{\Four}{\on{Four}}
\renewcommand{\H}{{\on{H}}}   
\newcommand{\R}{\on{R}\!}   

\newcommand{\DD}{\mathbb{D}}  
\newcommand{\D}{\on{D}}       

\newcommand{\select}[1]{{\it{#1}}}

\newcommand{\ov}[1]{\overline{#1}}

\renewcommand{\div}{\on{div}}

\renewcommand{\P}{{\on{P}}}

\newcommand{\<}{\langle}
\renewcommand{\>}{\rangle}
\newcommand{\trian}{\vartriangle}

\newcommand{\ev}{\on{ev}}

\newcommand{\Sph}{\on{Sph}}

\newcommand{\tboxtimes}{\tilde\boxtimes}

\newcommand{\act}{\on{act}}
\newcommand{\dimrel}{\on{dim.rel}}

\newcommand{\SL}{\on{SL}}
\newcommand{\RCov}{\on{RCov}}
\newcommand{\Cov}{\on{Cov}}
\newcommand{\diag}{\on{diag}}

\newcommand{\Mat}{\on{Mat}}
\newcommand{\Wald}{\on{{\cal W}ald}}

\newtheorem{Lm}{Lemma}
\newtheorem{Th}{Theorem}
\newtheorem{Pp}{Proposition}
\newtheorem{Cor}{Corolary}

\theoremstyle{remark}
\newtheorem{Rem}{Remark}
\newtheorem{Rems}{Remarks}

\theoremstyle{definition}
\newtheorem{Def}{Definition}

\newenvironment{Prf}{\par\noindent {\it Proof }}{\QED}
\newcommand{\Step}[1]{\par\noindent{\bf Step {#1}}.}

\begin{document}

\author{Sergey Lysenko}
\title{Whittaker and Bessel functors for $\GSp_4$}
\date{}
\maketitle
\begin{abstract}
\noindent{\scshape Abstract}\hskip 0.8 em 
The theory of Whittaker functors for $\GL_n$ is an essential technical tools in Gaitsgory's proof of the Vanishing Conjecture appearing in the geometric Langlands correspondence (\cite{G}). 
We define Whittaker functors for $\GSp_4$ and study their properties. These functors correspond to the maximal parabolic subgroup of $\GSp_4$, whose unipotent radical is not commutative. 
 
  We also study similar functors corresponding to the Siegel parabolic subgroup of $\GSp_4$, they are related with Bessel models for $\GSp_4$ and Waldspurger models for $\GL_2$.
 
  We define the Waldspurger category, which is a geometric counterpart of the Waldspurger module over the Hecke algebra of $\GL_2$. We prove a geometric version of the multiplicity one result for the Waldspurger models.
\end{abstract} 

\medskip\smallskip

{\centerline{\scshape 1. Introduction}}

\medskip\smallskip\noindent
1.1 Whittaker and Bessel models are of importance in the theory of automorphic representations of $\GSp_4$. This paper is the first in a series of two, where we study some phenomena corresponding to these models in the geometric Langlands program. 

  The theory of Whittaker functors for $\GL_n$ is an essential technical tool in Gaitsgory's proof of the Vanishing Conjecture appearing in the geometric Langlands correspondence (\cite{G}). First part of our results is an analog of this theory for  $\GSp_4$. 

 Let us first remind some facts about automorphic forms on $G=\Sp_4$. Let $X$
be a smooth projective absolutely irreducible curve over $\Fq$, $F=\Fq(X)$ and $\AA$ be the
adeles ring of $F$. Let $B$ be a Borel subgroup of $G$ and $U\subset B$ its unipotent
radical. For a character $\psi: U(F)\backslash U(\AA)\to\CC^*$ one has a global Whittaker
module over $G(\AA)$
$$
WM_{\psi}=\{f:U(F)\backslash G(\AA)\to\CC\mid\; f(ug)=\psi(u)f(g)\;\mbox{for}\; u\in
U(\AA),\; f\; \mbox{is smooth}\}
$$

 Let $\cA_{cusp}(G(F)\backslash G(\AA))$ be the space of cusp forms on
$G(F)\backslash G(\AA)$. The usual Whittaker operator $W_{\psi}: \cA_{cusp}(G(F)\backslash
G(\AA))\to WM_{\psi}$ is given by
$$
W_{\psi}(f)(g)=\int_{U(F)\backslash U(\AA)} f(ug)\psi(u^{-1})du,
$$
where $du$ is induced from a Haar measure on $U(\AA)$. Whence for $\GL_n$ (and generic
$\psi$) the operator $W_{\psi}$ is an injection, this is not always the case for more
general groups. There are cuspidal automorphic representations of $\Sp_4$ that
don't admit a
$\psi$-Whittaker model for any $\psi$.

 Recall that $\cA_{cusp}(G(F)\backslash G(\AA))$ decomposes as a direct sum
\begin{equation}
\label{decomp_1}
\cA_{cusp}(G(F)\backslash G(\AA))=I_3(H_2)\oplus I_4(H_3)\oplus I_5(H_4)
\end{equation}
in the notation of (\cite{R}, Sect. 1.3, p. 359), the summands being
$G(\AA)$-invariant\footnote{In \select{loc.cit.} $F$ is a number field, but
(\ref{decomp_1}) holds also over function fields.}. The decomposition is orthogonal with
respect to the scalar product
\begin{equation}
\label{scalar_product_classical}
\<f,h\>=\int_{G(F)\backslash G(\AA)} f(x)\ov{h(x)}dx,
\end{equation}
where $dx$ is induced from a Haar measure on $G(\AA)$. 

 For any $f\in \cA_{cusp}(G(F)\backslash G(\AA))$ its $\theta$-lifting to $\OO(2)(\AA)$ vanishes (\select{loc.cit.}, Corolary~2 to Theorem~I.2.1). By definition, $I_4(H_3)\oplus I_5(H_4)$ (resp., $I_5(H_4)$) are those cuspidal forms whose $\theta$-lifting to $\OO_4(\AA)$ (resp., to $\OO_4(\AA)$ and $\OO_6(\AA)$) vanishes. Here $\OO_{2r}$ is the orthogonal group defined by the hyperbolic quadratic form in a $2r$-dimensional space. 
 
 The space $I_5(H_4)$ is also the intersection of kernels of $W_{\psi}$ for all $\psi$. It is
known as the space of hyper-cuspidal forms on $G(F)\backslash G(\AA)$ (\cite{PSh},
Definition on p. 328). Another description of $I_5(H_4)$ is as follows. Let $P_1\subset G$
be the parabolic preserving a 1-dimensional isotropic subspace in the standard
representation $V$ of $G$, $U_1\subset P_1$ be its unipotent radical, $U_0$ the center of
$U_1$. Then $f\in \cA_{cusp}(G(F)\backslash G(\AA))$ lies in $I_5(H_4)$ if and only if
$$
\int_{U_0(F)\backslash U_0(\AA)} f(ug)du=0
$$
for all $g\in G(\AA)$.
 
  If $V'\subset V$ is a 2-dimensional subspace such that the symplectic form on $V$ restricts to a non degenerate form on $V'$ then let $H\subset G=\Sp(V)$ be the subgroup of those $g\in G$ that preserve and act trivially on $V'$. Then $f\in \cA_{cusp}(G(F)\backslash G(\AA))$ lies in 
$I_4(H_3)\oplus I_5(H_4)$ if and only if
$$
\int_{H(F)\backslash H(\AA)} f(hg)dh=0
$$
for all $g\in G(\AA)$ (\select{loc.cit.}, Section~3). Note that $H\,\iso\,\SL_2$.   
  
\medskip\noindent 
1.2 In the geometric setting we work with $G=\GSp_4$ (over an algebraically closed field
of characteristic $p>2$). For a scheme (or a stack $S$) write $\D(S)$ for the derived category of $\ell$-adic \'etale sheaves on $S$.  

 Let $\Bun_G$ be the stack of $G$-bundles on $X$. Inside of the triangulated category
$\D_{cusp}(\Bun_G)$ of cuspidal sheaves on $\Bun_G$ we single out a full triangulated
subcategory $\D_{hcusp}(\Bun_G)$ of hyper-cuspidal sheaves. Both they are preserved by
Hecke functors. So, a natural step in the geometric Langlands program for $G$ is to
understand the Hecke action on $\D_{hcusp}(\Bun_G)$ and on
$\D_{cusp}(\Bun_G)/\D_{hcusp}(\Bun_G)$. 

 The category $\D_{cusp}(\Bun_G)$ is equiped with the `scalar product', which is an
analogue of (\ref{scalar_product_classical}), it sends $K_1,K_2$ to $\R\Hom(K_1,K_2)$. The
(left and right) orthogonal complements $^{\perp}\D_{hcusp}(\Bun_G)$,
$\D_{hcusp}(\Bun_G)^{\perp}\subset \D_{cusp}(\Bun_G)$ are also preserved by Hecke functors.

\smallskip

 A $G$-bundle on $X$ is a triple: a rank 4 vector bundle $M$ on $X$, a line bundle $\cA$ on $X$, and a symplectic form $\wedge^2 M\to\cA$. Let $\alpha:\bar\cQ_1\to\Bun_G$ be the stack over $\Bun_G$ whose fibre over $(M,\cA)$ consists of all nonzero maps of coherent sheaves $\Omega\hook{} M$, where $\Omega$ is the canonical line bundle on $X$. 

 We introduce the notion of cuspidality and hyper-cuspidality on $\bar\cQ_1$, thus leading
to full triangulated subcategories
$\D_{hcusp}(\bar\cQ_1)\subset\D_{cusp}(\bar\cQ_1)\subset \D(\bar\cQ_1)$. 

\smallskip

 Then we describe $\D_{cusp}(\bar\cQ_1)/\D_{hcusp}(\bar\cQ_1)$ in terms
of geometric Whittaker models. Namely, we introduce a stack $\bar\cQ$ (it was
denoted by $\cY$ in \cite{Ly}) and a full triangulated subcategory $\D^W(\bar\cQ)\subset \D(\bar\cQ)$.
Our $\D^W(\bar\cQ)$ is a geometric analog of the space $WM_{\psi}$. 
  
  We define Whittaker functors that give rise to an equivalence of triangulated categories
\begin{equation}
\label{equvalence_intro}
W: \D_{cusp}(\bar\cQ_1)/\D_{hcusp}(\bar\cQ_1)\,\iso\, \D^W(\bar\cQ)
\end{equation}

The Hecke functor $\H^{\gamma}$ corresponding to the standard representation of the Langlands dual group $\check{G}\,\iso\, \GSp_4$ acts on all the categories mentioned in Sect.~1.2. Moreover, the equivalence (\ref{equvalence_intro}) commutes with $\H^{\gamma}$. The restriction functor 
$$
\alpha^*: \D_{cusp}(\Bun_G)/\D_{hcusp}(\Bun_G)\to
\D_{cusp}(\bar\cQ_1)/\D_{hcusp}(\bar\cQ_1)
$$ 
also commutes with $\H^{\gamma}$. As in the case of $\GL_n$ (\cite{G}, Theorem~7.9), the advantage of $\bar\cQ$ over $\Bun_G$ is that the functor $\H^{\gamma}: \D(\bar\cQ)\to\D(X\times\bar\cQ)$ is right-exact for the perverse t-structures.  

  The essential difference with $\GL_n$ case is that the Whittaker functor $W:\D(\bar\cQ_1)\to\D^W(\bar\cQ)$ is not exact for the perverse t-structures. We can only indicate  full triangulated subcategories $\D^W_{cusp}(\bar\cQ_1)\subset \D_{cusp}(\bar\cQ_1)\subset \D(\bar\cQ_1)$ such that the restriction of $W$ yields an equivalence 
$$
\D^W_{cusp}(\bar\cQ_1)\,\iso\, \D^W(\bar\cQ)
$$ 
of triangulated categories. Then (\ref{equvalence_intro}) follows from the fact that the natural inclusion functor induces an equivalence of triangulated categories
$$
\D^W_{cusp}(\bar\cQ_1)\,\iso\, \D_{cusp}(\bar\cQ_1)/\D_{hcusp}(\bar\cQ_1)
$$ 
This is the content of Sect.~2-6.  
 
\smallskip\noindent
1.3  The stack $\bar\cQ_1$ corresponds to the parabolic subgroup $P_1\subset G$. In Sect.~7 we define functors similar to the Whittaker ones for the Siegel parabolic subgroup $P\subset G$. They are related to Bessel models\footnote{Bessel models will be studied in the second paper of this series.} for $G$. 

 The general idea behind is that various Fourier coefficients of automorphic sheaves carry additional structure coming from the action of Hecke operators. 
 
  Let $\alpha_{\cZ}: \cZ_1\to\Bun_G$ be the stack whose fibre over $(M,\cA)$ is the scheme of isotropic subsheaves $L_2\subset M$, where $L_2$ is a locally free $\cO_X$-module of rank 2. 
The open substack $\Bun_P\subset \cZ_1$ is given by the condition that $L_2$ is a subbundle. 
Then $\Bun_P$ is the stack classifying: a rank 2 bundle $L_2$ on $X$, a line bundle $\cA$ on $X$, and an exact sequence $0\to\Sym^2 L_2\to ?\to\cA\to 0$.    
 
  Let $\cS_{ex}$ denote the stack classifying: a rank 2 vector bundle $L_2$ on $X$, a line bundle $\cA$ on $X$, and a map $\Sym^2 L_2\to \cA\otimes\Omega$. Write $\Bun_i$ for the stack of rank $i$ vector bundles on $X$. Then $\Bun_P$ and $\cS_{ex}$ are dual (generalized) vector bundles over $\Bun_2\times\Bun_1$, so we have the Fourier transform functor $\Four:\D(\Bun_P)\to \D(\cS_{ex})$.  
 
  For a complex $K\in\D(\Bun_G)$ its Fourier coefficient with respect to the Siegel parabolic is, by definition, $F_{\cS_{ex}}(K)=\Four(K\mid_{\Bun_P})$. If $K$ is a Hecke eigen-sheaf on $\Bun_G$ then $F_{\cS_{ex}}(K)$ satisfies some additional property (cf. Proposition~\ref{Pp_13_property}), which is a consequence of the following result.
  
   Let $\cZ_{2,ex}\to\cZ_1$ be the stack whose fibre over a point $(L_2\subset M, \cA)\in\cZ_1$ is the space $\Hom(\Sym^2 L_2, \cA\otimes\Omega)$. We define a full triangulated subcategory $\D^W(\cZ_{2,ex})\subset \D(\cZ_{2,ex})$ singled out by some equivariance condition. Then we establish an equivalence of triangulated categories
$$
WZ: \D(\cZ_1)\,\iso\, \D^W(\cZ_{2,ex}),
$$   
which is exact for perverse t-structures. The Hecke functor $\H^{\gamma}$ acts on both categories and commutes with this equivalence. Our $\D^W(\cZ_{2,ex})$ is a way to think about the Fourier coefficients $F_{\cS_{ex}}(K)$ \select{together with an action of Hecke operators}.  
 
  One also has a notion of hyper-cuspidality on $\cZ_1$ and $\cZ_{2,ex}$ leading to full triangulated subcategories $\D_{hcusp}(\cZ_1)\subset \D(\cZ_1)$ and $\D^W_{hcusp}(\cZ_2)\subset \D^W(\cZ_{2,ex})$ preserved by $H^{\gamma}$. The functor $WZ$ induces an equivalence 
$$
\D_{hcusp}(\cZ_1)\,\iso\, \D^W_{hcusp}(\cZ_2)
$$  
A complex $K\in\D(\Bun_G)$ is hyper-cuspidal if and only if $\alpha_{\cZ}^*K$ is hyper-cuspidal.
 
\medskip\noindent
1.4 In a sense, Bessel models for $G$ is a way to think about the Fourier coefficients $F_{\cS_{ex}}(K)$ of automorphic sheaves $K\in\D(\Bun_G)$ in terms of the Waldspurger models 
for $\GL_2$ (\cite{BFF}). This is our motivation for the study of these Walspurger models in Sect.~8, which is independent of the rest of this paper. 

 The following background result is due to Waldspurger (\cite{W}, Lemma~8). Set $F=\Fq((t))$ and $\cO=\Fq[[t]]$. Let $\tilde F$ be an \'etale $F$-algebra with $\dim_F(\tilde F)=2$ such that $\Fq$ is algebraically closed in $\tilde F$. Let $\tilde\cO$ be the integral closure of $\cO$ in $\tilde F$. We have two cases:
\begin{itemize}
\item $\tilde F\,\iso\, \Fq((t^{\frac{1}{2}}))$ \  (the nonsplit case)
\item $\tilde F\,\iso\, F\oplus F$  \  (the split case)
\end{itemize} 
 
  Write $\GL(\tilde F)$ for the automorphism group of the $F$-vector space $\tilde F$, and $\GL(\tilde\cO)\subset\GL(\tilde F)$ for the stabilizor of $\tilde\cO$. Fix a nonramified character
$\chi:\tilde F^*/\tilde\cO^*\to \Qlb^*$. Denote by $\chi_c: F^*/\cO^*\to\Qlb^*$ the restriction of $\chi$. The \select{Waldspurger module} is the vector space
\begin{multline*}
W\! A_{\chi}=\{f: \GL(\tilde F)/\GL(\tilde\cO)\to\Qlb\mid f(ux)=\chi(u)f(x)\;\;\mbox{for}\; u\in\tilde F^*,\\   f\; \mbox{is of compact support modulo} \;\;\tilde F^*\}
\end{multline*}
The Hecke algebra
\begin{multline*}
\H_{\chi_c}=\{h: \GL(\tilde\cO)\backslash \GL(\tilde F)/\GL(\tilde\cO)\to\Qlb\mid h(ux)=\chi_c(u)h(x)\;\;\mbox{for}\;\; u\in F^*,\\
h\;\;\mbox{is of compact support}\}
\end{multline*}
acts on $W\! A_{\chi}$ via
$$
h\in\H_{\chi_c}, \; f\in W\! A_{\chi} \;\to \; (h\ast f)(g)=\int_{\GL(\tilde F)} h(x)f(gx^{-1})dx,
$$ 
where $dx$ is the Haar measure of $\GL(\tilde F)$ such that the volume of $\GL(\tilde\cO)$ is one. Then $W\! A_{\chi}$ is \select{a free module of rank one} over $\H_{\chi_c}$ (mutilpicity one for Waldspurger model). 
 
  We prove a categorical version of this. Namely, the affine grassmanian $\Gr_{\tilde F}:=\GL(\tilde F)/\GL(\tilde\cO)$ can be viewed as an ind-scheme over $\Fq$ equiped with an action of the group scheme $\tilde F^*$. Pick a 1-dimensional $\Qlb$-vector space $\tilde E_{\tilde x}$ for each $\tilde x\in\Spec\tilde F$. We introduce \select{Waldspurger category} of those $\tilde\cO^*$-equivariant perverse sheaves on $\Gr_{\tilde F}$ that change under the action of each uniformizor $t_{\tilde x}\in\tilde F^*/\tilde\cO^*$ by $\tilde E_{\tilde x}$ (for each $\tilde x\in\Spec \tilde F$). This is a geometric counterpart of $W\! A_{\chi}$. 
  
  The nonramified Hecke algebra for $\GL_2$ also admits a geometric counetrpart, the category $\Sph(\Gr_{\tilde F})$ of $\GL_2(\tilde\cO)$-equivariant perverse sheaves on $\Gr_{\tilde F}$. This is a tensor category equivalent to the category of representations of $\GL_2$ (\cite{MV}). It acts on the Waldspurger category by convolutions. 
  
  Actually we work with a global version $\P^{\tilde E}(\Wald_{\pi}^x)$ of the Waldspurger category
(in geometric setting we replace $\Fq$ by an algebraically closed field $k$ of characteristic 
$p>2$). The input data for our definition of $\P^{\tilde E}(\Wald_{\pi}^x)$ is a two-sheeted covering $\pi: \tilde X\to X$ ramified at some divisor $D_{\pi}$ on $X$, a point $x\in X$, and a rank one local system $\tilde E$ on $\tilde X$. Here $\tilde X$ and $X$ are smooth projective curves over $k$ (with $X$ connected). 

 Objects of $\P^{\tilde E}(\Wald_{\pi}^x)$ are some perverse sheaves on a stack $\Wald_{\pi}^x$, which is a global model  of `the space' of $\tilde F^*$-orbits on $\Gr_{\tilde F}$. By definition, $\Wald_{\pi}^x$ classifies collections: a rank 2 vector bundle $L$ on $X$, a line bundle $\cB$ on $\pi^{-1}(X-x)$, and an isomorphism $\pi_*\cB\,\iso\, L\mid_{X-x}$.
 
 Our main result here is Theorem~\ref{Th_4} describing the action of $\Sph(\Gr_{\tilde F})$ on irreducible objects of $\P^{\tilde E}(\Wald_{\pi}^x)$. It implies the above cited multiplicity one for the Waldspurger models. This circle of ideas is very much inspired by \cite{FGV}. Note that, to the difference with the case of Whittaker categories studied in \select{loc.cit.}, the category $\P^{\tilde E}(\Wald_{\pi}^x)$ is not semi-simple.

\bigskip\medskip
\centerline{\scshape 2. Whittaker categories}

\bigskip\noindent
2.1 {\scshape Notation\ } Let $k$ denote an algebraically closed field of
characteristic $p>2$. All the schemes (or stacks) we consider are defined over $k$.
Let $X$ be a smooth projective connected curve. Fix a
prime $\ell\ne p$. For a scheme (or stack) $S$ write $\D(S)$ for the bounded
derived category of $\ell$-adic
\'etale sheaves on $S$, and $\P(S)\subset \D(S)$ for the category of perverse sheaves.

Fix a nontrivial character $\psi: \Fp\to\Qlb^*$ and denote by $\cL_{\psi}$ the
corresponding Artin-Shreier sheaf on $\A^1$. The Fourier transform functors will be always normalized to preserve perversity and purity. 

 Let $G=\GSp_4$, so $G$ is the quotient of $\Gm\times\Sp_4$ by the diagonally
embedded
$\{\pm 1\}$. Denote by $\check{G}$ the Langlands dual group to $G$ (over $\Qlb$).
We use the following notation from (\cite{Ly}, example 2 in the
appendix). The group $G$ is realized as the subgroup of $\GL(k^4)$ preserving up to a
scalar the bilinear form given by the matrix
$$
\left(
\begin{array}{cc}
0 & E_2\\
-E_2 & 0
\end{array}
\right),
$$
where $E_2$ is the unit matrix of $\GL_2$.

 Let $T$ be the  maximal torus of $G$ given by $\{(y_1,\ldots,y_4)\mid y_iy_{2+i} 
\mbox{ does not depend on  } i\}$. Let $\Lambda$ (resp., $\check{\Lambda}$) denote the
coweight (resp., weight) lattice of $T$. Write $V^{\check{\lambda}}$ for the
irreducible representation of $G$ of highest weight $\check{\lambda}$.

Let $\check{\epsilon}_i\in\check{\Lambda}$ 
be the character that sends a point of $T$ to $y_i$. 
 We have $\Lambda=\{(a_1,\ldots,a_4)\in\ZZ^4\mid a_i+a_{2+i} \mbox{ does not depend on }
i\}$ and
$$
\check{\Lambda}=\ZZ^4/\{\check{\epsilon}_1+
\check{\epsilon}_3-\check{\epsilon}_2-\check{\epsilon}_4\} 
$$

 Let $P_1\subset G$ be the parabolic subgroup preserving the isotropic subspace $ke_1$.
Let $P_2\subset G$ denote the 
Borel subgroup preserving the flag $ke_1\subset ke_1\oplus
ke_2$ of isotropic subspaces. Here $\{ e_i\}$ is the standard basis of $k^4$. 
Let $U_i$ be the unipotent radical of $P_i$ and $M_1=P_1/U_1$.

 The simple roots are $\check{\alpha}_1=\check{e}_1-\check{e}_2$ and
$\check{\alpha}_2=\check{e}_2-\check{e}_4$. The half sum of positive roots of $G$ is denoted by $\check{\rho}\in\check{\Lambda}$.

 Let $P\subset G$ denote the Siegel parabolic preserving the lagrangian
subspace $ke_1\oplus ke_2\subset k^4$. Let $U\subset P$ be its unipotent radical and
$M=P/U$. 

Set $\gamma=(1,1;0,0)\in\Lambda$, this is the dominant coweight corresponding to the
standard representation of $\check{G}\,\iso\,\GSp_4$. Fix fundamental weights
$\check{\omega}_1=(1,0,0,0)$ and $\check{\omega}_2=(1,1,0,0)$. 
So, $V^{\check{\omega}_1}$ is the standard representation. The orthogonal to the
coroot lattice is $\ZZ\check{\omega}_0$ with
$\check{\omega}_0=(1,0,1,0)$.  

\smallskip

 Note that the symplectic form $\wedge^2 V^{\check{\omega}_1}\to V^{\check{\omega}_0}$
induces an isomorphism $\det V^{\check{\omega}_1}\;\iso\; (V^{\check{\omega}_0})^{\otimes 2}$. 

\bigskip\noindent
2.2 {\scshape Hecke functor}\ \ 
Let $\Bun_G$ denote the stack of $G$-bundles on $X$. For a $G$-bundle
$\cF_G$ let
$M=V^{\check{\omega}_1}_{\cF_G}$, $\cW=V^{\check{\omega}_2}_{\cF_G}$ and $\cA=V^{\check{\omega}_0}_{\cF_G}$. In this way $\Bun_G$ becomes the stack classifying
the data: a line bundle $\cA$ on $X$, a vector bundle
$M$ of rank 4 on $X$ with a symplectic form $\wedge^2 M\to\cA$. The exact sequence
$$
0\to\cW\to \wedge^2 M\to\cA\to 0
$$
splits canonically.
 
 Denote by $\cH_G$ the stack of collections: $x\in X, \cF_G,\cF'_G\in\Bun_G$ and
$\cF_G\,\iso\,\cF'_G\mid_{X-x}$ such that $\cF_G$ is in the position $\gamma$ with respect
to $\cF'_G$. In other words, we have $\cA'=\cA(x)$, $\, M\subset M'$, the diagrams
commute
$$
\begin{array}{ccc}
\wedge^2 M' & \to & \cA'\\ 
\uparrow && \uparrow\\
\wedge^2 M & \to & \cA
\end{array}
$$
and
$$
\begin{array}{ccc}
\det M' & \iso & {\cA'}^2\\
\uparrow && \uparrow\\
\det M & \iso & \cA^2
\end{array}
$$
and $M/M'(-x)\subset M'/M'(-x)$ is a lagrangian subspace. 

 We have a diagram $\Bun_G\getsup{\gp}\cH_G\toup{\gq}\Bun_G$, where the map
$\gp$ (resp., $\gq$) sends the above collection to $\cF_G$ (resp., $\cF'_G$).
Let $\supp: \cH_G\to X$ be the map sending the above point to $x$.
Note that $\gq$ is smooth of relative dimension $1+\<\gamma, 2\check{\rho}\>$.
Let 
$$
\H:\D(\Bun_G)\to\D(X\times\Bun_G)
$$ 
denote the Hecke functor corresponding to $\gamma$, that is, 
$$
\H(K)=(\supp\times\gp)_!\gq^*K\otimes \Qlb(\frac{1}{2})[1]^{\otimes 1+\<\gamma, 
2\check{\rho}\>}
$$

\bigskip\noindent
2.3 {\scshape Drinfeld compactifications}\ \ 
 We fix a particular $T$-torsor on $X$ with trivial conductor $(\cF_T,\tilde\omega)$ 
by requiring $\cL^{\check{\omega}_1}_{\cF_T}\iso\Omega$. The pair
$(\cF_T,\tilde\omega)$ with this property is defined up to a unique isomorphism, and
we have $\cL^{\check{\omega}_2}_{\cF_T}\;\iso\;\Omega$ and
$\cL^{\check{\omega}_0}_{\cF_T}\;\iso\;\Omega^{-1}$.

 For $k=1,2,3$ define the stack $\bar\cQ_k$ as follows. It classifies a point
$\cF_G\in\Bun_G$ together with sections $t_1,\ldots, t_k$ satisfying Plucker
relations, where
$$
\begin{array}{l}
t_1: \Omega\hook{} M\\
t_2: \Omega\hook{} \cW\\
t_3: \Omega^{-1}\hook{} \cA
\end{array}
$$
 It is understood that Plucker relations are empty for $k=1$, and for $k=2,3$ they mean
that, at the generic point of $X$, the sections $t_1,\ldots,t_k$ come from a
$B$-structure on $\cF_G$. 
 
 Set $\bar\cQ=\bar\cQ_3$. Let also $\bar\cQ_{k,ex}$ be the stack defined in the same way as $\bar\cQ_k$ with the only difference that the last section $t_k$ is not necessairy an inclusion
(here `ex' stands for `extended'). So, $\bar\cQ_k\subset\bar\cQ_{k,ex}$ is an open substack.

 Denote by $\pi_{k+1,k}: \bar\cQ_{k+1}\to\bar\cQ_k$ and
$\pi_{k+1,k,ex}:\bar\cQ_{k+1,ex}\to \bar\cQ_k$ the natural forgetful maps.

 For each $k$ we have the diagram 
$$
\bar\cQ_{k,ex}\;\getsup{\gp_{k,ex}}\;
\bar\cQ_{k,ex}\times_{\Bun_G}\cH_G\;\toup{\gq_{k,ex}}\; \bar\cQ_{k,ex},
$$
where we used the map $\gp:\cH_G\to\Bun_G$ in the definition of the fibred product,
$\gp_{k,ex}$ is the projection, and $\gq_{k,ex}$ sends a point of
$\bar\cQ_{k,ex}\times_{\Bun_G}\cH_G$ to $(\cF'_G, t'_1, \ldots, t'_k)$ with $t'_i$
being the compositions
$$
\begin{array}{l}
t_1: \Omega\to M\hook{} M'\\
t_2: \Omega\to \cW\hook{}\cW'\\
t_3: \Omega^{-1}\to \cA\hook{}\cA'
\end{array}
$$
For $k=1,2,3$ we have the functor $\H^{\bar\cQ_{k,ex}}:\D(\bar\cQ_{k,ex})\to
\D(X\times\bar\cQ_{k,ex})$ given by
$$
\H^{\bar\cQ_{k,ex}}(K)=(\supp\times\gp_{k,ex})_!\gq_{k,ex}^*K\otimes
\Qlb(\frac{1}{2})[1]^{\otimes \<\gamma, 2\check{\rho}\>}
$$

 The restriction of $\gq_{k,ex}$ to $\bar\cQ_k\times_{\Bun_G}\cH_G$ factors through
$\bar\cQ_k\subset \bar\cQ_{k,ex}$. So, we also have diagrams
$$
\bar\cQ_k\;\getsup{\gp_k}\;
\bar\cQ_k\times_{\Bun_G}\cH_G\;\toup{\gq_k}\; \bar\cQ_k,
$$
where $\gp_k$ (resp., $\gq_k$) is the restriction of $\gp_{k,ex}$ (resp., of $\gq_{k,
ex}$). For $k=1,2,3$ denote by 
$$
\H^{\bar\cQ_k}:\D(\bar\cQ_k)\to
\D(X\times\bar\cQ_k)
$$ 
the functor given by
\begin{equation}
\label{def_Hecke_functor}
\H^{\bar\cQ_k}(K)=(\supp\times\gp_k)_!\gq_k^*K\otimes
\Qlb(\frac{1}{2})[1]^{\otimes \<\gamma, 2\check{\rho}\>}
\end{equation}

 The projection $\alpha:\bar\cQ_1\to\Bun_G$ fits into the diagram 
$$
\begin{array}{ccccc}
\bar\cQ_1 & \getsup{\gp_1} &
\bar\cQ_1\times_{\Bun_G}\cH_G & \toup{\gq_1} & \bar\cQ_1\\
\downarrow\lefteqn{\scriptstyle \alpha} && \downarrow && 
\downarrow\lefteqn{\scriptstyle \alpha}\\
\Bun_G & \getsup{\gp} & \cH_G & \toup{\gq} & \Bun_G,
\end{array}
$$
in which the left square is cartesian. So, $(\id\times\alpha)^*\comp \H\;\iso\;
\H^{\bar\cQ_1}\comp\alpha^*[1](\frac{1}{2})$ naturally. Over the open substack of $\Bun_G$
given by
$\Ext^1(\Omega, M)=0$, the map
$\alpha:\bar\cQ_1\to\Bun_G$ is smooth.

\bigskip\noindent
2.4 Let $\pi_{0,1,ex}:\bar\cQ_{0,ex}\to\bar\cQ_1$ be the vector bundle with fibre
consisting of all sections $t_0:\Omega\to\cA$. Let $i_0:\bar\cQ_1\to\bar\cQ_{0,ex}$ denote
the zero section and $j:\bar\cQ_0\subset\bar\cQ_{0,ex}$ its complement given by: $t_0$ is
an inclusion.

 We have the diagram 
$$
\bar\cQ_{0,ex}\;\toup{\gp_{0,ex}}\;
\bar\cQ_{0,ex}\times_{\Bun_G}\cH_G\;\toup{\gq_{0,ex}}\;\bar\cQ_{0,ex},
$$
where we used $\gp:\cH_G\to\Bun_G$ in the definition of the fibred product,
$\gp_{0,ex}$ is the projection, and $\gq_{0,ex}$ sends a point of
$\bar\cQ_{0,ex}\times_{\Bun_G}\cH_G$ to $(\cF'_G, t'_0,t'_1)$. Here, as above, $t'_i$
are the compositions 
$$
\begin{array}{l}
t_0:\Omega\to\cA\hook{}\cA'\\
t_1: \Omega\hook{} M\hook{} M'
\end{array}
$$  

 Restricting, one gets the diagram 
$\bar\cQ_0\;\toup{\gp_0}\;
\bar\cQ_0\times_{\Bun_G}\cH_G\;\toup{\gq_0}\;\bar\cQ_0,$.
The functors 
$$
\H^{\bar\cQ_{0,ex}}:\D(\bar\cQ_{0,ex})\to\D(X\times \bar\cQ_{0,ex})
$$ 
and
$\H^{\bar\cQ_0}:\D(\bar\cQ_0)\to\D(X\times \bar\cQ_0)$ are defined as in
(\ref{def_Hecke_functor}).

\begin{Rem} For any $K\in\D(\bar\cQ_{0,ex})$ we have a natural isomorphism of distinguished
triangles
$$
\begin{array}{ccccc}
j_!\H^{\bar\cQ_0}(j^*K) & \to & \H^{\bar\cQ_{0,ex}}(K) & \to &
(i_0)_*\H^{\bar\cQ_1}(i_0^*K)\\
\downarrow && \downarrow\lefteqn{\scriptstyle\id} && \downarrow \\
j_!j^*\,\H^{\bar\cQ_{0,ex}}(K) & \to & \H^{\bar\cQ_{0,ex}}(K) & \to &
(i_0)_*i^*_0\,\H^{\bar\cQ_{0,ex}}(K)
\end{array}
$$
\end{Rem}

\bigskip\noindent
2.5 {\scshape Categories to construct}\ 
We will introduce triangulated categories $\D^W(\bar\cQ_k)$ (resp.,
$\D^W(\bar\cQ_{k,ex})$) of sheaves on
$\bar\cQ_k$ (resp., on $\bar\cQ_{k,ex}
$) for $k=0,1,2,3$ (resp, for $k=0,2,3$). 

 Each $\D^W(\bar\cQ_k)$ will be a
full triangulated subcategory of $\D(\bar\cQ_k)$ defined by the condition that $K\in
\D^W(\bar\cQ_k)$ if its perverse cohomology belong to a certain Serre subcategory
$\P^W(\bar\cQ_k)$ singled out by some equivariance condition; and similarly for
$\D^W(\bar\cQ_{k,ex})$.

 Though we don't reflect this in the notation, all our equivariant categories (except 
$\D^W(\bar\cQ_1)$) will depend on the character $\psi$.

\bigskip\noindent
2.6 Let $y\in X$ be a closed point. 
For $k=1,2$ let $\bar\cQ_k^y\subset\bar\cQ_k$ be the open substack given by the
condition that neither of the maps $t_1,\ldots,t_k$ has zero at $y$.

 If $(\cF_G, t_1,\ldots,t_k)$ is a point of $\bar\cQ^y_k$ then over the formal disk
$D_y$ at $y$ we obtain a $P_k$-torsor $\cF_{P_k}$. Let $N_{k,y}\to \bar\cQ_k^y$ be
stack whose fibre over a point of $\bar\cQ_k^y$ is
$$ 
\H^0(D_y^*, \;\cF_{P_k}\times_{P_k} U_k)
$$
This is an ind-groupscheme over $\bar\cQ_k^y$, it can be represented as a union of
group schemes $^iN_{k,y}$ for $i\in\NN$, where $^iN_{k,y}\hook{} {^{i+1}N_{k,y}}$ is a
closed immersion, and $^iN_{k,y}/{^0N_{k,y}}$ is of finite type over $\bar\cQ^y_k$ for
$i>0$. We assume that the fibre of $^0N_{k,y}\to \bar\cQ_k^y$ is
$$
\H^0(D_y, \;\cF_{P_k}\times_{P_k} U_k)
$$

 Let $\cH_{k,y}\to\bar\cQ^y_k$ denote the stack over $\bar\cQ^y_k$ with fibre
$N_{k,y}/{^0N_{k,y}}$. This is an ind-scheme over $\bar\cQ^y_k$, and we have 
$$
\cH_{k,y}=\mathop{\cup}\limits_i {^i\cH_{k,y}},
$$ 
where $^i\pr_k: {^i\cH_{k,y}}\to \bar\cQ^y_k$ is the stack with fibre
$^iN_{k,y}/{^0N_{k,y}}$.

\bigskip\noindent
2.7 {\scshape Groupoids} \  \ 
As in (\cite{G}, sect. 4.3) one endows $\cH_{k,y}$ with the
structure of a groupoid over $\bar\cQ^y_k$. We denote by 
$$
^i\act_k: {^i\cH_{k,y}}\to
\bar\cQ^y_k
$$ 
the restriction of the action map.

 For $k=1,2$ define the open substack $\bar\cQ^y_{k+1,ex}\subset \bar\cQ_{k+1,ex}$ 
as $\bar\cQ^y_{k+1,ex}=\bar\cQ^y_k\times_{\bar\cQ_k}\bar\cQ_{k+1,ex}$.
The groupoid $\cH_{k,y}\to \bar\cQ^y_k\;$ "lifts" to $\bar\cQ^y_{k+1,ex}$. In other
words, 
$$
\cH_{k,y}\times_{\bar\cQ^y_k} \;\bar\cQ^y_{k+1,ex}
$$ 
has a structure of a groupoid
over $\bar\cQ^y_{k+1,ex}$ (we used the projections to define the above fibre product).
Moreover, the diagram is cartesian
$$
\begin{array}{ccc}
\cH_{k,y}\times_{\bar\cQ^y_k} \bar\cQ^y_{k+1,ex} & \;\,\toup{\act} &
\bar\cQ^y_{k+1,ex}\\
\downarrow\lefteqn{\scriptstyle \id\times\pi_{k+1,k,ex}} &&
\downarrow\lefteqn{\scriptstyle
\pi_{k+1,k,ex}}\\
\cH_{k,y} & \;\,\toup{\act} & \bar\cQ^y_k
\end{array}
$$

 Denote by 
$$
^i\act_{k,ex}: {^i\cH_{k,y}}\times_{\bar\cQ^y_k}
\bar\cQ^y_{k+1,ex}\to
\bar\cQ^y_{k+1,ex}
$$ 
the action map.

 Let $\bar\cQ^y_{0,ex}\subset \bar\cQ_{0,ex}$ be the preimage of $\bar\cQ^y_1$ under
$\pi_{0,1,ex}: \bar\cQ_{0,ex}\to\bar\cQ_1$. The groupoid $\cH_{1,y}\to \bar\cQ^y_1\;$
"lifts" to $\bar\cQ^y_{0,ex}$ in the same sense as above.   

\bigskip
\noindent
2.8 \  We single out the subgroupoid $\cH_{0,y}\subset \cH_{1,y}$ as follows. 

Let $U_0\subset U_1$ denote the center of $U_1$. The exact
sequence $1\to U_0\to U_1\to U_1/U_0\to 1$ does not split, we have $U_0\;\iso\; \Ga$
and $U_1/U_0\;\iso\; \Ga^2$.

 The stack $\Bun_{P_1}$ classifies: a $G$-torsor $\cF_G$ together with a line subbundle 
$L_1\subset M$. Note that $L_1$ is automatically isotropic and denote by
$L_{-1}\subset M$ its orthogonal complement. For such $\cF_{P_1}\in\Bun_{P_1}$ the
vector bundle
$\cF_{P_1}\times_{P_1} U_0$ is $L_1^2\otimes\cA^{-1}$. It is
understood that $P_1$ acts on $U_0$ adjointly.  

 By definition, the fibre of $\cH_{0,y}\to \bar\cQ^y_1$ is
$$
\H^0(D^*_y, \;\cF_{P_1}\times_{P_1} U_0)/\H^0(D_y, \;\cF_{P_1}\times_{P_1} U_0)\;\iso\; 
\Omega^2\otimes\cA^{-1}(\infty y)/ \Omega^2\otimes\cA^{-1}
$$

 Denote by $^i\cH_{0,y}\subset \cH_{0,y}$ the subgroupoid with fibre
$$
\Omega^2\otimes\cA^{-1}(i y)/ \Omega^2\otimes\cA^{-1}
$$
We write $^i\pr_0: {^i\cH_{0,y}}\to \bar\cQ^y_1$ for the projection and 
$^i\act_0: {^i\cH_{0,y}}\to \bar\cQ^y_1$ for the action map. 

Let also
$$
^i\act_{0,ex}: {^i\cH_{0,y}}\times_{\bar\cQ^y_1} \bar\cQ^y_{0,ex}\to  \bar\cQ^y_{0,ex}
$$
denote the action map.

\bigskip\noindent
2.9 {\scshape Characters}\ \ 
Let us construct a natural map
$$
\chi_{0,y}: \cH_{0,y}\times_{\bar\cQ^y_1} \bar\cQ^y_{0,ex}\to\A^1
$$
The element $t_0:\Omega\to\cA$ gives rise to a morphism
$$
\Omega^2\otimes\cA^{-1}(i y)/ \Omega^2\otimes\cA^{-1}\to \Omega(\infty y)/\Omega
$$
and we take the residue of the image of $g\in \cH_{0,y}$ under this map.

Let us construct for $k=1,2$ a natural map
$$
\chi_{k,y}: \cH_{k,y}\times_{\bar\cQ^y_k} \;\bar\cQ^y_{k+1,ex}\to \A^1
$$

\noindent
{\bf CASE $k=1$.} \\  
If $\cF_{P_1}$ is a $P_1$-torsor on a scheme given by $(L_1\subset L_{-1}\subset M)$
then the vector bundle $\cF_{P_1}\times_{P_1} U_1/U_0$ is $\Hom(L_{-1}/L_1,\; L_1)$. 

 Recall that a point of $\bar\cQ^y_1$ defines a $P_1$-torsor $\cF_{P_1}$ on $D_y$. Let
$\cE_1\to\bar\cQ^y_1$ be the stack whose fibre over a point of $\bar\cQ^y_1$ is
$$
\H^0(D_y^*,\; \cF_{P_1}\times_{P_1} U_1/U_0)/\H^0(D_y,\; \cF_{P_1}\times_{P_1} U_1/U_0)
\;\iso\; (L_{-1}/\Omega)^*\otimes (\Omega(\infty y)/\Omega)
$$

 We have a natural map $\cH_{1,y}\to \cE_1$ over $\bar\cQ^y_1$.  Given a point of
$\bar\cQ^y_{2,ex}$, over $D_y$ the section $t_2:\Omega\to\cW$ gives rise to a map
$s:\cO\to L_{-1}/\Omega$ such that $t_2=t_1\wedge s$. By definition, $\chi_{1,y}$ is
the residue of the pairing of $s$ with the image of $g\in\cH_{1,y}$ in $\cE_1$. 

\medskip\noindent
{\bf CASE $k=2$.}\\
Given a point of $\bar\cQ^y_2$ we obtain a $P_2$-torsor $\cF_{P_2}$ over $D_y$. Let
$U_{2,ab}$ be the abelinization of $U_2$ then 
\begin{equation}
\label{space_1}
\H^0(D_y^*,\; \cF_{P_2}\times_{P_2} U_{2,ab})/
\H^0(D_y,\; \cF_{P_2}\times_{P_2} U_{2,ab})\;\iso\;
\Omega(\infty y)/\Omega \oplus \cA^{-1}(\infty y)/\cA^{-1},
\end{equation}
where the two summands correspond to the simple roots of $G$. To define
$\chi_{2,y}$, we take the image of
$g\in\cH_{2,y}$ in (\ref{space_1}), pair it with $t_3:\Omega^{-1}\to\cA$ and take the
sum of residues. 
 
\medskip

 For $k=1,2$ write 
$$
^i\chi_{k,y}: {^i\cH_{k,y}}\times_{\bar\cQ^y_k}
\;\bar\cQ^y_{k+1,ex}\to
\A^1
$$ 
for the restriction of $\chi_{k,y}$, and similarly for $^i\chi_{0,y}$.

\bigskip\smallskip
\noindent
2.10 {\scshape Categories on $\bar\cQ^y_{k+1,ex}$}

\medskip\noindent
 For $k=1,2$ define the full subcategory $\P^W(\bar\cQ^y_{k+1,ex})\subset
\P(\bar\cQ^y_{k+1,ex})$ to consist of all perverse sheaves $K\in
\P(\bar\cQ^y_{k+1,ex})$ with the property:

 For any $i\in\NN$ there is an isomorphism on ${^i\cH_{k,y}}\times_{\bar\cQ^y_k}
\;\bar\cQ^y_{k+1,ex}$
$$
^i\chi_{k,y}^*(\cL_{\psi})\otimes \pr_2^*K\;\iso \;
{^i\act_{k,ex}^*}K
$$
whose restriction to the unit section $\bar\cQ^y_{k+1,ex}\subset {^i\cH_{k,y}}\times_{\bar\cQ^y_k}
\;\bar\cQ^y_{k+1,ex}$ is the identity map.
 
 Similarly, $\P^W(\bar\cQ^y_{0,ex})\subset \P(\bar\cQ^y_{0,ex})$ is the full
subcategory consisiting of perverse sheaves $K$ with the property:
 
 For any $i\in\NN$ there is an isomorphism on $^i\cH_{0,y}\times_{\bar\cQ^y_1}
\bar\cQ^y_{0,ex}$
$$
^i\chi_{0,y}^*(\cL_{\psi})\otimes \pr_2^*K\;\iso \;
{^i\act_{0,ex}^*}K
$$
whose restriction to the unit section $\bar\cQ^y_{0,ex}\subset
{^i\cH_{0,y}}\times_{\bar\cQ^y_1}
\;\bar\cQ^y_{0,ex}$ is the identity map.

 For $k=-1,1,2$ as in (\cite{G}, Sect.~4.7-4.8) one shows that
$\P^W(\bar\cQ^y_{k+1,ex})$ is a Serre subcategory of $\P(\bar\cQ^y_{k+1,ex})$. Then
$\D^W(\bar\cQ^y_{k+1,ex})\subset
\D(\bar\cQ^y_{k+1,ex})$ is a full triangulated subcategory consisting of objects whose
perverse cohomology belong to $\P^W(\bar\cQ^y_{k+1,ex})$.

\bigskip\noindent
2.11 In all the three cases $k=-1,1,2$ define $\P^W(\bar\cQ_{k+1,ex})\subset
\P(\bar\cQ_{k+1,ex})$ as the full subcategory consisting of $K\in \P(\bar\cQ_{k+1,ex})$
such that
$$
K\mid_{\bar\cQ^y_{k+1,ex}}\in \P^W(\bar\cQ^y_{k+1,ex})
$$ 
for any $y\in X$.
Then $\P^W(\bar\cQ_{k+1,ex})$ is a Serre subcategory of $\P(\bar\cQ_{k+1,ex})$. Set
$\D^W(\bar\cQ_{k+1,ex})$ to be the full triangulated subcategory of
$\D(\bar\cQ_{k+1,ex})$ generated by $\P^W(\bar\cQ_{k+1,ex})$. 

\smallskip

 For $k=-1,1,2$ we also have the categories $\D^W(\bar\cQ_{k+1})$
and $\P^W(\bar\cQ_{k+1})$ defined in a similar fashion, because the open substack
$\bar\cQ_{k+1}\subset \bar\cQ_{k+1,ex}$ is preserved by the action of the corresponding
groupoid.

\bigskip\noindent
2.12 Recall the vector bundle $\pi_{0,1,ex}:\bar\cQ_{0,ex}\to\bar\cQ_1$.
Let
$i_0:\bar\cQ_1\hook{} \bar\cQ_{0,ex}$ denote its zero section and 
$j:\bar\cQ_0\hook{} \bar\cQ_{0,ex}$ the complement to the zero section.

 Let $\D^W(\bar\cQ_1)\subset\D(\bar\cQ_1)$ be the full triangulated subcategory
consisiting of those $K\in \D(\bar\cQ_1)$ for which $(i_0)_!K\in \D^W(\bar\cQ_{0,ex})$.
The Serre subcategory $\P^W(\bar\cQ_1)\subset\P(\bar\cQ_1)$ is defined by the same
condition. 

 In other words, $K\in\P(\bar\cQ_1)$ lies in $\P^W(\bar\cQ_1)$ if and only if it is invariant
under the action of the groupoids $\cH_{0,y}$ for all $y\in X$.

\smallskip

 Note that any $K\in\D^W(\bar\cQ_{0,ex})$ fits into a distinguished triangle
$j_!j^*K\to K\to (i_0)_!i_0^*K$ with $i_0^*K\in \D^W(\bar\cQ_1)$ and
$j^*K\in\D^W(\bar\cQ_0)$.

\bigskip\noindent
2.13 {\scshape Stratifications}

\smallskip\noindent
 For $k=1,2$ stratify $\bar\cQ_k$ as follows. For a string of nonnegative integers
$\bar d=(d_1,\ldots,d_k)$ let $^{\bar d}\bar\cQ_k\subset\bar\cQ_k$ be the locally
closed substack given by: there exist $D_1\in X^{(d_1)},\ldots,D_k\in X^{(d_k)}$ such
that 
$$
t_i: \cL^{\check{\omega}_i}_{\cF_T}(D_i)\hook{} V^{\check{\omega}_i}_{\cF_G}
$$
is a subbundle for $i=1,\ldots,k$.
In other words, $\Omega(D_1)\subset M$ is a subbundle for $k=1$; and for $k=2$
there is one more condition: $\Omega(D_2)\subset\cW$ is a subbundle.
 
 Let $^{\bar d}\bar\cQ_k^y$ be the preimage of $\bar\cQ^y_k$ under  $^{\bar d}\bar\cQ_k\hook{} \bar\cQ_k$. The stack $^{\bar d}\bar\cQ_k^y$ is stable under the action of $\cH_{k,y}$ on $\bar\cQ_k^y$.
 For $k=1,2$ set 
$$
^{\bar d}\bar\cQ_{k+1,ex}={^{\bar
d}\bar\cQ_k}\times_{\bar\cQ_k}\bar\cQ_{k+1,ex}
$$ 
Set also 
$$
^{\bar d}\bar\cQ_{0,ex}={^{\bar
d}\bar\cQ_1}\times_{\bar\cQ_1}\bar\cQ_{0,ex}
$$ 
For $y\in X$ denote by
$^{\bar d}\bar\cQ_{k+1,ex}^y$ the preimage of $\bar\cQ^y_{k+1,ex}$ under 
$^{\bar d}\bar\cQ_{k+1,ex}\to \bar\cQ_{k+1,ex}$. For $k=1,2$ the stack $^{\bar
d}\bar\cQ_{k+1,ex}^y$ is stable under the action of $\cH_{k,y}$ on
$\bar\cQ_{k+1,ex}^y$.
 The stack $^{\bar d}\bar\cQ^y_{0,ex}$ is stable under the action of $\cH_{0,y}$ on
$\bar\cQ^y_{0,ex}$.

 Thus, following the same lines one defines the categories $\P^W(^{\bar
d}\bar\cQ_{k+1,ex}^y)$ and $\D^W(^{\bar
d}\bar\cQ_{k+1,ex}^y)$, and further $\P^W(^{\bar
d}\bar\cQ_{k+1,ex})$ and $\D^W(^{\bar
d}\bar\cQ_{k+1,ex})$ for $k=-1,1,2$. Similarly for $\P^W(^{\bar d}\bar\cQ_2)$ and
$\D^W(^{\bar d}\bar\cQ_2)$.

 By abuse of notation, write 
$$
i_0: {^{\bar d}\bar\cQ_1}\hook{} {^{\bar d}\bar\cQ_{0,ex}}
$$
for the natural closed immersion. Denote by 
$
\D^W(^{\bar d}\bar\cQ_1)\subset \D(^{\bar
d}\bar\cQ_1)
$ 
the full triangulated subcategory consisting of $K\in \D({^{\bar
d}\bar\cQ_1})$ such that 
$$
(i_0)_!K\in \D^W({^{\bar d}\bar\cQ_{0,ex}})
$$

As in (\cite{G},
Lemma~4.11) one shows the following (cf. also Appendix A).

\begin{Lm} 
\label{Lm_on_stratifications}
1) Let $k=-1,1,2$. The functors of $*$- and $!$-restriction map $\D^W(\bar\cQ_{k+1,ex})$ to $\D^W({^{\bar d}\bar\cQ_{k+1,ex}})$. The functors of $*$- and $!$-direct image map $\D^W({^{\bar d}\bar\cQ_{k+1,ex}})$ to $\D^W(\bar\cQ_{k+1,ex})$. 

\smallskip

 For $K\in \D(\bar\cQ_{k+1,ex})$ we have $K\in\D^W(\bar\cQ_{k+1,ex})$ if and only if
its $*$-restriction (or, equivalently, $!$-restriction) to $^{\bar d}\bar\cQ_{k+1,ex}$ lies in 
$\D^W({^{\bar d}\bar\cQ_{k+1,ex}})$ for any $\bar d$. 

\smallskip\noindent
2) Let $k=1,2$. For $K\in \D(\bar\cQ_k)$ we have $K\in\D^W(\bar\cQ_k)$ if and only if its $*$-restriction (or, equivalently, $!$-restriction) to $^{\bar d}\bar\cQ_k$ lies in $\D^W(^{\bar d}\bar\cQ_k)$ for any
$\bar d$. 
\end{Lm}

\bigskip\noindent
2.14 For $k=1,2$ define a closed substack $^{\bar d}\bar\cQ'_{k+1,ex}\hook{} {^{\bar
d}\bar\cQ_{k+1,ex}}$ by the conditions: $t_2$ comes from $\H^0(X,
\Omega^{-1}\otimes\cW(-2D_1))$ in both cases,  
and for $k=2$ we require in addition that
$t_3$ comes from $\H^0(X, \cA\otimes\Omega(-2D'_2))$, where we have put $D'_2=D_2-D_1$.

 Let us define for $k=1,2$ a natural map 
$$
^{\bar d}\chi_{k+1,ex}: {^{\bar d}\bar\cQ'_{k+1,ex}}\to\A^1
$$

\noindent
{\bf CASE $k=1$.} The stack $^{\bar d}\bar\cQ'_{2,ex}$ classifies collections: $D_1\in
X^{(d_1)}$,  a $P_1$-torsor on $X$ given by 
$$
(L_1\subset L_{-1}\subset M)
$$ 
with $L_1\;\iso\; \Omega(D_1)$, and a section $s: \cO(D_1)\to L_{-1}/L_1$. 

 The map
$^{\bar d}\chi_{2,ex}$ sends this collection to the class in $\Ext^1(
\cO(D_1), \Omega(D_1))\,\iso\, k$ of the pull-back of $0\to L_1\to L_2\to L_2/L_1\to 0$
under $s$. 
  
\medskip
\noindent
{\bf CASE $k=2$.} Note that $\Bun_P$ is the stack classifying: a rank 2 vector
bundle $L_2$ on $X$, a line bundle $\cA$ on $X$, and an exact sequence $0\to \Sym^2
L_2\to ?\to\cA\to 0$.   For such $\cF_P\in\Bun_P$ the vector bundle $\cF_P\times_P U$
is $(\Sym^2 L_2)\otimes\cA^{-1}$. 

\smallskip

 Therefore, the stack $^{\bar d}\bar\cQ'_{3,ex}$ classifies collections: $D_1\in
X^{(d_1)}$,  
$D'_2\in X^{(d_2-d_1)}$ with $D'_2\ge D_1$, two exact sequences
$$
0\to L_1\to L_2\to L_2/L_1\to 0
$$ 
and
\begin{equation}
\label{sequence_2}
0\to \Sym^2 L_2\to ?\to\cA\to 0
\end{equation}
with $L_1\;\iso\;\Omega(D_1)$ and $L_2/L_1\;\iso\;\cO(D'_2)$, and a
section $t_3: \Omega^{-1}(2D'_2)\to \cA$. 

 The map $^{\bar d}\chi_{3,ex}$ sends this collection to the sum of two numbers, the first
being defined as for $^{\bar d}\chi_{2,ex}$, and the second is the class in
$\Ext^1(\Omega^{-1}(2D'_2), \cO(2D'_2))\;\iso\; k$ of the pull-back of
\begin{equation}
\label{sequence_1}
0\to \Sym^2(L_2/L_1)\to ?\to\cA\to 0
\end{equation}
under $t_3$. Here (\ref{sequence_1}) is the push-forward of (\ref{sequence_2}) under
$\Sym^2 L_2\to
\Sym^2(L_2/L_1)$. 

\bigskip\noindent
2.15 For $k=1,2$ define the stack $^{\bar d}\cP_{k+1,ex}$ as follows.

 The stack $^{\bar d}\cP_{2,ex}$ classifies: $D_1\in X^{(d_1)}$, 
a rank 2 vector bundle $M_2$ on $X$ with section $s: \cO(D_1)\to M_2$. 

 The stack $^{\bar d}\cP_{3,ex}$ classifies: $D_1\in
X^{(d_1)}$, $D'_2\in X^{(d_2-d_1)}$ with $D'_2\ge D_1$, a line bundle $\cA$ on $X$ with
a section $\Omega^{-1}(2D'_2)\to\cA$. 

 In both cases we have a projection $\phi_{k+1,ex}: {^{\bar d}\bar\cQ'_{k+1,ex}}\to
{^{\bar d}\cP_{k+1,ex}}$. For $k=1$ it is given by $M_2=L_{-1}/L_1$.  

\bigskip\noindent  
2.16  For $\bar d=(d_1,d_2)$ let $^{\bar d}\bar\cQ'_2\subset {^{\bar d}\bar\cQ_2}$ be the closed substack given by $D'_2\ge D_1$, where $D'_2=D_2-D_1$. We have a natural map
$$
^{\bar d}\chi_2: {^{\bar d}\bar\cQ'_2}\to\A^1
$$
defined in the same way as $^{\bar d}\chi_{2,ex}$.

For $k=1,2$ and $\bar d=(d_1,\ldots,d_k)$ as above define the stack
$^{\bar d}\cP_k$ as follows.
The stack $^{\bar d}\cP_1$ 
classifies $D_1\in X^{(d_1)}$
and an exact sequence of vector bundles on $X$ 
$$
0\to L_1\to L_{-1}\to L_{-1}/L_1\to 0
$$ 
with $L_1\iso\Omega(D_1)$, where $L_{-1}/L_1$ is of rank 2. 

 The stack $^{\bar d}\cP_2$ 
classifies $D_1\in X^{(d_1)}$, $D'_2\in X^{(d_2-d_1)}$ with $D'_2\ge D_1$, a line bundle $\cA$ on $X$, and an exact sequence on $X$
$$
0 \to \cO(D'_2)\to M_2\to \cA(-D'_2)\to 0
$$

We have projections
$\phi_1: {^{\bar d}\bar\cQ_1}\to {^{\bar d}\cP_1}$ and $\phi_2: {^{\bar d}\bar\cQ'_2}\to {^{\bar d}\cP_2}$.

 As in (\cite{G}, Proposition 4.13) one proves

\begin{Lm} 
\label{Lm_2}
For $k=1,2$ and a string of nonnegative integers $\bar d=(d_1,\ldots,d_k)$ we have the following. \\
i) Any object $K\in\D^W({^{\bar
d}\bar\cQ_{k+1,ex}})$ is supported at ${^{\bar
d}\bar\cQ'_{k+1,ex}}$. The functor $K\mapsto {^{\bar d}\chi_{k+1,ex}^*}\cL_{\psi}\otimes
\phi_{k+1,ex}^*K$ provides an equivalence of categories
$\D(^{\bar d}\cP_{k+1,ex})\to \D^W({^{\bar
d}\bar\cQ_{k+1,ex}})$. 

\smallskip\noindent
ii) We have an equivalence of categories $\D(^{\bar d}\cP_k)\to \D^W({^{\bar d}\bar\cQ_k})$. 
For $k=1$ it is given by the functor $K\mapsto\phi_1^*K$,
whence for $k=2$ it is given by the functor $K\mapsto \phi_2^*K\otimes 
{^{\bar d}\chi_2^*}\cL_{\psi}$.  \QED
\end{Lm}

\bigskip\bigskip
\centerline{\scshape 3. Whittaker functors}

\bigskip\noindent
In this section we prove the following theorem.

\begin{Th}
\label{Th_1}
i) There is an equivalence of categories $W_{1,0,ex}:\D(\bar\cQ_1)\to
\D^W(\bar\cQ_{0,ex})$, which is $t$-exact,
and
$(\pi_{0,1,ex})_!$ is quasi-inverse to it. Moreover, for any $K\in
\D^W(\bar\cQ_{0,ex})$ the natural map $(\pi_{0,1,ex})_!K\to (\pi_{0,1,ex})_* K$
is an isomorphism.

\medskip\noindent
ii) For $k=1,2$ there is an equivalence of categories $W_{k,k+1,ex}:\D^W(\bar\cQ_k)\to
\D^W(\bar\cQ_{k+1,ex})$, which is $t$-exact,
and
$(\pi_{k+1,k,ex})_!$ is quasi-inverse to it. Moreover, for any $K\in
\D^W(\bar\cQ_{k+1,ex})$ the natural map $(\pi_{k+1,k,ex})_!K\to (\pi_{k+1,k,ex})_* K$
is an isomorphism.
\end{Th}

\medskip\noindent
3.1 First, we explain what the corresponding functors do on strata. For $k=1,2$ let $\bar d=(d_1,\ldots,d_k)$ be a string of nonnegative integers. 
Using Lemma~\ref{Lm_2}, define 
$$
^{\bar d}W_{k,k+1,ex}: \D^W(^{\bar d}\bar\cQ_k)\to \D^W(^{\bar d}\bar\cQ_{k+1,ex})
$$ 
as the composition
$$
\D^W(^{\bar d}\bar\cQ_k)\;\iso\;  \D(^{\bar d}\cP_k)\;
\toup{\Four}\; \D(^{\bar d}\cP_{k+1,ex})\;\iso\;
\D^W(^{\bar d}\bar\cQ_{k+1,ex})
$$
 
So, $^{\bar d}W_{k,k+1,ex}$ is an equivalence of triangulated categories and $t$-exact. It also follows from the standard properties of the Fourier transform that 
$(\pi_{k+1,k,ex})_!$ is quasi-inverse to $^{\bar d}W_{k,k+1,ex}$, and we have 
$(\pi_{k+1,k,ex})_!K\;\iso\; (\pi_{k+1,k,ex})_*K$ for any $K\in \D^W(^{\bar d}\bar\cQ_{k+1,ex})$.

\medskip\noindent
3.2  For $\bar d=(d_1)$ define the functor 
$$
^{\bar d}W_{1,0,ex}: \D(^{\bar d}\bar\cQ_1)\to \D^W(^{\bar d}\bar\cQ_{0,ex})
$$
 as follows.  Let $^{\bar d}\cE\to {^{\bar d}\bar\cQ_1}$ be the stack whose fibre over a point of $^{\bar d}\bar\cQ_1$ is the stack of exact sequences $0\to \Omega^2(2D_1)\otimes\cA^{-1}\to ?\to\cO\to 0$. This is a groupoid over $^{\bar d}\bar\cQ_1$, let $^{\bar d}\act: {^{\bar d}\cE}\to {^{\bar d}\bar\cQ_1}$ denote the action.

 Let $^{\bar d}\bar\cQ'_{0,ex}\subset {^{\bar d}\bar\cQ_{0,ex}}$ be the closed substack given by: $t_0$ comes from $\H^0(X, \cA\otimes\Omega^{-1}(-2D_1))$.
Any object of $\D^W(^{\bar d}\bar\cQ_{0,ex})$ is supported on $^{\bar d}\bar\cQ'_{0,ex}$. As in Appendix A.2, let
$$
^{\bar d}W_{1,0,ex}(K)=\Four(^{\bar d}\act^*K)[\dimrel](\frac{\dimrel}{2}),
$$
where $\dimrel$ is the relative dimension of $^{\bar d}\cE\to  {^{\bar d}\bar\cQ_1}$.
This functor satisfies the same properties as $^{\bar d}W_{k,k+1,ex}$ in Sect.~3.1.

\medskip
\noindent
3.3 For $k=1,2$ we single out the subgroupoids $\cH'_{k,y}\subset \cH_{k,y}$ as
follows.  

\medskip\noindent
{\bf CASE $k=1$.} \ 
We let $\cH'_{1,y}=\cH_{1,y}$ and $^i\cH'_{1,y}={^i\cH_{1,y}}$. 
 Recall the map $\cE_1\to\bar\cQ^y_1$ defined in Sect.~2.9. Write $\cE_1$ as a union
of vector bundles $^i\cE_1\to \bar\cQ^y_1$ with fibre
$$
(L_{-1}/\Omega)^*\otimes (\Omega(i y)/\Omega)
$$
The fibre of $^i\cE_1^*\to \bar\cQ^y_1$ is $(L_{-1}/\Omega)\otimes(\cO/\cO(-iy))$.

\medskip\noindent
{\bf CASE $k=2$.} \   A point of $\bar\cQ^y_2$ gives rise to a $P_2$-bundle
$\cF_{P_2}$ on
$D_y$ given by $(L_1\subset L_2\subset M)$ with $L_1\;\iso\;\Omega\mid_{D_y}$ and
$L_2/L_1\;\iso\;\cO\mid_{D_y}$. The $P_2$-bundle $\cF_{P_2}$ gives rise to a $P$-bundle
$\cF_P=\cF_{P_2}\times_{P_2} P$ on $D_y$. 

\smallskip 

By definition, the fibre of 
$
\cH'_{2,y}\to\bar\cQ^y_2
$
is
$$
\H^0(D^*_y, \;\cF_P\times_P U)/\H^0(D_y, \;\cF_P\times_P U)\;\iso\; 
(\Sym^2 L_2)\otimes(\cA^{-1}(\infty y)/\cA^{-1})
$$
Let $^i\cH'_{2,y}\subset \cH'_{2,y}$ be the subgroupoid with fibre
$$
(\Sym^2 L_2)\otimes(\cA^{-1}(iy)/\cA^{-1})
$$
Let $\cE_2\to \bar\cQ^y_2$ be the stack with fibre 
$$
(\Sym^2(L_2/L_1))\otimes(\cA^{-1}(\infty y)/\cA^{-1})\;\iso\; \cA^{-1}(\infty
y)/\cA^{-1}
$$
This is a union of vector bundles $^i\cE_2\to \bar\cQ^y_2$ with fibre
$\cA^{-1}(i y)/\cA^{-1}$.  The fibre of $^i\cE_2^*\to \bar\cQ^y_2$ is $\cA\otimes(\Omega/\Omega(-iy))$.

\bigskip
\noindent
3.4 For $k=1,2$ we have a natural map $\cH'_{k,y}\to \cE_k$ over $\bar\cQ^y_k$. Without loss of generality we may assume that the image of $^i\cH'_{k,y}$ in $\cE_k$ is $^i\cE_k$. The corresponding map $^ip_k: {^i\cH'_{k,y}}\to {^i\cE_k}$ is smooth with contractible fibres, we denote by $d_{i,k}$ its relative dimension. From (\cite{G}, Lemma~4.8) we get

\begin{Lm} 
\label{Lm_full_subcategory}
The functor $K\mapsto {^ip_k^*}K[d_{i,k}]$ is t-exact and identifies $\D({^i\cE_k})$ with
a full triangulated subcategory of $\D(^i\cH'_{k,y})$. \QED
\end{Lm}
 
 For $i'\ge i$ we have $^i\cE_k\hook{} {^{i'}\cE_k}$ is a subbundle.
Denote by $\pr_{i',i}: \,{^{i'}\cE^*_k}\to{^i\cE^*_k}$ the dual map.
  
\begin{Lm} 
\label{Lm_3}
For each $i\ge 0$ we have a natural map $f_i: \bar\cQ^y_{k+1,ex}\to {^i\cE^*_k}$
over $\bar\cQ^y_k$.
For $i'\ge i$ the composition 
$$
\bar\cQ^y_{k+1,ex}\,\toup{f_{i'}}\, {^{i'}\cE^*_k}\;\toup{\pr_{i',i}}\;
 {^i\cE^*_k}
$$ 
equals $f_i$. For each open substack of finite type $U\subset \bar\cQ^y_k$ there is an integer $i(U)$ such that over the preimage of $U$, the map $f_i: \bar\cQ^y_{k+1,ex}\to {^i\cE^*_k}$ is a closed embedding for every $i\ge i(U)$.
\end{Lm}
\begin{Prf}\\
{\bf CASE $k=1$.} \    Given a point of $\bar\cQ^y_{2,ex}$, over $D_y$ the section $t_2:\Omega\to \cW$ yields a map $s:\cO\to L_{-1}/\Omega$ such that $t_2=t_1\wedge s$. Now $f_i$ sends a point of $\bar\cQ^y_{2,ex}$ to the image of $s$ in $^i\cE^*_1$.

 Let $^i\cV\to \bar\cQ^y_1$ be the vector bundle with fibre
$\Hom(\Omega, \cW/\cW(-iy))$.  Given a point of $\bar\cQ^y_1$, we have a subbundle $\Omega\otimes (L_{-1}/\Omega)\mid_{D_y}\subset \cW\mid_{D_y}$ over $D_y$. Therefore, 
$$
(L_{-1}/\Omega)\otimes (\Omega/\Omega(-iy))\hook{}
\cW/\cW(-iy)
$$
So, we have a natural closed embedding $^i\cE^*_1\to {^i\cV}$ over $\bar\cQ^y_1$.

 Let $i(U)$ be such that the vector space $\Hom(\Omega, \cW(-iy))$ is zero for any point of $U$. Then for $i\ge i(U)$ the natural map $\bar\cQ^y_{2,ex}\to {^i\cV}$ is a closed embedding over $U$. So, $\bar\cQ^y_{2,ex}\to  {^i\cE^*_1}$ is a closed embedding over $U$ for $i\ge i(U)$.

\medskip\noindent
{\bf CASE $k=2$.} \    The map $f_i$ sends a point of $\bar\cQ^y_{3,ex}$ to the image of $t_3\in\H^0(X, \cA\otimes\Omega)$ in $\cA\otimes(\Omega/\Omega(-iy))$.

 Let $i(U)$ be such that the vector space $\H^0(X, \cA\otimes\Omega(-iy))$ is zero for any point of $U$. Then for $i\ge i(U)$ the map $f_i$ is a closed embedding.
\end{Prf}

\bigskip\noindent
3.5 For $k=1,2$ let $^i\act'_k: {^i\cH'_{k,y}}\to
\bar\cQ^y_k$
denote the action map. It is smooth, and we denote by $a_{i,k}$ its relative dimension.

Define the functor 
$$
W^{y,i}_{k,k+1,ex}: \D^W(\bar\cQ^y_k)\to \D(^i\cE^*_k)
$$
as follows. Given $K\in \D^W(\bar\cQ^y_k)$, from Lemma~\ref{Lm_full_subcategory} we learn that there exists $\tilde K\in \D(^i\cE_k)$ and an isomorphism
\begin{equation}
\label{iso_tilde_K}
h: {^ip_k^*\tilde K}[d_{i,k}](\frac{d_{i,k}}{2})\;\iso\;
(^i\act'_k)^*K[a_{i,k}](\frac{a_{i,k}}{2})
\end{equation}
The pair $(\tilde K,h)$ is defined up to a unique isomorphism. Set $W^{y,i}_{k,k+1,ex}(K)=\Four(\tilde K)$.

 By construction, the functor $W^{y,i}_{k,k+1,ex}$ is $t$-exact.

\medskip

\begin{Rem}  We could replace $^i\cH'_{k,y}$ by any subgroupoid
$^i\cH'_{k,y} \subset \cH'_{k,y}$ of finite type over $\bar\cQ^y_k$
such that the image of $^i\cH'_{k,y}$ in $\cE_k$ is $^i\cE_k$. The corresponding functors $W^{y,i}_{k,k+1,ex}:\D^W(\bar\cQ^y_k)\to \D(^i\cE^*_k)$ would be naturally isomorphic. Thus, the functors $W^{y,i}_{k,k+1,ex}$ do not depend on the choice of the group subschemes $^iN_{k,y}$ inside of $N_{k,y}$.
\end{Rem}

\smallskip

 Using the above remark together with appendix A.3, one shows that for $i'\ge i$ we have an isomorphism of functors $(\pr_{i',i})_!\comp W^{y,i'}_{k,k+1,ex}\;\iso\; W^{y,i}_{k;k+1,ex}$.

\begin{Lm} For $k=1,2$ let $K\in\D^W(\bar\cQ^y_k)$. 
For any open substack of finite type $U\subset \bar\cQ^y_k$ and any 
integer $i$ large enough (in particular, $i\ge i(U)$ of Lemma~\ref{Lm_3}), over the
preimage of $U$, the complex $W^{y,i}_{k,k+1,ex}(K)$ is supported on $\bar\cQ^y_{k+1,ex}\subset
{^i\cE^*_k}$.
\end{Lm}

\begin{Prf}
Since $U$ is contained in a finite number of strata $^{\bar d}\bar\cQ_k$, we are easily reduced to the case where $U\subset {^{\bar d}\bar\cQ^y_k}$ for some $\bar d$, and $K$ is the extension by zero from $U$. 

\smallskip
\noindent
{\bf CASE $k=1$.}   There is $i'(U)$ such that for any point of $U$ given by $D_1\in (X-y)^{(d_1)}$,  $(\Omega(D_1)\subset L_{-1}\subset M)\in\Bun_{P_1}$ we have
$\H^1(X, \; (L_{-1}/\Omega(D_1))^*\otimes\Omega(D_1+iy))=0$. 
So, for any $i\ge i'(U)$ the natural map 
$$
^i\cE_1\to \Ext^1(L_{-1}/\Omega(D_1), \Omega(D_1))
$$ 
is surjective over $U$. 
 If $i\ge i(U),i'(U)$ then $W^{y,i}_{1,2,ex}(K)$ is supported at $^{\bar d}\bar\cQ^y_{2,ex}$ and is isomorphic to $^{\bar d}W_{1,2,ex}(K)$.

\smallskip\noindent
{\bf CASE $k=2$.}  Recall that a point of $U$ is given by a collection: $D_1\in (X-y)^{(d_1)}$, $D'_2\in (X-y)^{(d_2-d_1)}$ with $D'_2\ge D_1$ and $(L_1\subset L_2\subset L_{-1}\subset M)\in\Bun_{P_2}$ with $L_1\,\iso\,\Omega(D_1)$ and $L_2/L_1\,\iso\,\cO(D'_2)$. There is $i'(U)$ such that for any point of $U$ as above we have
$$
\H^1(X, (\Sym^2(L_2/L_1))\otimes \cA^{-1}(iy))=0
$$
This implies that for $i\ge i'(U)$ the natural map 
$$
^i\cE_2\to \Ext^1(\cA, \,\Sym^2 (L_2/L_1))
$$ 
is surjective over $U$. If $i\ge i(U),i'(U)$ then $W^{y,i}_{2,3,ex}(K)$ is supported at $^{\bar d}\bar\cQ^y_{3,ex}$ and is isomorphic to $^{\bar d}W_{2,3,ex}(K)$.
\end{Prf}

\bigskip

 Thus, we get a well-defined functor $W^y_{k,k+1,ex}: \D^W(\bar\cQ^y_k)\to \D(\bar\cQ^y_{k+1,ex})$, it is $t$-exact by construction.

  Given $K\in \D^W(\bar\cQ^y_k)$, the $*$-restriction of
$W^y_{k,k+1,ex}(K)$ to $^{\bar d}\bar\cQ^y_{k+1,ex}$ is naturally isomorphic to 
$^{\bar d}W_{k,k+1,ex}$ applied to the $*$-restriction $K\mid_{^{\bar d}\bar\cQ^y_k}$.
By 1) of Lemma~\ref{Lm_on_stratifications}, we conclude that the image of $W^y_{k,k+1,ex}$ lies in $\D^W(\bar\cQ^y_{k+1,ex})$.

\begin{Pp} For $k=1,2$ the functor $K\mapsto (\pi_{k+1,k,ex})_!K$ maps 
$\D^W(\bar\cQ^y_{k+1,ex})$ to $\D^W(\bar\cQ^y_k)$ and is quasi-inverse to $W^y_{k,k+1,ex}$.
Moreover, for $K\in \D^W(\bar\cQ^y_{k+1,ex})$ the natural map $(\pi_{k+1,k,ex})_!K\to
(\pi_{k+1,k,ex})_*K$ is an isomorphism.
\end{Pp}

\smallskip

\begin{Prf}
First, let us show that for $K\in \D^W(\bar\cQ^y_k)$ we have $(\pi_{k+1,k,ex})_!W^y_{k,k+1,ex}(K)\,\iso\, K$ naturally.
 Indeed, over an open substack of finite type $U\subset \bar\cQ^y_k$ and $i$ large enough we have $W^y_{k,k+1,ex}(K)=\Four(\tilde K)$ and 
$$
(\pi_{k+1,k,ex})_!W^y_{k,k+1,ex}(K)\;\iso\; i_U^*\tilde K[d_{i,k}-a_{i,k}](\frac{d_{i,k}-a_{i,k}}{2}),
$$
where $\tilde K$ is that of (\ref{iso_tilde_K}), and $i_U: U\to {^i\cE_k}$ is the zero section. The equivariance property of $K$ implies that the RHS of the above formula is identified with $K\mid_U$.

\smallskip

 The fact that $K\mapsto (\pi_{k+1,k,ex})_!K$ maps $\D^W(\bar\cQ^y_{k+1,ex})$ to $\D^W(\bar\cQ^y_k)$ follows from Appendix A.1.

\smallskip

 Now let us show that for $K\in \D^W(\bar\cQ^y_{k+1,ex})$ we have
\begin{equation}
\label{iso_quasi_inverse}
W^y_{k,k+1,ex}(\pi_{k+1,k,ex})_!K\;\iso\; K
\end{equation} 
naturally. To establish this isomorphism over
the preimage of an open substack of finite type $U\subset \bar\cQ^y_k$, fix an integer $i$
large enough with respect to $U$.

 The groupoid $^i\cH'_{k,y}\to \bar\cQ^y_k$ lifts to $^i\cE^*_k$ in the sense of A.1.
In particular, we have a cartesian square
$$
\begin{array}{ccc}
^i\cH'_{k,y}\times_{\bar\cQ^y_k} {^i\cE^*_k} & \toup{\act} & ^i\cE^*_k\\
\downarrow\lefteqn{\scriptstyle{\id\times\pi_{\cE}}} &&
\downarrow\lefteqn{\scriptstyle{\pi_{\cE}}}\\ ^i\cH'_{k,y} & \toup{^i\act'_k} & \bar\cQ^y_k,
\end{array}
$$
where we used the projections to define the fibred product, and $\pi_{\cE}$ also denotes 
the projection. We may start with $K\in \D(^i\cE^*_k)$ that satisfies the
equivariance property $\act^* K\,\iso\, \pr_2^*K\otimes \chi^*\cL_{\psi}$, where $\chi$ is
the composition
$$
^i\cH'_{k,y}\times_{\bar\cQ^y_k} {^i\cE^*_k} \to {^i\cE_k}\times_{\bar\cQ^y_k} 
{^i\cE^*_k}\toup{\ev} \A^1
$$
(Actually, for $k=2$ the complex $K$ satisfies a stronger equivariance property with respect
to the action of $\cH_{2,y}$, which we don't need for the moment.) 

 Looking at one more cartesian square
$$
\begin{array}{ccc}
^i\cH'_{k,y}\times_{\bar\cQ^y_k} {^i\cE^*_k} & \to & ^i\cE_k\times_{\bar\cQ^y_k}
{^i\cE^*_k}\\
\downarrow\lefteqn{\scriptstyle{\id\times\pi_{\cE}}} && \downarrow\\
^i\cH'_{k,y} & \toup{^ip_k} & ^i\cE_k
\end{array}
$$
we obtain 
$$
(^i\act'_k)^*(\pi_{\cE})_!K\;\iso\; ^ip_k^*
\Four(K)[d_{i,k}-a_{i,k}](\frac{d_{i,k}-a_{i,k}}{2})
$$
We have used the fact that the rank
of the vector bundle $^i\cE_k\to \bar\cQ^y_k$ is $a_{i,k}-d_{i,k}$. The isomorphism
(\ref{iso_quasi_inverse}) over the preimage of $U$ follows.

 The above diagrams also show that
$$
(^i\act'_k)^*(\pi_{\cE})_!K\;\iso\;(^i\act'_k)^*(\pi_{\cE})_*K,
$$ 
because $!$- and
$*$-Fourier transforms coincide. So, $(\pi_{k+1,k,ex})_!K\to
(\pi_{k+1,k,ex})_*K$ is an isomorphism.
\end{Prf}

\bigskip

 Now arguing as in (\cite{G}, 5.11) one finishes the proof of Theorem~\ref{Th_1} ii).

\bigskip\noindent
3.6 The proof of Theorem~\ref{Th_1} i) is similar. First, let
$^i\cE_0={^i\cH_{0,y}}$. The action map $^i\act_0: {^i\cH_{0,y}}\to\bar\cQ^y_1$ is smooth,
denote by $a_{i,0}$ its relative dimension. For $i\ge 0$ define the functors 
$$
W^{y,i}_{1,0,ex}: \D(\bar\cQ^y_1)\to \D(^i\cE^*_0)
$$
by $W^{y,i}_{1,0,ex}(K)=\Four(^i\act_0^*K)[a_{i,0}](\frac{a_{i,0}}{2})$. As in Sect.~3.5,
this gives rise to a functor
$W^y_{1,0,ex}: \D(\bar\cQ^y_1)\to \D^W(\bar\cQ^y_{0,ex})$ and so on. The details are left to
the reader. \QED

\bigskip\bigskip
\centerline{\scshape 4. Cuspidality}

\bigskip\noindent
4.1 Recall the notion of cuspidality on $\Bun_G$. For a proper parabolic $Q\subset G$
let $M_Q$ be its Levi quotient. We have a diagram of natural maps
$$
\Bun_{M_Q}\getsup{\alpha_Q} \Bun_Q\toup{\beta_Q}\Bun_G
$$ 
The constant term functor $\CT_Q:
\D(\Bun_G)\to \D(\Bun_{M_Q})$ is defined as $\CT_Q(K)=(\alpha_Q)_!\beta_Q^*K$.

 A complex $K\in \D(\Bun_G)$ is \select{cuspidal} if $\CT_Q(K)=0$ for any standard proper
parabolic $P_2\subset Q\subset G$. It suffices to check this condition for $Q=P_1$ and $Q=P$. 
 
 Denote by $\D_{cusp}(\Bun_G)\subset\D(\Bun_G)$ the full triangulated subcategory consisting
of cuspidal objects. Similarly, for a scheme of parameters $S$, one defines
$\D_{cusp}(S\times \Bun_G)$.

\medskip\noindent
4.2 Let us introduce the notion of cuspidality on $\bar\cQ_k$ for $k=1,2,3$.

 The stack $\Bun_{M_1}$ classifies pairs: a line bundle $L_1$ on $X$ and a rank 2 bundle
$M_2$ on $X$. The projection $\Bun_{P_1}\to\Bun_{M_1}$ sends $(L_1\subset L_{-1}\subset
M)\in\Bun_{P_1}$ to $(L_1, M_2=L_{-1}/L_1)$. 

 The stack $\Bun_M$ classifies pairs: a line bundle $\cA$ on $X$ and a rank 2 bundle $L_2$
on $X$. The projection $\Bun_P\to\Bun_M$ sends a collection $(\cA, L_2, 0\to \Sym^2 L_2\to
?\to\cA\to 0)$ to $(\cA, L_2)$.

 For $k=1,2$ consider the natural diagram
$$
\begin{array}{ccccc}
\bar\cQ_k & \getsup{\beta_P^k} & \bar\cQ^P_k & \toup{\alpha_P^k} &
\bar\cQ^M_k\\
\downarrow && \downarrow && \downarrow \\
\Bun_G & \getsup{\beta_P} & \Bun_P & \toup{\alpha_P} & \Bun_M,
\end{array}
$$
where the right square is cartesian, and 
the stack $\bar\cQ^M_k$ classifies collections: an
$M$-torsor $(\cA, L_2)$ on $X$ together with sections $t_1,\ldots,t_k$, where
$$
\begin{array}{l}
t_1:\Omega\hook{} L_2\\
t_2:\Omega\hook{} \wedge^2 L_2
\end{array}
$$

 The constant term functor $\CT_P^{\bar\cQ_k}:\D(\bar\cQ_k)\to \D(\bar\cQ^M_k)$ is
defined as $\CT_P^{\bar\cQ_k}(K)=(\alpha_P^k)_!(\beta_P^k)^*K$.

 Consider the natural diagram
$$
\begin{array}{ccccc}
\bar\cQ_1 & \getsup{\beta_{P_1}^1} & \bar\cQ^{P_1}_1 & \toup{\alpha_{P_1}^1} &
\bar\cQ^{M_1}_1\\
\downarrow && \downarrow && \downarrow \\
\Bun_G & \getsup{\beta_{P_1}} & \Bun_{P_1} & \toup{\alpha_{P_1}} & \Bun_{M_1},
\end{array}
$$
where the right square is cartesian, and the stack $\bar\cQ^{M_1}_1$ classifies collections:
a $M_1$-torsor $(L_1, M_2)$ on $X$ together with section $t_1: \Omega\hook{} L_1$. 

\smallskip

 The constant term functor $\CT_{P_1}^{\bar\cQ_1}:\D(\bar\cQ_1)\to \D(\bar\cQ^{M_1}_1)$ is
defined as $\CT_{P_1}^{\bar\cQ_1}(K)=(\alpha_{P_1}^1)_!(\beta_{P_1}^1)^*K$.

\begin{Def} i) An object $K\in\D(\bar\cQ_1)$ is \select{cuspidal} if $\CT_P^{\bar\cQ_1}(K)=0$
and
$\CT_{P_1}^{\bar\cQ_1}(K)=0$.\\
ii) An object $K\in\D(\bar\cQ_2)$ is \select{cuspidal} if $\CT_P^{\bar\cQ_2}(K)=0$.\\
iii) Any object $K\in \D(\bar\cQ_3)$ is \select{cuspidal}.
\end{Def}

\noindent
4.3 For $k=1,2$ denote by $W_{k,k+1}: \D^W(\bar\cQ_k)\to \D^W(\bar\cQ_{k+1})$ the functor
$W_{k,k+1,ex}$ followed by the restriction to $\bar\cQ_{k+1}\subset \bar\cQ_{k+1,ex}$.

\begin{Pp} 
\label{Pp_cuspidal_to_cuspidal}
i) The functor $W_{k,k+1}: \D^W(\bar\cQ_k)\to \D^W(\bar\cQ_{k+1})$ maps cuspidal
objects to cuspidal.\\
ii) If $K\in\D^W(\bar\cQ_k)$ is cuspidal then the $*$-restriction of $W_{k,k+1,ex}(K)$ to
$\bar\cQ_{k+1,ex}-\bar\cQ_{k+1}$ vanishes.
\end{Pp}
\begin{Prf}
ii) Note that $\bar\cQ_{k+1,ex}-\bar\cQ_{k+1}$ is isomorphic to $\bar\cQ_k$, the zero
section of the bundle $\pi_{k+1,k,ex}:\bar\cQ_{k+1,ex}\to\bar\cQ_k$.
We will calculate the $*$-restriction $W_{k,k+1,ex}(K)\mid_{^{\bar d}\bar\cQ_k}$ for any
stratum $^{\bar d}\bar\cQ_k\subset \bar\cQ_k$.

\smallskip

 Let $\psi_1: {^{\bar d}\bar\cQ_1}\to \bar\cQ^{M_1}_1$ be the map that sends $((L_1\subset
L_{-1}\subset M), \;\Omega\hook{t_1} L_1)$ to 
$$
(M_2=L_{-1}/L_1, \; \Omega\hook{t_1} L_1)
$$
 Let $\psi_2: {^{\bar d}\bar\cQ_2}\to
\bar\cQ^M_2$ be the map that sends $((L_1\subset L_2\subset L_{-1}\subset M), \; 
\Omega\hook{t_1} L_1, \; \Omega\hook{t_2} \wedge^2 L_2)$ to 
$$
(\cA, L_2, \; 
\Omega\hook{t_1} L_2, \;\Omega\hook{t_2} \wedge^2 L_2)
$$ 

 Using Lemma~\ref{Lm_2}, one shows the following:
\begin{itemize}
\item for $K\in \D^W(\bar\cQ_1)$ we have
$W_{1,2,ex}(K)\mid_{^{\bar d}\bar\cQ_1}\;\iso\; \psi_1^*\CT_{P_1}^{\bar\cQ_1}(K)$ up to
a cohomological shift and a twist;

\item for $K\in \D^W(\bar\cQ_2)$ we have
$W_{2,3,ex}(K)\mid_{^{\bar d}\bar\cQ_2}\;\iso\; \psi_2^*\CT_P^{\bar\cQ_2}(K)$ up to
a cohomological shift and a twist.
\end{itemize}
Part ii) follows. 

\begin{Rem}
\label{Rem_iff_cuspidality}
Actually, we showed that for $K\in\D^W(\bar\cQ_k)$ the condition
$W_{k,k+1,ex}(K)\mid_{\bar\cQ_k}=0$ is equivalent to $\CT_{P_1}^{\bar\cQ_1}(K)=0$ for $k=1$
(resp., to $\CT_P^{\bar\cQ_2}(K)=0$ for $k=2$). Indeed, as $\bar d$ ranges over strings of
nonnegative integers $\bar d=(d_1,\ldots,d_k)$, the images of $\psi_k$ form a stratification
of the corresponding stack.
\end{Rem}

\medskip\noindent
i) Let $\bar\cQ^M_{2,ex}$ denote the stack classifying $(\cA,L_2)\in\Bun_M$ and sections
$\Omega\hook{t_1} L_2$, $\;\Omega\toup{t_2}\wedge^2 L_2$. As in Sect.~4.2, we have the
diagram
$$
\begin{array}{ccccc}
\bar\cQ_{2,ex} & \getsup{\beta_P^{2,ex}} & \bar\cQ^P_{2,ex} & \toup{\alpha_P^{2,ex}} &
\bar\cQ^M_{2,ex}\\
\downarrow && \downarrow && \downarrow \\
\Bun_G & \getsup{\beta_P} & \Bun_P & \toup{\alpha_P} & \Bun_M,
\end{array}
$$
where the right square is cartesian. Let $\CT^{\bar\cQ_{2,ex}}_P: \D(\bar\cQ_{2,ex})\to
\D(\bar\cQ^M_{2,ex})$ denote the functor $K\mapsto (\alpha_P^{2,ex})_!(\beta_P^{2,ex})^*K$.
 Proceeding as in Sect.~2-3, one introduces the category $\D^W(\bar\cQ^M_{2,ex})$ and the
functor 
$$
W^M_{1,2,ex}: \D(\bar\cQ^M_1)\to \D^W(\bar\cQ^M_{2,ex}),
$$ 
which is also an equivalence of categories.

 One checks that $\CT^{\bar\cQ_{2,ex}}_P$ sends $\D^W(\bar\cQ_{2,ex})$ to
$\D^W(\bar\cQ^M_{2,ex})$. Let us only indicate that the groupoid
$\cH_{1,y}\times_{\bar\cQ^y_1}\bar\cQ^y_{2,ex}\to \bar\cQ^y_{2,ex}$ lifts to
$$
\bar\cQ^{P,y}_{2,ex}=\bar\cQ^P_{2,ex}\times_{\bar\cQ_{2,ex}}
\bar\cQ^y_{2,ex}
$$
 We claim that there is a natural isomorphism of functors from $\D^W(\bar\cQ_1)$ to
$\D^W(\bar\cQ^M_{2,ex})$
\begin{equation}
\label{iso_functors_CT}
\CT^{\bar\cQ_{2,ex}}_P\comp
W_{1,2,ex}\;\iso\; W^M_{1,2,ex}\comp \CT^{\bar\cQ_1}_P
\end{equation}

 The functor $\CT^{\bar\cQ_1}_P$ admits a right adjoint, which will be denoted by
$\Eis^{\bar\cQ_1}_P$, it sends $K$ to $(\beta^1_P)_*(\alpha^1_P)^!K$. Actually,
$\CT^{\bar\cQ_1}_P$ maps $\D(\bar\cQ^M_1)$ to $\D^W(\bar\cQ_1)$.

 Similarly, $\CT^{\bar\cQ_{2,ex}}_P$ admits a right adjoint functor
$$
\Eis^{\bar\cQ_{2,ex}}_P: \D^W(\bar\cQ^M_{2,ex})\to \D^W(\bar\cQ_{2,ex})
$$ 
that sends $K$ to $(\beta^{2,ex}_P)_*(\alpha^{2,ex}_P)^!K$.

\smallskip

 We have the following diagram, where the right square is cartesian
$$
\begin{array}{ccccc}
\bar\cQ_1 & \getsup{\beta^1_P} & \bar\cQ^P_1 & \toup{\alpha^1_P} & \bar\cQ^M_1\\
\uparrow\lefteqn{\scriptstyle \pi_{2,1,ex}}
 && \uparrow && \uparrow\lefteqn{\scriptstyle \pi^M_{2,1,ex}}\\
\bar\cQ_{2,ex} & \getsup{\beta_P^{2,ex}} & \bar\cQ^P_{2,ex} & \toup{\alpha_P^{2,ex}} &
\bar\cQ^M_{2,ex}
\end{array}
$$ 
It follows that $(\pi_{2,1,ex})_*\comp \Eis^{\bar\cQ_{2,ex}}_P\;\iso\;
\Eis^{\bar\cQ_1}_P\comp (\pi^M_{2,1,ex})_*$ naturally. Passing to left adjoint functors, we
get the isomorphism (\ref{iso_functors_CT}).

 So, if $K\in\D^W(\bar\cQ_1)$ is cuspidal then $\CT^{\bar\cQ_{2,ex}}_P
W_{1,2,ex}(K)=0$. By i), the complex $W_{1,2,ex}(K)$ is the extension by zero from
$\bar\cQ_2$, so 
$\CT^{\bar\cQ_2}_P
W_{1,2}(K)=0$ and $W_{1,2}(K)$ is cuspidal. 
\end{Prf}

\bigskip
Recall the notation $\bar\cQ=\bar\cQ_3$. Let $W:\D^W(\bar\cQ_1)\to \D^W(\bar\cQ)$ be the functor $W_{2,3}\comp W_{1,2}$. Exactly as in (\cite{G}, Theorem~6.4), one derives from Proposition~\ref{Pp_cuspidal_to_cuspidal} the following corolary.

\begin{Cor} 
\label{Cor_cuspidal}
For $k=1,2$ let $K_1, K_2\in\D^W(\bar\cQ_k)$ be two objects with $K_1$ cuspidal.  Then the map $\Hom_{\D^W(\bar\cQ_k)}(K_1,K_2)\to
\Hom_{\D^W(\bar\cQ_{k+1})}(W_{k,k+1}(K_1),\; W_{k,k+1}(K_2))$ is an isomorphism.
So, for $k=1$
$$
\Hom_{\D^W(\bar\cQ_1)}(K_1,K_2)\to
\Hom_{\D^W(\bar\cQ)}(W(K_1),\; W(K_2))
$$
is also an isomorphism. \QED
\end{Cor}

\medskip\noindent
4.4 We also have the following analog of (\cite{G}, Theorem~6.9). For $k=1,2,3$ let
$\D^W_{cusp}(\bar\cQ_k)\subset \D^W(\bar\cQ_k)$ denote the full subcategory consisting of cuspidal objects. This is a triangulated subcategory.

\begin{Th} For $k=1,2$ the functor $W_{k,k+1}$ induces an equivalence of triangulated categories $\D^W_{cusp}(\bar\cQ_k)\to \D^W_{cusp}(\bar\cQ_{k+1})$. In particular,
$W: \D^W_{cusp}(\bar\cQ_1)\to \D^W(\bar\cQ)$ is an equivalence.
\end{Th}
\begin{Prf}
We know by Proposition~\ref{Pp_cuspidal_to_cuspidal} that $W_{k,k+1}$ maps cuspidal objects to cuspidal. 
Let 
$$
W^{-1}_{k,k+1}: \D^W_{cusp}(\bar\cQ_{k+1})\to
\D^W(\bar\cQ_k)
$$
be the functor sending $K$ to  $(\pi_{k+1,k,ex})_!K'$, where $K'$ is the extension by zero of $K$ to $\bar\cQ_{k+1,ex}$. 

\smallskip

 If $K\in  \D^W_{cusp}(\bar\cQ_{k+1})$ then the complex $W^{-1}_{k,k+1}(K)$ is cuspidal. Indeed, for $k=2$ the assertion follows from Remark~\ref{Rem_iff_cuspidality}. For $k=1$ set $F=W^{-1}_{1,2}(K)$.
We have $\CT^{\bar\cQ_1}_{P_1} F=0$ by Remark~\ref{Rem_iff_cuspidality}. Further, 
$W^M_{1,2,ex}\CT^{\bar\cQ_1}_P(F)=0$ by (\ref{iso_functors_CT}).  Since the functor 
$W^M_{1,2,ex}$ is an equivalence, we get $\CT^{\bar\cQ_1}_P(F)=0$.

 Let us show that $W^{-1}_{k,k+1}: \D^W_{cusp}(\bar\cQ_{k+1})\to \D^W_{cusp}(\bar\cQ_k)$ is quasi-inverse to $W_{k,k+1}$. From ii) of Theorem~\ref{Th_1} we conclude that $W_{k,k+1}\comp W^{-1}_{k,k+1}\;\iso\; \id_{\D^W_{cusp}(\bar\cQ_{k+1})}$ naturally, and there is a natural adjunction map 
$W^{-1}_{k,k+1}\comp W_{k,k+1}\to \id_{\D^W_{cusp}(\bar\cQ_k)}$.

 For $K\in \D^W_{cusp}(\bar\cQ_k)$ consider a distinguished triangle
$$
W^{-1}_{k,k+1} W_{k,k+1}(K)\to K\to K'
$$
We have $W_{k,k+1}(K')=0$ and $K'$ is cuspidal. Hence, $K'=0$ by 
Corolary~\ref{Cor_cuspidal}.
\end{Prf}

\bigskip\bigskip
\centerline{\scshape 5. Hecke functors}

\bigskip\noindent
5.1 \  Recall the Hecke functors $\H, \H^{\bar\cQ_{k,ex}}$ and $\H^{\bar\cQ_k}$ introduced in Sect.~2.2-2.4. 

\begin{Pp} 
\label{Pp_Hecke_preserves}
The functor $\H^{\bar\cQ_{k,ex}}$ sends $\D^W(\bar\cQ_{k,ex})$
to $\D^W(X\times\bar\cQ_{k,ex})$.
The functor $\H^{\bar\cQ_k}$ sends $\D^W(\bar\cQ_k)$ to
$\D^W(X\times\bar\cQ_k)$.
\end{Pp}

\begin{Prf}\   
Let $_x\H^{\bar\cQ_{k,ex}}$ denote the functor $\H^{\bar\cQ_{k,ex}}$ followed by
$*$-restriction to $x\times \bar\cQ_{k,ex}\subset X\times \bar\cQ_{k,ex}$. To simplify the
notation, we will show that $_x\H^{\bar\cQ_{k+1,ex}}$ preserves the category
$\D^W(\bar\cQ_{k+1,ex})$ for $k=1,2$. The other cases are treated similarly.

 Let $y\in X$ be distinct from $x$. Let $_x\cH_G$ be the preimage of $x$ under $\supp:\cH_G\to X$. We have a well-defined functor
$_x\H^{\bar\cQ^y_{k+1,ex}}: \D(\bar\cQ^y_{k+1,ex})\to \D(\bar\cQ^y_{k+1,ex})$. Let us show
that it preserves the subcategory $\D^W(\bar\cQ^y_{k+1,ex})$. Indeed, the groupoid
$$
\cH_{k,y}\times_{\bar\cQ^y_k}\bar\cQ^y_{k+1,ex}\to \bar\cQ^y_{k+1,ex}
$$ 
lifts to
$\bar\cQ^y_{k+1,ex}\times_{\Bun_G}{_x\cH_G}$ with respect to both $\gp_{k+1,ex}$ and
$\gq_{k+1,ex}$, so that we have diagrams

\smallskip

$$
\begin{array}{ccccc}
\cH_{k,y}\times_{\bar\cQ^y_k}\bar\cQ^y_{k+1,ex} & \getsup{\gp_{\cZ}} & \cZ &
\toup{\gq_{\cZ}} &
\cH_{k,y}\times_{\bar\cQ^y_k}\bar\cQ^y_{k+1,ex}\\
\downarrow\lefteqn{\scriptstyle\pr} && \downarrow\lefteqn{\scriptstyle\pr}
&& \downarrow\lefteqn{\scriptstyle\pr}\\
\bar\cQ^y_{k+1,ex} & \getsup{\gp_{k+1,ex}} & \bar\cQ_{k+1,ex}\times_{\Bun_G}{_x\cH_G} &
\toup{\gq_{k+1,ex}} & \cQ^y_{k+1,ex}
\end{array}
$$ 

and

$$
\begin{array}{ccccc}
\cH_{k,y}\times_{\bar\cQ^y_k}\bar\cQ^y_{k+1,ex} & \getsup{\gp_{\cZ}} & \cZ &
\toup{\gq_{\cZ}} &
\cH_{k,y}\times_{\bar\cQ^y_k}\bar\cQ^y_{k+1,ex}\\
\downarrow\lefteqn{\scriptstyle\act} && \downarrow\lefteqn{\scriptstyle\act}
&& \downarrow\lefteqn{\scriptstyle\act}\\
\bar\cQ^y_{k+1,ex} & \getsup{\gp_{k+1,ex}} & \bar\cQ_{k+1,ex}\times_{\Bun_G}{_x\cH_G} &
\toup{\gq_{k+1,ex}} & \cQ^y_{k+1,ex}
\end{array}
$$ 

\smallskip
\noindent
in both of which both squares are cartesian. Moreover, the compositions
$$
\cZ\;\toup{\gp_{\cZ}}\; \cH_{k,y}\times_{\bar\cQ^y_k}\bar\cQ^y_{k+1,ex} \;\toup{\chi_{k,y}}
\;\A^1
$$
and
$$
\cZ\;\toup{\gq_{\cZ}}\; \cH_{k,y}\times_{\bar\cQ^y_k}\bar\cQ^y_{k+1,ex} \;\toup{\chi_{k,y}}
\;\A^1
$$
coincide. Thus, $_x\H^{\bar\cQ^y_{k+1,ex}}$ preserves the equivariance condition.
Using the following remark, one finishes the proof.

\begin{Rem} Let $x_1,\ldots,x_n$ be a finite collection of points of $X$. Let
$K\in \D(\bar\cQ_{k,ex})$ be such that its restriction to $\bar\cQ_{k,ex}^y$ lies in
$\D^W(\bar\cQ^y_{k,ex})$ for any $y\ne x_i$. Then $K\in \D^W(\bar\cQ_{k,ex})$. 
Indeed, as $y$ ranges over points of $X$ different from $x_i$, the union of
$\bar\cQ^y_{k,ex}$ is $\bar\cQ_{k,ex}$.
Similar statement for $\D^W(\bar\cQ_k)$ holds. 
\end{Rem}
\noindent
\end{Prf}(Proposition~\ref{Pp_Hecke_preserves})

\bigskip\noindent
5.2 \  Similarly to the $\GL_n$ case, Hecke functors and Whittaker functors
commute with each other. The proof of the following result mimics that of (\cite{G},
Proposition~7.6).

\begin{Pp}
\label{Pp_Hecke_Whittaker_commute}
 i) For $k=1,2$ there is a natural isomorphism of functors
$$
\H^{\bar\cQ_{k+1,ex}}\comp W_{k,k+1,ex}\,\iso\, (\id\times W_{k,k+1,ex})\comp
\H^{\bar\cQ_k}: \; \D^W(\bar\cQ_k)\to \D^W(X\times \bar\cQ_{k+1,ex})
$$
ii) There is a natural isomorphism of functors
$$
\H^{\bar\cQ_{0,ex}}\comp W_{1,0,ex}\,\iso\, (\id\times W_{1,0,ex})\comp \H^{\bar\cQ_1}:
\; \D(\bar\cQ_1)\to \D^W(X\times\bar\cQ_{0,ex})
$$
\end{Pp}
\begin{Prf}
i) To simplify the notation, we replace the functors $\H^{\bar\cQ_k},
\H^{\bar\cQ_{k,ex}}$ by $_x\H^{\bar\cQ_k},\; {_x\H^{\bar\cQ_{k,ex}}}$.
 In view of Theorem~\ref{Th_1}, it suffices to show that for $K\in \D^W(\bar\cQ_{k+1,ex})$
we have 
$$
_x\H^{\bar\cQ_k}((\pi_{k+1,k,ex})_!K)\;\iso\; (\pi_{k+1,k,ex})_!
{_x\H^{\bar\cQ_{k+1,ex}}}(K)
$$

 For $x\in X$ let $\bar\cQ_{k+1,ex,x}$ be the stack defined in the same way as
$\bar\cQ_{k+1,ex}$ with the difference that the last map $t_{k+1}$ is allowed to have a
pole of order 2 at $x$ for $k=1$ (resp., of order 1 at $x$ for $k=2$). 

 Write
$_x\cH^{\bar\cQ_k}$ for the preimage of $x\times\bar\cQ_k$ under $\supp\times\gp_k:
\bar\cQ_k\times_{\Bun_G}\cH_G\to X\times\bar\cQ_k$. 
 We have a diagram
$$
\begin{array}{ccccc}
\bar\cQ_{k+1,ex,x} & \getsup{\gp_{k+1,ex,x}} & _x\cH^{\bar\cQ_{k+1,ex,x}} &
\toup{\gq_{k+1,ex}} & \bar\cQ_{k+1,ex}\\
\downarrow && \downarrow && \downarrow\lefteqn{\scriptstyle\pi_{k+1,k,ex}}\\
\bar\cQ_k & \getsup{\gp_k} & _x\cH^{\bar\cQ_k} & \toup{\gq_k} & \bar\cQ_k,
\end{array}
$$
where the stack $_x\cH^{\bar\cQ_{k+1,ex,x}}$ is defined by the condition that the right
square is cartesian, and $\gp_{k+1,ex,x}$ is the natural map.

 It suffices to show that for $K\in \D^W(\bar\cQ_{k+1,ex})$ the complex
$(\gp_{k+1,ex,x})_!\gq_{k+1,ex}^*K$ is supported on $\bar\cQ_{k+1,ex}\subset
\bar\cQ_{k+1,ex,x}$. This direct image will verify an appropriate equivariance condition
on $\bar\cQ_{k+1,ex,x}$. So, our assertion is verified stratum by stratum using
an analog of Lemma~\ref{Lm_2}. 

 Part ii) is proved similarly.
\end{Prf}

\bigskip\bigskip

\centerline{\scshape 6. Hyper-cuspidality}

\bigskip\noindent
6.1 Recall that $U_0$ denotes the center of $U_1$. Set $P_0=P_1/U_0$. We have a diagram of natural maps 
$\Bun_{P_0}\getsup{\alpha_{P_0}}\Bun_{P_1}\toup{\beta_{P_1}}\Bun_G$. Define the constant
term functor 
$$
\CT_{P_0}: \D(\Bun_G)\to\D(\Bun_{P_0})
$$ 
as $\CT_{P_0}(K)=(\alpha_{P_0})_!\beta_{P_1}^*(K)$. The following is a geometric
version of (\cite{PSh}, Definition on p. 328).

\begin{Def} i) A complex $K\in\D(\Bun_G)$ is \select{hyper-cuspidal} if $\CT_{P_0}(K)=0$.\\
ii) A complex $K\in\D(\bar\cQ_1)$ is \select{hyper-cuspidal} if $W_{1,0,ex}(K)$ is the
extension by zero from $\bar\cQ_0$.
\end{Def}

 Denote by $\D_{hcusp}(\Bun_G)\subset\D(\Bun_G)$ and by
$\D_{hcusp}(\bar\cQ_1)\subset\D(\bar\cQ_1)$ the full triangulated
subcategories consisting of hyper-cuspidal objects. Similarly, we have
$\D_{hcusp}(S\times\Bun_G)$ for a scheme of parameters $S$.

 If $f:S_1\to S_2$ is a morphism of schemes then for the map $f\times\id:
S_1\times\Bun_G\to S_2\times\Bun_G$ the functors $(f\times\id)_!$ and $(f\times\id)^*$
preserve hyper-cuspidality (and cuspidality). The same is true for the functor $\D(S)\times
\D(S\times\Bun_G)\to \D(S\times\Bun_G)$ of the tensor product along $S$.

\begin{Pp}
\label{Pp_hcusp_inclusion}
 In both cases $\D_{hcusp}(\Bun_G)\subset\D_{cusp}(\Bun_G)$ and
$\D_{hcusp}(\bar\cQ_1)\subset\D_{cusp}(\bar\cQ_1)$ is a full
triangulated subcategory.
\end{Pp}
\begin{Prf} Let $K\in\D_{hcusp}(\Bun_G)$. It is clear that $\CT_{P_1}(K)=0$. Let us show that $\CT_P(K)=0$.
We have a diagram
$$
\begin{array}{ccccc}
\Bun_G & \gets & \Bun_P & \gets & \Bun_{P_2}\\
&& \downarrow && \downarrow\\
&& \Bun_M & \gets & \Bun_{B(M)},
\end{array}
$$
where $B(M)\subset M$ is a Borel subgroup, the square is cartesian, and the composition in the top line is 
$\beta_{P_2}$. The right vertical arrow factors as $\Bun_{P_2}\toup{\delta}
\Bun_{P_2/U_0}\to\Bun_{B(M)}$. So, it is enough to show that
$\delta_!\beta_{P_2}^*(K)=0$.

 Since we have the following diagram, where the square is cartesian
$$
\begin{array}{ccccc}
\Bun_G & \gets & \Bun_{P_1} & \gets & \Bun_{P_2}\\
&& \downarrow && \downarrow\lefteqn{\scriptstyle\delta}\\
&& \Bun_{P_0} & \gets & \Bun_{P_2/U_0},
\end{array}
$$
the first assertion follows. 

 For sheaves on $\bar\cQ_1$ the proof is similar.
\end{Prf}

\bigskip\noindent
6.2 Recall that for each dominant coweight $\lambda$ of $G$ we have the Hecke functor
$\H^{\lambda}_G: \D(\Bun_G)\to\D(X\times\Bun_G)$ normalized to commute with Verdier duality
(cf. \cite{BG}, Sect.~2.1.4 for the precise definition). In our notation $\H^{\gamma}_G=\H$.
It is well-known that the subcategory $\D_{cusp}(\Bun_G)\subset\D(\Bun_G)$ is preserved
by Hecke functors. That is, each $\H^{\lambda}_G$ sends $\D_{cusp}(\Bun_G)$ to the category
$\D_{cusp}(X\times\Bun_G)$.

\begin{Pp} The subcategory $\D_{hcusp}(\Bun_G)\subset \D(\Bun_G)$ is preserved by Hecke
functors.
\end{Pp}
\begin{Prf}
\Step 1 Let us show that $\H$ preserves $\D_{hcusp}(\Bun_G)$.
One may introduce a version of stacks $\bar\cQ_1$ and $\bar\cQ_{0,ex}$, where instead
of a fixed $T$-torsor with trivial conductor $(\cF_T,\tilde\omega)$ one considers all of
them as additional parameter. In other words, the stack $\bar\cQ_1$ would classify
$\cF_G\in\Bun_G$, a line bundle $\cB$ on $X$ and a section $t_1:\cB\hook{} M$; the stack
$\bar\cQ_{0,ex}$ would classify the data just above together with
$\cB^2\otimes\Omega^{-1}\to\cA$.  

An analog of Theorem~\ref{Th_1} would hold in this setting.
Then for $K\in\D(\Bun_G)$ hyper-cuspidality would be equivalent to
requiring that $W_{1,0,ex}(\alpha^*K)$ is the extension by zero from $\bar\cQ_0$. Our
assertion follows from an analog of Proposition~\ref{Pp_Hecke_Whittaker_commute} ii) in this
situation.

\smallskip

\Step 2 Recall that Hecke functors can be composed in the following way. For $G$-dominant
coweights $\lambda_1,\lambda_2$ the functor 
$$
\H^{\lambda_1}_G\star \H^{\lambda_2}_G:
\D(\Bun_G)\to \D(X\times\Bun_G)
$$ 
is defined as 
$$
\H^{\lambda_1}_G\star \H^{\lambda_2}_G(K)=(\trian^*_X\boxtimes\id)
((\id\boxtimes \H^{\lambda_1})\comp\H^{\lambda_2}_G(K))[-1](\frac{-1}{2})
$$
It is known (\cite{BG} Sect. 2.1.6 and \cite{BD}) that there is a canonical
isomorphism functorial in
$K$ 
$$
\H^{\lambda_1}_G\star \H^{\lambda_2}_G(K)\;\iso\;
\mathop{\oplus}\limits_{\lambda\in\Lambda^+}\H^{\lambda}_G(K)\otimes\Hom_{\check{G}}(V^{\lambda},
V^{\lambda_1}\otimes V^{\lambda_2})
$$

 The group of coweight of $G$ orthogonal to all roots
is free abelian of rank 1. It is easy to see that Hecke functors corresponding to both
generators $\pm \omega$ of this group preserve $\D_{hcusp}(\Bun_G)$. One checks that any
irreducible representation $V^{\lambda}$ of $\check{G}$ appears in $(V^{\gamma})^{\otimes
k}\otimes V^{r\omega}$ for some $k\ge 0$ and $r\in\ZZ$. Thus, our assertion follows from the
fact that the subcategory $\D_{hcusp}(\Bun_G)$ is saturated: a direct summand of an object
of $\D_{hcusp}(\Bun_G)$ is again an object of $\D_{hcusp}(\Bun_G)$.
\end{Prf}

\medskip

 Clearly, the functor $\H^{\bar\cQ_1}$ preserves the subcategory $\D_{hcusp}(\bar\cQ_1)$,
and $\alpha^*$ sends $\D_{hcusp}(\Bun_G)$ to $\D_{hcusp}(\bar\cQ_1)$. 

\begin{Pp} 
\label{Pp_cusp_and_hcusp}
We have equivalences of triangulated categories

\medskip\noindent
i) $\D(\bar\cQ_1)/\D_{hcusp}(\bar\cQ_1)\;\iso\; \D^W(\bar\cQ_1)$

\medskip\noindent
ii) $\D_{cusp}(\bar\cQ_1)/\D_{hcusp}(\bar\cQ_1)\;\iso\; \D^W_{cusp}(\bar\cQ_1)$.
\end{Pp}

\begin{Rem}
\label{Rem_quotients}
 Let $F: D\to D'$ be a triangulated functor between triangulated categories. If
$F$ admits a fully faithfull right adjoint functor $F':D'\to D$ then $F$ induces an
equivalence of triangulated categories $D/\Ker F\to D'$. 
\end{Rem}

\medskip

\begin{Prf}\select{of Proposition~\ref{Pp_cusp_and_hcusp}}\\
i) Recall the closed immersion $i_0:\bar\cQ_1\hook{}\bar\cQ_{0,ex}$ and its complement
$j: \bar\cQ_0\hook{} \bar\cQ_{0,ex}$. The functor $i_0^*: \D^W(\bar\cQ_{0,ex})\to
\D^W(\bar\cQ_1)$ adimts a right adjoint $(i_0)_*$, which is fully faithfull. 
The category $\D^W(\bar\cQ_0)$ is embedded in $\D^W(\bar\cQ_{0,ex})$ fully faithfully by
$j_!$. By Remark~\ref{Rem_quotients}, 
$i^*$ induces an equivalence of triangulated categories
$$
\D^W(\bar\cQ_{0,ex})/\D^W(\bar\cQ_0)\;\iso\; \D^W(\bar\cQ_1)
$$
So, the functor $i_0^*\comp W_{0,1,ex}$ induces an equivalence i).

\smallskip\noindent
ii) For $K\in\D_{cusp}(\bar\cQ_1)$ let us show that $i_0^* W_{0,1,ex}(K)$ is cuspidal.
We have a distinguished triangle in $\D(\bar\cQ_1)$
$$
(\pi_{0,1})_! j^* W_{0,1,ex}(K)\to K \to i_0^*W_{0,1,ex}(K)
$$
Since $(\pi_{0,1})_! j^* W_{0,1,ex}(K)$ is hyper-cuspidal, it is cuspidal 
by Proposition~\ref{Pp_hcusp_inclusion}. So, $i_0^*W_{0,1,ex}(K)$ is also
cuspidal.

 We conclude that $i_0^*\comp W_{0,1,ex}$ induces a functor
$F: \D_{cusp}(\bar\cQ_1)/\D_{hcusp}(\bar\cQ_1)\;\to\; \D^W_{cusp}(\bar\cQ_1)$. Let
$F^{-1}$ denote the composition 
$$
\D^W_{cusp}(\bar\cQ_1)\to \D_{cusp}(\bar\cQ_1)\to
\D_{cusp}(\bar\cQ_1)/\D_{hcusp}(\bar\cQ_1)
$$
We claim that $F$ and $F^{-1}$ are quasi-inverse to each other. Indeed, the above
distinguished triangle shows that $\id\,\iso\, F^{-1}\comp F$. 
Since for $K\in \D^W(\bar\cQ_1)$ we have $W_{0,1,ex}(K)\iso (i_0)_*K$ naturally, it follows
that $\id\,\iso\, F\comp F^{-1}$. 
\end{Prf}

\bigskip\noindent
6.3 \  If $D'$ is a triangulated category and $D\subset D'$ is a full triangulated
subcategory, we write $D^{\perp}\subset D'$ for the full subcategory consisting of $K\in D'$
such that
$\Hom_{D'}(L,K)=0$ for all $L\in D$. Then $D^{\perp}\subset D'$ is a full triangulated
subcategory, and the composition $D^{\perp}\to D'\to D'/D$ is fully faithfull (cf. \cite{V},
Proposition 2.3.3, p.128). 

 Consider the subcategory $\D_{hcusp}(\Bun_G)^{\perp}\subset \D_{cusp}(\Bun_G)$. Let $_x
\H^{\lambda}_G:\D(\Bun_G)\to\D(\Bun_G)$ denote the functor $\H^{\lambda}_G$ followed by
$*$-restriction to $x\times\Bun_G\hook{} X\times\Bun_G$. Since Hecke functors admit left and
right adjoint functors (cf.\cite{BG}, 3.2.4), it follows that $\D_{hcusp}(\Bun_G)^{\perp}$
is preserved by all functors $_x\H^{\lambda}_G$.

\bigskip\bigskip
\centerline{\scshape 7. More Whittaker type functors}

\bigskip\noindent
7.1 Let $\cZ_1$ be the stack of collections: $(M,\cA)\in\Bun_G$ together with an isotropic
subsheaf $L_2\subset M$, where $L_2\in\Bun_2$. The stack $\cZ_1$ is nothing but $\Bunt_P$ in
the notation of (\cite{BG}, 1.3.6).

 Let $\pi_{2,1,ex}:\cZ_{2,ex}\to\cZ_1$ be the stack over $\cZ_1$ with fibre 
consisting of all maps
\begin{equation}
\label{map_s}
s:\Omega^{-1}\to\cA\otimes\Sym^2 L_2^*
\end{equation}
Let $\pi_{2,1}:\cZ_2\to \cZ_1$ be the open substack of $\cZ_{2,ex}$ given by the condition:
$s$ is injective.

 For $k=1,2$ we have the diagram 
$$
\cZ_k\getsup{\gp_k} \cZ_k\times_{\Bun_G}\cH_G \toup{\gq_k}\cZ_k,
$$
where we used the map $\gp:\cH_G\to \Bun_G$ in the definition of the fibred product, $\gp_k$
is the projection, and $\gq_k$ sends a point of $\cZ_k\times_{\Bun_G}\cH_G$ to $\cF'_G$
equiped with an isotropic subsheaf (and for $k=2$ a section $s'$) that are the compositions
$$
\begin{array}{l}
L_2\hook{} M\hook{} M'\\
s': \Omega^{-1}\hook{} \cA\otimes\Sym^2 L_2^*\hook{} \cA'\otimes\Sym^2 L_2^*
\end{array}
$$

 For $k=1,2$ we have the functor $\H^{\cZ_k}:\D(\cZ_k)\to\D(X\times\cZ_k)$ given by
$$
\H^{\cZ_k}(K)=(\supp\times\gp_k)_!\gq_k^*K\otimes \Qlb(\frac{1}{2})[1]^{\otimes \<\gamma,
2\check{\rho}\>} 
$$
Similarly, one defines the functor $\H^{\cZ_{2,ex}}:\D(\cZ_{2,ex})\to\D(X\times\cZ_{2,ex})$.

 The projection $\alpha_{\cZ}:\cZ_1\to \Bun_G$ fits into the diagram
$$
\begin{array}{ccccc}
\cZ_1 & \getsup{\gp_1} & \cZ_1\times_{\Bun_G}\cH_G & \toup{\gq_1} & \cZ_1\\
\downarrow\lefteqn{\scriptstyle \alpha_{\cZ}} && \downarrow &&
\downarrow\lefteqn{\scriptstyle
\alpha_{\cZ}}\\
\Bun_G & \getsup{\gp} & \cH_G & \toup{\gq} & \Bun_G
\end{array}
$$
So, $(\id\times\alpha_{\cZ})^*\comp \H\,\iso\, \H^{\cZ_1}\comp
\alpha_{\cZ}^*[1](\frac{1}{2})$ naturally.

 In this normalization the Hecke property on $\cZ_k$ (for $k=1,2$) with respect to $\H^{\cZ_k}$ 
and a given  local system $W$ on $X$ writes
$$
\H^{\cZ_k}(K)\,\iso\, W\boxtimes K [3](\frac{3}{2})
$$

\medskip\noindent
7.2 One defines the category $\D^W(\cZ_{2,ex})$ as in Sect. 2.10-2.11. Let us just indicate
its description on strata (they are equivariant under the corresponding groupoids).

 For $d\ge 0$ let $^d\cZ_1\subset\cZ_1$ be the locally closed substack given by: there is a
subbundle $L'_2\subset M$ such that $L_2\subset L'_2$ is a subsheaf with $d=\deg(L'_2/L_2)$.
The stack $^d\cZ_1$ classifies collections: a modification of rank 2 bundles $L_2\subset
L'_2$ on $X$, and an exact sequence $0\to\Sym^2L'_2\to ?\to\cA\to 0$, where $\cA$ is a line
bundle on $X$. 

Let $^d\cZ_{2,ex}=\cZ_{2,ex}\times_{\cZ_1}{^d\cZ_1}$. An analog of
Lemma~\ref{Lm_on_stratifications} holds for this stratification of $\cZ_{2,ex}$, so it
suffices to describe the categories $\D^W(^d\cZ_{2,ex})$ for each $d$.

 Let $^d\cZ'_{2,ex}\hook{} {^d\cZ_{2,ex}}$ be the closed substack given by the condition:
$s$ factors as 
$$
\Omega^{-1}\to \cA\otimes\Sym^2 L'^*_2\hook{} \cA\otimes\Sym^2 L^*_2
$$
Let $^d\chi_{2,ex}: {^d\cZ'_{2,ex}}\to\A^1$ be the map that pairs $s$ with the extension
$0\to \Sym^2 L'_2\to ?\to\cA\to 0$. 

 Let $^d\cP_{2,ex}$ be the stack classifying: a modification of rank 2 bundles $L_2\subset
L'_2$ on $X$ with $d=\deg(L'_2/L_2)$, a line bundle $\cA$ on $X$ and a section
$s: \Sym^2L'_2\to \Omega\otimes\cA$. Let
$$
\phi_{2,ex}:{^d\cZ'_{2,ex}}\to {^d\cP_{2,ex}}
$$ 
be the projection.

\begin{Lm} 
\label{Lm_useful_in_fact}
Any object of $\D^W(^d\cZ_{2,ex})$ is supported at $^d\cZ'_{2,ex}$. The functor  
$$
^dJ(K)={^d\chi_{2,ex}^*}\cL_{\psi}\otimes \phi^*_{2,ex}K[1](\frac{1}{2})^{\otimes\dimrel}
$$ 
provides an equivalence of categories $^dJ: \D(^d\cP_{2,ex})\,\iso\, \D^W(^d\cZ_{2,ex})$. Here $\dimrel$ is a function of a connected component of $^d\cZ'_{2,ex}$ given by $\dimrel=-\chi(\cA^{-1}\otimes\Sym^2 L'_2)$. 
\QED
\end{Lm}

 One mimics the proof of Theorem~\ref{Th_1} to get

\begin{Th} 
\label{Th_3}
There is an equivalence of categories $W\! Z_{1,2,ex}: \D(\cZ_1)\,\iso\,
\D^W(\cZ_{2,ex})$, which is $t$-exact, and $(\pi_{2,1,ex})_!$ is quasi-inverse to it.
Moreover, for any $K\in\D^W(\cZ_{2,ex})$ the natural map $(\pi_{2,1,ex})_!K\to
(\pi_{2,1,ex})_*K$ is an isomorphism. \QED
\end{Th}

 Let us just explain what this functor does on strata. We have the functor 
$$
^dW\! Z_{1,2,ex}: \D(^d\cZ_1)\to \D^W(^d\cZ_{2,ex})
$$
defined as the composition
$$
\D(^d\cZ_1)\toup{\Four} \D(^d\cP_{2,ex})\,\toup{^dJ}\,\D^W(^d\cZ_{2,ex})
$$

 If $K\in\D(\cZ_1)$ is the extension by zero from $^d\cZ_1$ then $W\! Z_{1,2,ex}(K)$ is the
extension by zero of $^dW\! Z_{1,2,ex}(K)$ under $^d\cZ_{2,ex}\hook{} \cZ_{2,ex}$.

\bigskip\noindent
7.3 {\scshape Relation to hyper-cuspidality}  \  Denote by $\cZ^{P_0}_1$ the stack of
collections: $P_0$-torsor on $X$, that is, an exact sequence $0\to L_1\to L_{-1}\to
L_{-1}/L_1\to 0$ on $X$ with $L_1\in\Bun_1$, $L_{-1}\in\Bun_3$; and an isotropic subsheaf
$L_2\subset L_{-1}$ with $L_2\in\Bun_2$. Here `isotropic' means that the composition
$\wedge^2 L_2\to \wedge^2 L_{-1}\to \wedge^2(L_{-1}/L_1)$ vanishes.

  Denote by
$$
\pi^{P_0}_{2,1,ex}:\cZ_{2,ex}^{P_0}\to \cZ_1^{P_0}
$$ 
the stack over $\cZ_1^{P_0}$ with fibre consisting of all maps (\ref{map_s}), where
$\cA=\det(L_{-1}/L_1)$.
 
 We have a natural diagram
$$
\begin{array}{ccccc}
\cZ_{2,ex} & \getsup{\beta^{2,ex}_{P_1}} & \cZ_{2,ex}^{P_1} & \toup{\alpha_{P_0}^{2,ex}} &
\cZ_{2,ex}^{P_0}\\
\downarrow && \downarrow && \downarrow\lefteqn{\scriptstyle \pi^{P_0}_{2,1,ex}}\\
\cZ_1 & \getsup{\beta^1_{P_1}} & \cZ_1^{P_1} & \toup{\alpha_{P_0}^1} & \cZ_1^{P_0}\\
\downarrow && \downarrow && \downarrow \\
\Bun_G & \getsup{\beta_{P_1}} & \Bun_{P_1} & \toup{\alpha_{P_0}} & \Bun_{P_0},
\end{array}
$$
where both right squares are cartesian (thus defining the stacks in the middle column).

The constant term functor 
$$
\CT^{\cZ_1}_{P_0}:
\D(\cZ_1)\to \D(\cZ_1^{P_0})
$$ 
is defined by
$\CT^{\cZ_1}_{P_0}(K)=(\alpha_{P_0}^{\cZ})_!(\beta^{\cZ}_{P_1})^*K$. Similarly,
$\CT^{\cZ_{2,ex}}_{P_0}: \D(\cZ_{2,ex})\to \D(\cZ_{2,ex}^{P_0})$ is defined as
$$
\CT^{\cZ_{2,ex}}_{P_0}(K)=(\alpha_{P_0}^{2,ex})_!(\beta^{2,ex}_{P_1})^*K
$$

\begin{Def} A complex $K\in \D(\cZ_1)$ (resp., $K\in \D^W(\cZ_{2,ex})$) is \select{hyper-cuspidal} if
$\CT^{\cZ_1}_{P_0}(K)=0$ (resp., $\CT^{\cZ_{2,ex}}_{P_0}(K)=0$). We denote by $\D_{hcusp}(\cZ_1)\subset \D(\cZ_1)$ and $\D^W_{hcusp}(\cZ_2)\subset \D^W(\cZ_{2,ex})$ the full
triangulated subcategories of hyper-cuspidal objects.
\end{Def}

 Clearly, $K\in\D(\Bun_G)$ is hyper-cuspidal iff $\alpha^*_{\cZ}K\in \D(\cZ_1)$ is
hyper-cuspidal. The following is easy to prove.

\begin{Pp} 1) A complex $K\in\D^W(\cZ_{2,ex})$ is hyper-cuspidal if and only if the following holds:
for any $k$-point $z=(L_2\subset M, \,s:\Sym^2 L_2\to  \cA\otimes\Omega)$ such that $L_2$ has a rank 1 isotropic subbundle (with respect to the form $s$) we have $K_z=0$. 

\smallskip\noindent
2) The functor $WZ_{1,2,ex}: \D(\cZ_1)\,\iso\,
\D^W(\cZ_{2,ex})$ induces an equivalence of triangulated categories
$$
\D_{hcusp}(\cZ_1)\,\iso\,\D^W_{hcusp}(\cZ_2)\eqno{\square}
$$
\end{Pp}

\begin{Rems} i) For each integer $d$ we have a closed substack $Y_d\hook{}\cZ_{2,ex}$ given by the condition that $L_2$ admits an isotropic rank 1 subbundle (with respect to $s$) of degree $\ge d$. We have $Y_d\subset Y_{d-1}\subset\ldots$. A complex $K\in\D^W(\cZ_{2,ex})$ is hyper-cuspidal if and only if its $*$-restriction to each $Y_d$ vanishes.\\
ii) If $s:\Sym^2 L_2\to\cA\otimes\Omega$ is such that $L_2$ has no rank 1 isotropic subbundles then the form $s$ is generically nondegenerate, that is, $L_2\hook{} L_2^*\otimes\cA\otimes\Omega$ is an inclusion.
\end{Rems}

\medskip

  Hecke functors preserve our equivariance conditions as well as hyper-cuspidality.
Moreover, they commute with $WZ_{1,2,ex}$, namely as in Sect.~5 one proves

\begin{Pp} 1) The functor $\H^{\cZ_{2,ex}}$ sends $\D^W(\cZ_{2,ex})$ to
$\D^W(X\times\cZ_{2,ex})$ and $\D^W_{hcusp}(\cZ_2)$ to $\D^W_{hcusp}(X\times \cZ_2)$.

\smallskip\noindent
2) The functor $\H^{\cZ_1}$ sends $\D_{hcusp}(\cZ_1)$ to $\D_{hcusp}(X\times\cZ_1)$.

\smallskip\noindent
3) We have a canonical isomorphism of functors $\H^{\cZ_{2,ex}}\comp 
W\! Z_{1,2,ex}\,\iso\,(\id\times W\! Z_{1,2,ex})\comp\H^{\cZ_1}$ from $\D(\cZ_1)$ to
$\D^W(X\times\cZ_{2,ex})$. \QED
\end{Pp}

\medskip\noindent
7.4  \ In this subsection we prove the following generalization of (\cite{G}, Theorem~7.9).

\begin{Pp} The functor $\H^{\cZ_2}: \D(\cZ_2)\to \D(X\times\cZ_2)$ is right-exact for the perverse t-structures.
\end{Pp}

  Let $\pi_{3,2}: \cZ_3\to\cZ_2$ denote the stack classifying $(L_2\subset M, \; s:\Omega^{-1}\hook{}\cA\otimes\Sym^2 L_2^*)\in\cZ_2$ together with a line subbundle $L_1\subset L_2$ such that 
$$
\H^1(X, L_1^{-1}\otimes (L_2/L_1))=0
$$ 
The projection $\pi_{3,2}$ is smooth and surjective. Consider the diagram
$$
\begin{array}{ccccc}
\cZ_3 & \getsup{\gp_3} & \cZ_3\times_{\Bun_G}\cH_G & \toup{\gq_3} & \cZ_3\\
\downarrow\lefteqn{\scriptstyle \pi_{3,2}} &&  \downarrow
&& \downarrow\lefteqn{\scriptstyle \pi_{3,2}}\\
\cZ_2 & \getsup{\gp_2} & \cZ_2\times_{\Bun_G}\cH_G & \toup{\gq_2} & \cZ_2,
\end{array}
$$
where the left square is cartesian. Define $\H^{\cZ_3}: \D(\cZ_3)\to \D(X\times\cZ_3)$ by 
$$
\H^{\cZ_3}(K)=(\supp\times\gp_3)_!\gq_3^*K\otimes \Qlb(\frac{1}{2})[1]^{\otimes \<\gamma,
2\check{\rho}\>} 
$$
For $K\in\D(\cZ_2)$ we have $\pi_{3,2}^*\H^{\cZ_2}(K)[\dim]\;\iso\; \H^{\cZ_3}(\pi_{3,2}^*K)[\dim]$, where $\dim$ is a funtion of a connected component of 
$\cZ_3$, namely the relative dimension of the corresponding component over $\cZ_2$. 
Since $\pi_{3,2}^*[\dim]$ is exact, it suffices to show that $\H^{\cZ_3}$ is right-exact.

 For $\bar d=(d_1,d_2)$ with $0\le d_1\le d_2$ denote by $^{\bar d}\cZ_3\subset \cZ_3$ the locally closed substack given by the condition that there exist a diagram
$$
\begin{array}{ccccc}
\bar L_1 & \subset & \bar L_2 & \subset & M\\
 \cup && \cup\\
L_1 & \subset & L_2,
\end{array}
$$
where $\bar L_k\subset M$ is a subbundle of rank $k$ with $\deg(\bar L_k/L_k)=d_k$. The stacks $^{\bar d}\cZ_3$ form a stratification of $\cZ_3$. 

 For $x\in X$ let $_x\cH_G\subset\cH_G$ denote the preimage of $x$ under $\supp: \cH_G\to X$. The following is straightforward. 
 
\begin{Lm} 
\label{Lm_intersection_is_contained}
For a $k$-point of $\cZ_3$ let $D$ be the effective divisor such that $s: \Omega^{-1}(D)\hook{} \cA\otimes\Sym^2 L_2^*$ is a subbundle. Then the fibre of $\gq_3: \cZ_3\times_{\Bun_G}\cH_G\to \cZ_3$ over this point is contained in 
$$
\mathop{\cup}\limits_{x\in \supp(D)} \cZ_3\times_{\Bun_G}{_x\cH_G} 
\eqno{\square}
$$
\end{Lm}

  Given $\bar d=(d_1,d_2)$ and $\bar d'=(d'_1,d'_2)$ denote by $^{\bar d,\bar d'}\cZ_3\subset \cZ_3\times_{\Bun_G}\cH_G$ the intersection 
$$
(\gp_3)^{-1}(^{\bar d}\cZ_3)\cap (\gq_3)^{-1}(^{\bar d'}\cZ_3)
$$ 
For $x\in X$ let $_x^{\bar d,\bar d'}\cZ_3$ denote the intersection of $^{\bar d,\bar d'}\cZ_3$ with $\cZ_3\times_{\Bun_G}{_x\cH_G}$.  Combining Lemma~\ref{Lm_intersection_is_contained}
with (\cite{G}, Lemma~7.11), we are reduced to the following statement.

\begin{Lm} For any $\bar d,\bar d'$ and $x\in X$ the sum of (the maximum of) the dimensions of fibres of maps in the diagram
\begin{equation}
\label{diag_dim_fibres}
\cZ_3 \;\getsup{\gp_3}\; {_x^{\bar d,\bar d'}\cZ_3} \;\toup{\gq_3}\; \cZ_3
\end{equation}
does not exceed $\<\gamma, 2\check{\rho}\>=3$. 
\end{Lm}
\begin{Prf}
A point of $_x^{\bar d,\bar d'}\cZ_3$ gives rise to the diagram
$$
\begin{array}{ccccc}
\bar L'_1 & \subset & \bar L'_2 & \subset & M'\\
\cup && \cup &&\cup \\
\bar L_1 & \subset & \bar L_2 & \subset & M\\
 \cup && \cup\\
L_1 & \subset & L_2,
\end{array}
$$
with $d_k=\deg(\bar L_k/L_k)$ and $d'_k=\deg(\bar L'_k/L_k)$. We must examine the cases:

\smallskip\noindent
1) $\bar d=\bar d'$. In this case a fibre of $\gq_3$ is a point, because $\bar L'_2$ generates a lagrangian subspace in $M'/M'(-x)$. A fibre of $\gp_3$ is 3-dimensional.

\smallskip\noindent
2) $d'_1=d_1$, $d'_2=d_2+1$.  Then a fibre of $\gq_3$ is 1-dimensional, becuase $M$ must contain $\bar L'_1$. A fibre of $\gp_3$ is 2-dimensional. 

\smallskip\noindent
3) $d'_1=d_1+1$, $d'_2=d_2+1$. Then $\bar L'_2=\bar L_2+\bar L'_1$ and $\bar L'_1=\bar L_1(x)$. A fibre of $\gq_3$ is 2-dimensional, a fibre of $\gp_3$ is 1-dimensional. 

\smallskip\noindent
4) $d'_1=d_1+1$, $d'_2=d_2+2$. Then $\bar L'_2=\bar L_2(x)$. A fibre of $\gp_3$ is a point, because $M'=M+\bar L'_2$. A fibre of $\gq_3$ is 3-dimenisonal.
\end{Prf}

\medskip\smallskip\noindent
7.5  {\scshape Hecke functors on $\cP_2$}

\smallskip\noindent
Recall the stack $^d\cP_2$ classifying collections: a modification of rank 2 bundles $(L_2\subset L'_2)$ on $X$ with $d=\deg(L'_2/L_2)$, $\cA\in\Bun_1$, and a section $s:\Omega^{-1}\hook{} \cA\otimes\Sym^2 L'^*_2$. Lemma~\ref{Lm_useful_in_fact} yields the equivalence of categories
$^dJ: \D(^d\cP_2)\,\iso\, \D^W(^d\cZ_2)$.

 We are going to define for $i=0,1,2$ the functors 
$$
_i\H^{\cP}: \D(^{d+i}\cP_2)\to \D(X\times{^d\cP_2})
$$
which, by construction, will satisfy the following property.

\begin{Pp} 
\label{Pp_Hecke_functors_cP}
Let $K\in\D^W(\cZ_2)$,  $^dK\in \D(^d\cP_2)$ and $^dF\in\D(X\times {^d\cP_2})$. Assume given for each $d$ isomorphisms  
$$
^dJ(^dK)\;\iso\; K\mid_{^d\cZ'_2}\;\;\; \mbox{and}\;\;\; 
^dJ(^dF)\;\iso\;\H^{\cZ_2}(K)\mid_{X\times {^d\cZ'_2}},
$$
where we used the $*$-restrictions. Then $^dF$ is an extension of objects $_i\H^{\cP}(^{d+i}K)$ ($i=0,1,2$) in the triangulated category $\D(X\times {^d\cP_2})$. More precisely, there exist distinguished triangles in $\D(X\times {^d\cP_2})$
$$
C\to {^dF}\to {_2\H^{\cP}(^{d+2}K)}\;\;\;\;\,\mbox{and}\,\;\;\;\;\; 
_0\H^{\cP}(^dK)\to C\to {_1\H^{\cP}(^{d+1}K)}
$$ 
\end{Pp}

\smallskip\noindent
7.5.1 Let $\delta_0: X\times {^d\cP_2}\to {^d\cP_2}$ be the map sending $(x\in X, L_2\subset L'_2, \cA, s)$ to $(L_2\subset L'_2, \cA(x), s')$, where $s'$ is the composition
 $$
 \Sym^2 L'_2\toup{s}\cA\otimes\Omega\hook{} \cA(x)\otimes\Omega
 $$
Set $_0\H^{\cP}(S)=\delta_0^*S$. Since $\delta_0$ is quasi-finite, $_0\H^{\cP}$ is right exact for the perverse t-structures.  
Consider the diagram 
$$
X\times {^d\cP_2} \; \getsup{\id\times \delta_0}\; X\times {^d\cP_2} \; \toup{\delta_2}\; {^{d+2}\cP_2},
$$
where $\delta_2$ sends $(x\in X, L_2\subset L'_2, 
\cA, s)$ to $(L_2\subset L'_2(x), \cA(2x), s)$. Note that $\id\times\delta_0$ is a closed immersion. Set
$$
_2\H^{\cP}(S)=(\id\times\delta_0)_!\delta_2^* S
$$ 
Since $\delta_2$ is quasi-finite, $_2\H^{\cP}$ is right exact for the perverse t-structures. 
 
Let $^d\cH_{\cP}$ denote the stack of collections: $\cA\in\Bun_1$, modifications of rank 2 vector bundles $L_2\subset L'_2\subset L''_2$ with $d=\deg(L'_2/L_2)$, where $L''_2/L'_2$ is a torsion sheaf of length one supported at $x\in X$, and a commutative diagram
\begin{equation}
\label{diag_H_cP}
\begin{array}{ccc}
\Sym^2 L''_2 & \to & \cA\otimes\Omega(x)\\
\cup && \cup\\
\Sym^2 L'_2 & \toup{s} &\cA\otimes\Omega
\end{array}
\end{equation}
with $s\ne 0$. 

 The existence of the latter diagram means that $L''_2/L'_2$ is an isotropic subspace of $L'_2(x)/L'_2$ equiped with the form $s: \Sym^2 (L'_2(x)/L'_2)\to (\cA(2x)/\cA(x))\otimes\Omega$.  
 
 We have the diagram
 $$
 X\times {^d\cP_2} \;\getsup{\supp\times\gp_{\cP}} \; {^d\cH_{\cP}} \;\toup{\gq_{\cP}} \;{^{d+1}\cP_2},
 $$
 where $\gp_{\cP}$ sends a point of $^d\cH_{\cP}$ to $(L_2\subset L'_2, \cA, s)$. The map 
 $\gq_{\cP}$ sends a point of $^d\cH_{\cP}$ to 
 $$
 (L_2\subset L''_2, \cA(x), s)
 $$ 
 The map $\supp: {^d\cH_{\cP}}\to X$ sends a point of $^d\cH_{\cP}$ as above to $x$. 
We set
$$
_1\H^{\cP}(S)=(\supp\times\gp_{\cP})_!\gq_{\cP}^*S[1](\frac{1}{2})
$$    
 
\begin{Prf}\select{of Proposition~\ref{Pp_Hecke_functors_cP}}

\smallskip\noindent
Recall the diagram we used to define the functor $\H^{\cZ_2}$
$$
X\times \cZ_2\,\getsup{\supp\times\gp_2} \,\cZ_2\times_{\Bun_G}\cH_G\, \toup{\gq_2} \,\cZ_2
$$
Given $d\le d'$ set $^{d,d'}\cZ_2=\gq_2^{-1}({^{d'}\!\cZ'_2})\cap \gp_2^{-1}({^d\cZ'_2})$. We will calculate the direct image under $\supp\times\gp_2$ with respect
to the corresponding stratification of $\cZ_2\times_{\Bun_G}\cH_G$.  Let $u,v$ denote the maps in the induced diagram
$$
X\times {^d\cZ'_2} \;\getsup{u}\; {^{d,d'}\cZ_2}\; \toup{v}\; {^{d'}\!\cZ'_2}
$$
Let $\tilde K$ be the $*$-restriction of $\gq_2^*K$ to $^{d,d'}\cZ_2$. A point of ${^{d,d'}\cZ_2}$ gives rise to the diagram
$$
\begin{array}{ccc}
L''_2 & \subset &  M'\\
\cup && \cup\\
L'_2 & \subset & M\\
\cup\\
L_2
\end{array}
$$
with $d=\deg(L'_2/L_2)$ and $d'=\deg(L''_2/L_2)$. We must examine three cases:

\smallskip\noindent
1) $d=d'$. Then $L''_2=L'_2$, and a fibre of $u$ admits a free transitive action of the geometric fibre $(\cA^{-1}\otimes\Sym^2 L'_2)_x$. The complex $\tilde K$ is constant along the fibres of $u$. So, $_0\H^{\cP}(^dK)$ is the contribution in $^dF$ of the stratum $^{d,d}\cZ_2$ 

\smallskip\noindent
2) $d'=d+2$. For a $k$-point of $X\times{^d\cZ'_2}$ we have $L''_2=L_2(x)$ and $M'=L''_2+M$ in the above diagram. If $\Sym^2 L'_2\to \cA\otimes\Omega$ does not factor through $\cA\otimes\Omega(-x)$ then the fibre of $u$ over this point is empty, otherwise this fibre is a point scheme. In the second case the extension $0\to \Sym^2 L''_2\to ?\to\cA(x)\to 0$ is the push-forward of $0\to (\Sym^2 L'_2)(x)\to ?\to \cA(x)\to 0$ under 
$$
(\Sym^2 L'_2)(x)\hook{} (\Sym^2 L'_2)(2x)
$$
So, the contribution of $^{d,d+2}\cZ_2$ in $^dF$ is $_2\H^{\cP}(^{d+2}K)$. 

\smallskip\noindent
3) $d'=d+1$. Fix a $k$-point of $X\times {^d\cZ'_2}$ and denote by $\bar Y$ the corresponding fibre of $u$. 
 
 Let $Y$ be the scheme of $L''_2$  such that $L'_2\subset L''_2\subset L'_2(x)$ gives rise to the diagram (\ref{diag_H_cP}). Note that $Y\,\iso\,\PP^1$ if the form on $L'_2(x)/L'_2$ is zero, $Y$ is a point if the kernel of the corresponding form is 1-dimensional, and $Y$ consists of two points if the form on $L'_2(x)/L'_2$ is non degenerate. 

   The fibres of the projection $\bar Y\to Y$ are isomorphic to $\AA^1$. More precisely, the 1-dimensional space $\cA_x^{-1}\otimes\Sym^2 (L'_2/L''_2(-x))$ acts on a fibre freely and transitively. 
   
    To see that the restriction $\tilde K\mid_{\bar Y}$ is constant along the fibres of $\bar Y\to Y$, note that the morphism $\cA^{-1}\otimes(\Sym^2 L''_2)(-x)\to\Omega$ factors as
$$
\cA^{-1}\otimes(\Sym^2 L''_2)(-x)\hook{}\cN\to\Omega,
$$ 
where $\cN$ is  the upper modification of $\cA^{-1}\otimes(\Sym^2 L''_2)(-x)$ defined by the 1-dimensional subspace $\cA^{-1}_x\otimes\Sym^2(L'_2/L''_2(-x))$ in the geometric fibre $\cA^{-1}\otimes\Sym^2(L''_2)_x$. 

 \smallskip
 
 It easily follows that the contribution of $^{d,d+1}\cZ_2$ in $^dF$ is $_1\H^{\cP}(^{d+1}K)$.
\end{Prf}  
 
\smallskip

\begin{Cor} 
\label{Cor_Hecke_strange}
Let $K\in\D^W(\cZ_2)$ and $^dK\in \D(^d\cP_2)$ equiped with isomorphisms
$^dJ(^dK)\,\iso\, K\mid_{^d\cZ'_2}$. Assume 
$$
\H^{\cZ_2}(K)\,\iso\, W\boxtimes K[3](\frac{3}{2})
$$ 
for a local system $W$ on $X$. Then for each $d$ the complex 
$W\boxtimes {^dK}[3](\frac{3}{2})$ is an extension of objects $_i\H^{\cP}(^{d+i}K)$ $(i=0,1,2)$ in the triangulated category $\D(X\times{^d\cP_2})$. \QED
\end{Cor} 
 
\bigskip\noindent
7.6 {\scshape Hecke functors on $\bar\cS$}

\medskip\noindent
Let $\bar\cS$ denote the stack classifying $L_2\in\Bun_2$, $\cA\in\Bun_1$ and an inclusion of coherent sheaves $s:\Omega^{-1}\hook{}\cA\otimes\Sym^2 L_2^*$. Define the following Hecke operators for $i=0,1,2$ 
$$
_i\H^{\bar\cS}: \D(\bar\cS)\to \D(X\times\bar\cS)
$$

 Let $\delta_0: X\times \bar\cS\to \bar\cS$ be the map sending $(x\in X, L_2, \cA, s)$ to $(L_2, \cA(x), s')$, where $s'$ is the composition
 $$
 \Sym^2 L_2\toup{s}\cA\otimes\Omega\hook{} \cA(x)\otimes\Omega
 $$
Set $_0\H^{\bar\cS}(K)=\delta_0^*K$. Since $\delta_0$ is quasi-finite, $_0\H^{\bar\cS}$ is right exact for the perverse t-structures.

Consider the diagram 
$$
X\times \bar\cS \; \getsup{\id\times \delta_0}\; X\times \bar\cS \; \toup{\delta_2}\; \bar\cS,
$$
where $\delta_2$ sends $(x\in X, L_2, 
\cA, s)$ to $(L_2(x), \cA(2x), s)$. Note that $\id\times\delta_0$ is a closed immersion. Set
$$
_2\H^{\bar\cS}(K)=(\id\times\delta_0)_!\delta_2^* K
$$ 
 
Let $\cH_{\bar\cS}$ denote the stack of collections: $\cA\in\Bun_1$, modifications of rank 2 vector bundles $L_2\subset L'_2$, with $\div(L'_2/L_2)=x$, and a commutative diagram
\begin{equation}
\label{diag_H_cS}
\begin{array}{ccc}
\Sym^2 L'_2 & \to & \cA\otimes\Omega(x)\\
\cup && \cup\\
\Sym^2 L_2 & \toup{s} &\cA\otimes\Omega
\end{array}
\end{equation}
with $s\ne 0$. The existence of the latter diagram means that $L'_2/L_2$ is an isotropic subspace of $L_2(x)/L_2$ equiped with the form $s: \Sym^2 (L_2(x)/L_2)\to (\cA(2x)/\cA(x))\otimes\Omega$.  
 
 We have the diagram
 $$
 X\times \bar\cS \;\getsup{\supp\times\gp_{\bar\cS}} \; \cH_{\bar\cS} \;\toup{\gq_{\bar\cS}} \;\bar\cS,
 $$
 where $\gp_{\bar\cS}$ sends a point of $\cH_{\bar\cS}$ to $(L_2, \cA, s)$. The map 
 $\gq_{\bar\cS}$ sends a point of $\cH_{\bar\cS}$ to 
 $(L'_2, \cA(x), s)$. 
 The map $\supp: \cH_{\bar\cS}\to X$ sends a point  (\ref{diag_H_cS}) to $x$. 
Set
$$
_1\H^{\bar\cS}(K)=(\supp\times\gp_{\bar\cS})_!\gq_{\bar\cS}^*K[1](\frac{1}{2})
$$    
 
\smallskip\noindent
7.6.1  Define the functor $F_{\bar \cS}: \D(\Bun_G)\to\D(\bar\cS)$ as follows. Given $K\in\D(\Bun_G)$ set 
\begin{equation}
\label{eq_K_1}
K_1=\alpha_{\cZ}^*K[\dimrel](\frac{\dimrel}{2}),
\end{equation}
where $\dimrel$ is the relative dimension of the corresponding connected component of $\cZ_1$ over $\Bun_G$. Let
$K_P$ denote the restriction of $K_1$ to the open substack $\Bun_P\subset\cZ_1$. Set
$$
F_{\bar\cS}(K)=\Four(K_P)\mid_{\bar\cS}
$$

\begin{Pp} 
\label{Pp_13_property}
Let $K\in\D(\Bun_G)$ be a Hecke eigen-sheaf corresponding to a $\check{G}$-local system $W_{\check{G}}$ on $X$. Set $F=F_{\bar\cS}(K)$. Consider the local systems $W=W^{\check{\omega}_1}_{\check{G}}$ and $W^0=W^{\check{\omega}_0}_{\check{G}}$.  
Then \\
1) there exist distinguished triangles in $\D(X\times\bar\cS)$
$$
C\to W\boxtimes F \to {_2\H^{\bar\cS}(F)[3](\frac{3}{2})}\;\;\;\;\;\mbox{and}\;\;\;\;\;
_0\H^{\bar\cS}(F)[-3](\frac{-3}{2})\to C\to {_1\H^{\bar\cS}(F)}
$$ 
2) For $\delta_2:X\times\bar\cS\to\bar\cS$ we have $\delta_2^*F\,\iso\, W^0\boxtimes F$. 
\end{Pp}
\begin{Prf}
1) Let $K_1\in\D(\cZ_1)$ be given by (\ref{eq_K_1}). Recall that a point of $^d\cZ_1$ is given by 
$$
(\cA\in\Bun_2, \; L_2\subset L'_2, \; 0\to\Sym^2 L'_2\to ?\to\cA\to 0)
$$ 
Let $\tau^P :{^d\cZ_1}\to\Bun_P$ be the map forgetting $L_2$. Calculation of dimensions shows that for the $*$-restriction 
$$
K_1\mid_{^d\cZ_1}\;\iso\;(\tau^P)^*K_P[3d](\frac{3d}{2})
$$
canonically. Recall that a point of $^d\cP_2$ is given by $(\cA\in\Bun_1, L_2\subset L'_2, s:\Sym^2 L'_2\to\cA\otimes\Omega)$. Let $\tau: {^d\cP_2}\to\bar\cS$ be the map forgetting $L_2$.

Set $K_2=WZ_{1,2}(K_1)$. For each $d$ set $^dK_2=\tau^*F_{\bar\cS}(K)[3d](\frac{3d}{2})$. Then for the $*$-restriction we have canonically 
$$
K_2\mid_{^d\cZ_2}\;\iso\; ^dJ(^dK_2)
$$
One easily checks that for $i=0,1,2$ we have canonical isomorphisms of functors 
$$
(\id\times\tau)^*\comp {_i\H^{\bar\cS}}\,\iso\, {_i\H^{\cP}}\comp\tau^*
$$ 
from $\D(\bar\cS)$ to $\D(X\times {^d\cP_2})$. 

 By Corolary~\ref{Cor_Hecke_strange}, for each $d$ the complex $W\boxtimes {^dK_2}[3](\frac{3}{2})$ is an extension of objects $_i\H^{\cP}(^{d+i}K_2)$ $(i=0,1,2)$ in $\D(X\times {^d\cP_2})$. Specifying to $d=0$, one gets the desired assertion. 
\end{Prf}

\bigskip\smallskip\noindent
7.7  {\scshape The stacks $^{rss}\cS\subset\cS\subset\bar\cS$}

\medskip\noindent
7.7.1  Let $\cS$ denote the stack classifying $L\in\Bun_2$, $\cC\in\Bun_1$ and a map $\Sym^2 L\to\cC$ inducing an inclusion of coherent sheaves $L\hook{} L^*\otimes\cC$. The map $\cS\to\bar\cS$ given by $L_2=L$, $\cA=\cC\otimes\Omega^{-1}$ is an open immersion.  

 Since the open substack $\cS\subset\bar\cS$ is defined by a condition at the generic point of $X$, the Hecke operators $_i\H^{\bar\cS}$ preserve this open substack, we denote the corresponding functors by
$$
_i\H^{\cS}: \D(\cS)\to\D(X\times\cS)   \;\;\;\; (i=0,1,2)
$$

 Let $\cS_d$ denote the union of those components of $\cS$ for which $2\deg\cC-2\deg L=d$. Note that $d\ge 0$ is even. 
 
  The nonramified two-sheeted Galois coverings $\tilde X\to X$ are in bijection with $\H^1_{et}(X, \ZZ/2\ZZ)$, and also in bijection with the isomorphism classes of pairs $(\cE,\kappa)$, where $\cE$ is a line bundle on $X$ and $\kappa: \cE^{\otimes 2}\,\iso\,\cO$.   

\begin{Lm} 
\label{Lm_stack_cS_0}
The stack $\cS_0$ classifies pairs: a two-sheeted nonramified covering $\tilde X\to X$ and a line bundle $\cB$ on $\tilde X$.  
\end{Lm}
\begin{Prf}
1) Given a Galois covering $\pi:\tilde X\to X$ of degree 2 and a line bundle $\cB$ on $\tilde X$ set $L=\pi_*\cB$. Let $\sigma$ be the nontrivial automorphism of $\tilde X$ over $X$. Let $\cE$ denote the anti-invariants of $\pi_*\cO_{\tilde X}$ under $\sigma$, so $\pi_*\cO_{\tilde X}\,\iso\, \cO_X\oplus\cE$. Note that $\pi^*\cE\,\iso\,\cO_{\tilde X}$ is equiped with the nontrivial descent data 
$$
\sigma^*\cO_{\tilde X}=\cO_{\tilde X}\toup{-1}\cO_{\tilde X}
$$  

 Set $\cC=(\det L)\otimes\cE$, so $\pi^*\cC$ identifies with $\cB\otimes\sigma^*\cB$ equiped with the natural descent data 
$$
\sigma^*(\cB\otimes\sigma^*\cB)=\cB\otimes\sigma^*\cB\toup{\id} \cB\otimes\sigma^*\cB
$$
We have canonically $\pi^*L\,\iso\, \cB\oplus\sigma^*\cB$, where $\sigma$ acts on $\cB\oplus\sigma^*\cB$ naturally. The projection
$$
\Sym^2(\cB\oplus\sigma^*\cB)\to \cB\otimes\sigma^*\cB
$$
with the natural descent data gives rise to a map $\Sym^2 L\to\cC$, which is a point of $\cS_0$. 

\medskip\noindent
2) On the other side, let $s: \Sym^2 L\to\cC$ be a point of $\cS_0$. Let $\tilde X\subset \PP(L)$ be the two-sheeted covering of $X$ whose fibre over $x\in X$ is the set of isotropic subspaces in $(L)_x$. Let $\cB$ be the line bundle on $\tilde X$ whose fibre at $V\subset L_x$ is $V$ itself. 
For $\pi:\tilde X\to X$ we get $\pi_*\cB\,\iso\, L$ canonically. 

 Let $\sigma$ be the nontrivial automorphism of $\tilde X$ over $X$. Let $\cE$ denote the $\sigma$-anti-invariants in $\pi_*\cO$. By 1), we have the symmetric form $\Sym^2 L\to \cE\otimes\det L$, which is also a point of $\cS_0$. Let $E$ denote the kernel of $\Sym^2 L\to \cE\otimes\det L$. Let us show that the composition 
$$
E\to \Sym^2 L\to\cC
$$ 
vanishes. It suffices to prove this after applying $\pi^*$, but $\pi^*E\,\iso\, \cB^{\otimes 2}\oplus \sigma^*\cB^{\otimes 2}$. So, we get a map $\tau_L$ included into the commutative diagram
$$
\begin{array}{ccc}
\Sym^2 L & \to & \cE\otimes\det L\\
\downarrow & \swarrow\lefteqn{\scriptstyle\tau_L}\\
\cC
\end{array}
$$
Since both symmetric forms on $L$ are everywhere nondegenerate, $\tau_L$ is an isomorphism.
\end{Prf}

 \medskip
  
\begin{Rem}
\label{Rem_symm_form_description}
 i) A version holds for a curve which may be not complete.\\
ii) If in the above lemma $\cB=\cO$ on $\tilde X$ then $\ZZ/2\ZZ$ acts on $\pi_*\cB=L$. So, $L\,\iso\, \cO\oplus\cE$, where $\cE$ is the line bundle of anti-invariants. The map $\Sym^2 L\to\cE\otimes\det L$ becomes $\cO\oplus\cE\oplus\cE^{\otimes 2}\to \cO$, it is given by $(1,0,-\kappa)$.  The curve $\tilde X$ can be recovered from $(\cE,\kappa)$ as $\{e\in\cE\mid \kappa(e^2)=1\}$. 
\end{Rem}  

\medskip\noindent
7.7.2  Let $^{rss}X^{(d)}\subset X^{(d)}$ be the open subscheme of divisors of the form $x_1+\ldots+x_d$ with $x_i$ pairwise distinct. Denote by $^{rss}\cS_d\subset \cS_d$ the preimage of $^{rss}X^{(d)}$ under the map $\cS_d\to X^{(d)}$ sending a point of $\cS_d$ to $\div(L^*\otimes\cC/L)$. Set 
$$
 \RCov^d=\Bun_1\times_{\Bun_1} {^{rss}X^{(d)}},
$$ 
where the map $^{rss}X^{(d)}\to \Bun_1$ sends $D$ to $\cO_X(-D)$, and the map $\Bun_1\to\Bun_1$ takes a line bundle to its tensor square.  

 It is understood that $^{rss}X^{(0)}=\Spec k$ and the point $^{rss}X^{(0)}\to\Bun_1$ is $\cO_X$. 
  
\begin{Pp} 
\label{Pp_structure_of_rssS}
The two-sheeted coverings $\pi:\tilde X\to X$ ramified exactly at $D\in {^{rss}X^{(d)}}$ (with $\tilde X$ assumed smooth) form an algebraic stack that can be identified with $\RCov^d$.

 The stack $^{rss}\cS_d$ classifies collections: $D\in {^{rss}X^{(d)}}$, a two-sheeted covering $\pi:\tilde X\to X$ ramified exactly at $D$, and a line bundle $\cB$ on $\tilde X$.
\end{Pp}
\begin{Prf}  
1) Given a two-sheeted (ramified) covering $\pi:\tilde X\to X$ and a line bundle $\cB$ on $\tilde X$ set $L=\pi_*\cB$. Let $\sigma$ be the nontrivial automorphism of $\tilde X$ over $X$.  
Let $x_1,\ldots,x_d\in X$ be the points of the ramification and $\tilde x_1,\ldots,\tilde x_d$ their preimages. 
 
  We have a canonical inclusion $\pi^* L\hook{} \cB\oplus \sigma^*\cB$, actually 
$$
\pi^* L=\{v\in \cB\oplus \sigma^*\cB\mid \;\mbox{the image of}\; v \;\mbox{in}\; (\cB\oplus\sigma^*\cB)_{\tilde x_i}\;\mbox{lies in}\; \cB_{\tilde x_i}\toup{\diag} (\cB\oplus\sigma^*\cB)_{\tilde x_i}\;\mbox{for all}\; i\}
$$ 
In particular, $ \pi^*(\det L)\,\iso\, \cB\otimes\sigma^*\cB(-\tilde x_1-\ldots-\tilde x_d)$. 

 Let $\cE$ denote the $\sigma$-anti-invariants in $\pi_*\cO$, so $\pi_*\cO\;\iso\;\cO\oplus\cE$. 
Clearly, $\pi^*\cE\,\iso\, \cO(-\tilde x_1-\ldots-\tilde x_d)$, and $\sigma$ acts on $\pi^*\cE$ as $-1$. This yields an isomorphism 
$$
\kappa: \cE^{\otimes 2}\,\iso\, \cO(-x_1-\ldots-x_d)
$$ 

 The diagram
$$
\begin{array}{ccc}
\pi^*\Sym^2 L & \subset & \cB^{\otimes 2}\oplus (\cB\otimes\sigma^*\cB) \oplus (\sigma^*\cB)^{\otimes 2}\\
&& \downarrow\\
\pi^*(\cE\otimes\det L) & \subset &  \cB\otimes\sigma^*\cB
\end{array}
$$
shows that $\pi^*\Sym^2 L\to \pi^*(\cE\otimes\det L)(2\tilde x_1+\ldots+2\tilde x_d)$ is regular and surjective. This map is compatible with the descent data, so gives rise to a regular surjective map
\begin{equation}
\label{form_s_ramified}
s: \Sym^2 L\to (\cE\otimes\det L)(x_1+\ldots+x_d)
\end{equation}
For each $x_i$ on the fibre $L_{x_i}=\cB/\cB(-2\tilde x_i)$ we get a symmetric form whose kernel is exactly $\cB(-\tilde x_i)/\cB(-2\tilde x_i)$. Further, $s$ induces an inclusion 
$$
L\hook{} (\cE\otimes L^*\otimes\det L)(x_1+\ldots+x_d)\,\iso\, L\otimes\cE(x_1+\ldots+x_d)
$$
and the quotient $(L\otimes\cE(x_1+\ldots+x_d))/L$ is of length $d$.

 For each $x_i$ there is a base in $L\otimes \hat\cO_{x_i}$ and in $\cE\otimes \hat\cO_{x_i}$ such that the matrix of $s: \Sym^2 (\hat\cO_{x_i}^2)\to \hat\cO_{x_i}(x_i)$ over the formal disk at $x_i$ becomes 
$$
\left(
\begin{array}{cc} 
t^{-1} & 0\\
0 & 1
\end{array}
\right),
$$
where $t\in\hat\cO_{x_i}$ is a local parameter. In other words, 
$$
(\cE\otimes L(x_1+\ldots+x_d))/L\;\iso\; \cO/\cO(-x_1)\oplus\ldots\oplus\cO/\cO(-x_d)
$$ 
  
\medskip\noindent
2) On the other side, let $s:\Sym^2 L\to\cC$ be a k-point of $\cS$ with $(L^*\otimes\cC)/L\,\iso\, \cO_{x_1}\oplus\ldots\oplus\cO_{x_d}$. Set $D=x_1+\ldots+x_d$. Note that $s$ is surjective. 
Let $\cI\subset \Sym^* L$ denote the homogeneous ideal generated by the image of 
$s^*: \cC^*\otimes(\det L)^2\hook{} \Sym^2 L$. Let $\tilde X\subset \PP(L)$ denote the closed subscheme given by $\cI$. Over $X-D$ this is exactly the curve of isotropic subspaces in $L$, as in Lemma~\ref{Lm_stack_cS_0}. Write $\pi:\tilde X\to X$ for the projection. 

 We claim that $\tilde X$ is smooth. To check this in the neighbourhood of $x_i$, pick a base $e_1,e_2$ in $L\otimes\hat\cO_{x_i}$ and $e\in\cC\otimes\hat\cO_{x_i}$ such that the matrix of $s:\Sym^2 (L\otimes\hat\cO_{x_i})\to \cC\otimes\hat\cO_{x_i}$ in these bases becomes 
$$
\left(
\begin{array}{cc} 
1 & 0\\
0 & t
\end{array}
\right),
$$
where $t\in\hat\cO_{x_i}$ is a local parameter. Then $\tilde X\times_X \Spec\hat\cO_{x_i}$ identifies with the closed subscheme $Y\hook{i}\PP^1\times\Spec\hat\cO_{x_i}$ given by $u_1^2+tu_2^2=0$, where $u_1,u_1$ are the homogeneous coordinates on $\PP^1$. The ring
$$
\hat\cO_{x_i}[u_1/u_2]/((u_1/u_2)^2+t)
$$ 
is a standard ramified extension of $\hat\cO_{x_i}$ of degree 2. The scheme $\tilde X\times_X \Spec\hat\cO_{x_i}$ is regular, so $\tilde X$ is smooth. Note that $\pi^{-1}(x_i)=:\tilde x_i$ are exactly the ramification points of $\pi$. 

 Let $\cB$ be the restriction of $\cO_{\PP(L)}(1)$ to $\tilde X$, it is equiped with $\pi^*L\to \cB$. Let us check that the induced map $L\to\pi_*\cB$ is an isomorphism. This is easy over $X-D$. 
Let $i: Y\hook{} \PP^1\times\Spec\hat\cO_{x_i}$ be as above, $\xi: \PP^1\times\Spec\hat\cO_{x_i}\to \Spec\hat\cO_{x_i}$ be the projection. We must check that 
$$
\xi_*\cO(1)\to \xi_*i_*i^*\cO(1)
$$ 
is an isomorphism. Define $\cV$ by the exact sequence $0\to\cV\to \cO(1)\to i_*i^*\cO(1)\to 0$ on $\PP^1\times\Spec\hat\cO_{x_i}$. It suffices to show that $\R^1\xi_*\cV=0$. But this is easily checked fibrewise over $\hat\cO_{x_i}$. 

 Let $\sigma$ be the nontrivial automorphism of $\tilde X$ over $X$. Let $\cE$ denote the $\sigma$-anti-invariants in $\pi_*\cO$, so $\pi_*\cO\,\iso\,\cO\oplus\cE$. By 1), we have
$\cE^{\otimes 2}\,\iso\,\cO(-x_1-\ldots-x_d)$ canonically, and $L$ is equiped with the form (\ref{form_s_ramified}). Define the vector bundle $E$ on $X$ by the exact sequence
$$
0\to E\to\Sym^2 L\toup{s} (\cE\otimes\det L)(x_1+\ldots+x_d)
$$
As in Lemma~\ref{Lm_stack_cS_0}, one checks that the composition $E\to\Sym^2 L\to\cC$ vanishes, and the induced map $(\cE\otimes\det L)(x_1+\ldots+x_d)\to\cC$ is an isomorphism.
\end{Prf}   

\bigskip\smallskip\noindent
7.8 {\scshape  Local version $\cS^{loc}$ of the stack $\cS$}

\medskip\smallskip\noindent
7.8.1 Set $\cO=k[[t]]$ and  $F=k((t))$. Let $\cS^{loc}$ denote the stack classifying: a free $\cO$-module $L$ of rank 2, a free $\cO$-module $\cC$ of rank 1, and a map $s:\Sym^2 L\to\cC$ inducing an inclusion $L\hook{} L^*\otimes\cC$.  

 Set $\Sym_+(\cO)=\{B\in \Mat_2(\cO)\mid {^tB}=B,\;  \det B\ne 0\}$. This is a $k$-scheme not of finite type. Further, $\GL(2,\cO)\times\cO^*$ is a group scheme over $k$ (not of finite type), and $\cS^{loc}$ identifies with the stack quotient of $\Sym_+(\cO)$ by the action of
$\GL(2,\cO)\times\cO^*$ given by $B\mapsto AB(^tA)\epsilon$, $(A,\epsilon)\in \GL(2,\cO)\times\cO^*$. 
 
 Given a $k$-point $(L,\cC, s)$ of $\cS^{loc}$, there exist bases $e_1,e_2\in L$ and $e\in\cC$ such that the matrix of $s$ in these bases is $\diag(t^a, t^b)$ 
for some $a\ge b\ge 0$ and 
$$
L^*\otimes\cC/L\,\iso\, \cO/t^a\cO\oplus\cO/t^b\cO
$$ 
It follows that two $k$-points $(L,\cC,s)$ and $(L',\cC',s')$ are isomorphic if and only if the $\cO$-modules $L^*\otimes\cC/L$ and $L'^*\otimes\cC'/L'$ are isomorphic. 
We identify the set of isomorphism classes of $k$-points of $\cS^{loc}$ with 
$$
\Phi=\{(a,b)\in\ZZ^2\mid a\ge b\ge 0\}
$$  

 For a closed point $x\in X$ a choice of an isomorphism $\cO\,\iso\,\hat\cO_{X,x}$ yields a map
$\cS\to\cS^{loc}$ given by the restriction of $(L,\cC,s)$ under $\Spec\hat\cO_{X,x}\to X$.

\medskip\noindent
7.8.2 \  Denote by $\Cov_F$ the $k$-stack associating to a scheme $S$ the groupoid of pairs $(S', \; \pi)$, where $S'$ is a scheme, and $\pi: S'\to S\times\Spec F$ is an \'etale covering of degree 2.  

  The stack $\Cov_F$ has (up to isomorphism) two $k$-points $(\Spec F',\pi)$, where the $F$-algebra $F'$ is one of the following
\begin{itemize}
\item $\;\; F'\;\iso\; k((t^{\frac{1}{2}}))$ (anisotropic case)
\item $\;\;  F'\;\iso\; F\oplus F$ (hyperbolic case)  
\end{itemize}    
 
  Given an $S$-point $(S',\pi)$ of $\Cov_F$, consider the rank 2 vector bundle $L=\pi_*\cO_{S'}$ on $S\times\Spec F$. Let $\sigma$ be the nontrivial automorphism of $S'$ over $S\times \Spec F$. We have $L=\cO_S\oplus\cE$, where $\cE$ denotes $\sigma$-anti-invariants in $L$. We have a canonical isomorphism $\kappa: \cE^{\otimes 2}\;\iso\; \cO_{S\times\Spec F}$. As in Remark~\ref{Rem_symm_form_description}, $L$ is equiped with a symmetric form 
$$
\begin{array}{ccc}
\Sym^2 L & \iso & \cO_{S\times\Spec F}\oplus\cE\oplus\cE^{\otimes 2}\\
\downarrow\lefteqn{\scriptstyle s} & \swarrow\lefteqn{\scriptstyle (1,0,-\kappa)}\\
\cO_{S\times\Spec F}, 
\end{array}
$$
The form $s$ in non degenerate, that is, induces an isomorphism $L\,\iso\, L^*$ of $\cO_{S\times\Spec F}$-modules. 

 For a $k$-point of $\Cov_F$, the symmetric form on $L$ is either hyperbolic or anisotropic, this explains our terminology (\cite{MVW}, ch.~1). In the anisotrpic case $\kappa(t^{\frac{1}{2}}\otimes t^{\frac{1}{2}})=t$. 

  It is easy to find a $\AA^1$-point of $\Cov_F$ such that over $\Gm\subset\AA^1$ we get the hyperbolic point of $\Cov_F$ and over $0\in\AA^1$ we get the anisotropic point. 
 
 We have a morphism of stacks $\cS^{loc}\to\Cov_F$ defined as follows. If $(L,\cC,\Sym^2 L\toup{s}\cC)$ is a $S$-point of $\cS^{loc}$ then we have an isomorphism of vector bundles $L\,\iso\, (L^*\otimes\cC)\mid_{S\times\Spec F}$ over $S\times\Spec F$. Define $S'\subset \PP(L)\mid_{S\times\Spec F}$ as the closed subscheme corresponding to the homogeneous ideal in $\Sym^2 (L\otimes F)$ generated by the image of 
$$
\cC^*\otimes(\det L)^2\otimes F \to \Sym^2 (L\otimes F)
$$ 
Then $\pi:S'\to S\times\Spec F$ is a point of $\Cov_F$.  
 
  The image of the $k$-point $(a,b)\in\Phi$ under $\cS^{loc}\to \Cov_F$ is anysotropic if $a-b$ is odd and hyperbolic otherwise. 
  
\bigskip\noindent
7.9  For $d,k\ge 0$ let 
$$
\xi_{d,k}:  {^{rss}\cS_{d, k}}\to {^{rss}\cS_d}
$$ 
denote the stack over $^{rss}\cS_d$ classifying: a point of $^{rss}\cS_d$ given by $(L,\cC, \Sym^2 L\toup{s}\cC)$, a subsheaf $L'\subset L$, where $L/L'$ is a torsion sheaf of length $k$ on $X$, such that the composition
$$
\Sym^2 L'\to\Sym^2 L\toup{s} \cC
$$
is surjective.

 We have a morphism of stack $^{rss}\cS_{d,k}\times X^{(m)}\to \cS$ sending $(L'\subset L,\cC, \Sym^2 L\toup{s}\cC, \, D'\in X^{(m)})$ to $(L',\cC(D'),s')$, where $s'$ is the composition 
$$
\Sym^2 L'\to\Sym^2 L\toup{s}\cC\hook{}\cC(D')
$$
 
\begin{Pp} The stacks $^{rss}\cS_{d, k}\times X^{(m)}$ form a stratification of $\cS$. 
\end{Pp}
\begin{Prf}
Recall that a point of $\cS$ is given by $(L,\cC,\Sym^2 L\toup{s}\cC)$. Let $\cS^0\subset \cS$ be the open substack given by the condition that $s$ is surjective. Sratifying $\cS$ by 
length of the cokernel of $s$, we are reduced to show that $^{rss}\cS_{d,k}$ form a stratification of $\cS^0$.
  
   Let $(L',\cC,s)$ be a $k$-point of $\cS^0$. Set $D=\div((L'^*\otimes\cC)/L')$ and write $D=\sum d_x x$. The restriction of $(L'^*\otimes\cC)/L'$ to $\Spec\hat\cO_{X,x}$ is isomorphic to $\cO/t^{d_x}_x\cO$, where $t_x$ is a local parameter at $x$. There is a unique subsheaf
$L'\subset L\subset L'^*\otimes\cC$ such that $s$ extends to a map $\Sym^2 L\to\cC$ yielding
$$
L'\subset L\subset L^*\otimes\cC\subset L'^*\otimes\cC,
$$
and 
$$
(L^*\otimes\cC)/L\mid_{\Spec\hat\cO_{X,x}}\,\iso\, \left\{
\begin{array}{rl}
0, & \mbox{if}\;\; d_x\;\;\mbox{is even}\\ 
k, & \mbox{if}\;\; d_x\;\; \mbox{is odd}
\end{array}
\right.
$$  
Our assertion follows.
\end{Prf}

\bigskip\noindent
7.10 {\scshape The stack $\cS_{\pi}$}

\medskip\noindent
Fix a $k$-point of $\RCov^d$ given by $D_{\pi}\in {^{rss}X^{(d)}}$ and $\pi:\tilde X\to X$ ramified exactly at $D_{\pi}$. 

 Given a point $(L,\cC, \Sym^2 L\toup{s}\cC)$ of $\cS$, set 
$$
D=\div(L^*\otimes\cC/L)
$$ 
and let $\pi_L: \tilde X_L\to X-D$ denote the corresponding two-sheeted covering defined as in 
Lemma~\ref{Lm_stack_cS_0}. Denote by $\cS_{\pi}$ the stack classifying: a point $(L,\cC, \Sym^2 L\toup{s}\cC)$ of $\cS$ together with an isomorphism over $X-D$
$$
\begin{array}{ccc}
\tilde X_L & \iso & \pi^{-1}(X-D)\\
\downarrow\lefteqn{\scriptstyle\pi_L} & \swarrow\\
X-D
\end{array}
$$
(note that $D_{\pi}$ does not intersect $X-D$, because $\pi_L$ is unramified). 
 
\medskip\noindent
7.10.1 Let $\tilde E$ be a rank one local system on $\tilde X$. We are going to define the category $\P^{\tilde E}(\cS_{\pi})$ of $\tilde E$-equivariant perverse sheaves on $\cS_{\pi}$. 

 Let $(X\times\cS_{\pi})^0\subset X\times\cS_{\pi}$ be the open substack of those $x\in X$, $(L,\cC, \Sym^2 L\to\cC)\in\cS_{\pi}$, for which the map $L\to L^*\otimes\cC$ is an isomorphism over the formal disk around $x\in X$. 
 
 Let $(\tilde X\times\cS_{\pi})^0$ denote the preimage of $(X\times\cS_{\pi})^0$ under 
 $$
 \pi\times\id: \tilde X\times\cS_{\pi}\to X\times\cS_{\pi}
 $$ 
Write $\cH_{\cS_{\pi}}$ for the stack classifying: a point of $(X\times\cS_{\pi})^0$ given by $(L,\cC,\Sym^2 L\toup{s}\cC)\in\cS_{\pi}$, $x\in X$ together with a commutative diagram 
$$
\begin{array}{ccc}
\Sym^2 L'  & \to & \cC(x)\\    
\cup && \cup\\
\Sym^2 L & \toup{s} & \cC,
\end{array}
$$ 
where $L\subset L'\subset L(x)$ is an upper modification of $L$ with $x=\div(L'/L)$. 

 We have a diagram
$$
(\tilde X\times\cS_{\pi})^0 \getsup{\supp\times \gp_{\cS_{\pi}}} \cH_{\cS_{\pi}} \toup{\gq_{\cS_{\pi}}} \cS_{\pi},
$$ 
where $\supp\times \gp_{\cS_{\pi}}$ sends a point of $\cH_{\cS_{\pi}}$ to  $(L,\cC,\Sym^2 L\toup{s}\cC)\in\cS_{\pi}$ together with the point $\tilde x\in\tilde X$ corresponding to the isotropic subspace $L'/L\subset L(x)/L$, so $\pi(\tilde x)=x$. Actually, $\supp\times \gp_{\cS_{\pi}}$ is an isomorphism. The map $\gq_{\cS_{\pi}}$ sends the above point to 
$$
(L',\cC(x), \Sym^2 L'\to\cC(x))\in\cS_{\pi}
$$  

The following is a version of the \select{Waldspurger category} that will be introduced in Sect.~8.  
\begin{Def}  Let $\P^{\tilde E}(\cS_{\pi})$ be the category, whose objects are pairs: a perverse sheaf $F$ on $\cS_{\pi}$ and an isomorphism 
$$
(\supp\times \gp_{\cS_{\pi}})_!\gq_{\cS_{\pi}}^*F\,\iso\, \tilde E\boxtimes F
$$ 
over $(\tilde X\times\cS_{\pi})^0$. The morphisms in $\P(\cS_{\pi})$ are the maps of the corresponding perverse sheaves compatible with the equivariance isomorphisms.
\end{Def}

\bigskip\smallskip
\centerline{\scshape 8. Waldspurger model for $\GL_2$}  
 
\medskip\smallskip\noindent 
8.1 Fix a $k$-point of $\RCov^d$ given by $D_{\pi}\in {^{rss}X^{(d)}}$ and $\pi:\tilde X\to X$ ramified exactly at $D_{\pi}$. Denote by $\sigma$ the nontrivial automorphism of $\tilde X$ over $X$, let $\cE$ be the $\sigma$-anti-invariants in $\pi_*\cO_{\tilde X}$.

Fix a $k$-point $x\in X$, write $\cO_x$ for the completed local ring of $X$ at $x$ and $F_x$ for its fraction field. Write $\tilde F_x$ for the \'etale $F_x$-algebra of regular functions on $\tilde X\times_X \Spec F_x$. If $x\in D_{\pi}$ then $\tilde F_x$ is anysotropic otherwise it is hyperbolic (cf. Sect.~7.8.2).  Denote by $\tilde\cO_x$ the ring of regular functions on $\tilde X\times_X \Spec \cO_x$. 
    
\begin{Def}  Let $\Wald_{\pi}^{x,loc}$ denote the stack classifying: a free $\cO_x$-module $L$ of rank 2, a free $\tilde F_x$-module $\cB$ of rank 1 together with an isomorphism $\xi: L\otimes_{\cO_x} F_x\,\iso\, \cB$ of $F_x$-modules.  
\end{Def}

  Let $\GL(\tilde F_x)$ denote the group of automorphisms of the $F_x$-linear vector space $\tilde F_x$, let $\GL(\tilde\cO_x)\subset \GL(\tilde F_x)$ be the stabilizor of $\tilde\cO_x$. 
Then  $\Wald_{\pi}^{x,loc}$ identifies with the stack quotient of the affine grassmanian 
$\Gr_{\tilde F_x}:=\GL(\tilde F_x)/\GL(\tilde\cO_x)$ by the group ind-scheme $\tilde F_x^*$.  
 
  A choice of a base in the free $\cO_x$-module $\tilde\cO_x$ yields isomorphisms $\GL(\tilde F_x)\,\iso\,\GL_2(F_x)$, $\GL(\tilde\cO_x)\,\iso\,\GL_2(\cO_x)$, and an inclusion $\tilde F_x^*\hook{} \GL_2(F_x)$. 

  For a $k$-point of $\Wald_{\pi}^{x,loc}$ consider the set of free $\tilde\cO_x$-submodules of rank one $\cB_{ex}\subset \cB$ such that $\xi(L)\subset \cB_{ex}$. This set contains a unique minimal element that we denote by $\cB_{ex}$. 
  
  In both split ($x\notin D_{\pi}$) and nonsplit ($x\in D_{\pi}$) case the isomorphism classes of $k$-points of $\Wald_{\pi}^{x,loc}$ are indexed by non negative integers $m\ge 0$, the corresponding point is given by $\deg(\cB_{ex}/L)=m$. Denote by $\Gr_{\tilde F_x}^m$ the $\tilde F_x^*$-orbit on $\Gr_{\tilde F_x}$ corresponding to $m\ge 0$. 

 In matrix terms, in the split case $\tilde\cO_x\,\iso\, \cO_x\oplus\cO_x$ has a distinguished (defined up to permutation) base $\{(1,0), (0,1)\}$ over $\cO_x$. This base yields an inclusion $\tilde F_x^*\hook{}\GL_2(F_x)$ whose image is the set of diagonal matrices. Then $\tilde F_x^*$-orbit on $\GL_2(F_x)/\GL_2(\cO_x)$ corresponding to $m\ge 0$ is given by the matrix 
$$\left(
\begin{array}{ccc}
 t^m & 1\\  
 0 & 1
 \end{array}
 \right),
 $$  
where $t\in\cO_x$ is a local parameter (cf. \cite{BFF}, Sect.~1). 

 In the nonsplit case the lattice $\cO_x\oplus \cO_xt^{m+\frac{1}{2}}\subset \tilde F_x$ is a representative for the $\tilde F_x^*$-orbit on $\Gr_{\tilde F_x}$ corresponding to $m\ge 0$. Here $t\in\cO_x$ is a local parameter. 

\medskip\smallskip\noindent 
8.2   In the same manner as in \cite{FGV} we can consider the following global model of $\Wald_{\pi}^{x,loc}$.

\begin{Def} Let $\Wald_{\pi}^x$ denote the stack classifying: a rank 2 vector bundle $L$ on $X$, a line bundle $\cB$ on $\pi^{-1}(X-x)$ and an isomorphism $L\,\iso\, \pi_*\cB$ over $X-x$.  
\end{Def}  
  
  As in Proposition~\ref{Pp_structure_of_rssS}, a point of $\Wald_{\pi}^x$ gives rise to a map
$$
s: \Sym^2 L\to (\cE\otimes\det L)(D_{\pi}+\infty x)
$$
Write  $\Wald_{\pi}^{x, \le m}\hook{} \Wald_{\pi}^x$ for the closed substack given by the condition that 
\begin{equation}
\label{map_Waldsp_glob}
s: \Sym^2 L\to (\cE\otimes\det L)(D_{\pi}+mx)
\end{equation} 
is regular. 

\begin{Lm} The stack $\Wald_{\pi}^{x, \le m}$ is algebraic, so $\Wald_{\pi}^x$ is an inductive limit of algebraic stacks.  
\end{Lm}
\begin{Prf}
Set $\ov{\RCov}^d=\Bun_1\times_{\Bun_1} X^{(d)}$, where the map $X^{(d)}\to\Bun_1$ sends $D$ to $\cO_X(-D)$ and $\Bun_1\to\Bun_1$ takes a line bundle to its tensor square. We have a map $\cS_d\to \ov{\RCov}^d$ sending $(L, \Sym^2 L\to \cC)$ to $(\det L)\otimes \cC^{-1}$ equiped with $(\det L)^{\otimes 2}\otimes \cC^{\otimes -2}\,\iso\, \cO(-D), \; D\in X^{(d)}$. 

 For $d=\deg D_{\pi}+2m$ consider the $k$-point $(\cE(-mx), (\cE(-mx))^{\otimes 2}\,\iso\, \cO(-D_{\pi}-2mx))$ of $\ov{\RCov}^d$. Then $\Wald^{x, \le m}_{\pi}$ is the fibre of $\cS_d\to \ov{\RCov}^d$ over this $k$-point. 
\end{Prf}

\medskip\smallskip

 Denote by 
\begin{equation}
\label{open_substack_Wald_x_m}
\Wald_{\pi}^{x, m}\hook{} \Wald_{\pi}^{x, \le m}
\end{equation}
the open substack given by the condition that (\ref{map_Waldsp_glob}) is surjective. 
The stack $\Wald_{\pi}^{x, m}$ classifies collections: a line bundle $\cB_{ex}$ on $\tilde X$, for which we set $L_{ex}=\pi_*\cB_{ex}$, and a lower modification $L\subset L_{ex}$ of vector bundles on $X$ such that the composition is surjective
$$
\Sym^2 L\to \Sym^2 L_{ex}\toup{s}\cC  
$$ 
and $\div(L_{ex}/L)=mx$. Here we have denoted $\cC=(\cE\otimes\det L_{ex})(D_{\pi})$, so 
$(L_{ex},\cC, \Sym^2 L_{ex}\toup{s}\cC)$ is the point of $^{rss}\cS$ corresponding to $\cB_{ex}$. 

 Another way to say is that the stratum $\Wald_{\pi}^{x, m}$ is given by fixing an extension of $\cB$ to a line bundle $\cB_{ex}$ on $\tilde X$ such that for $L_{ex}:=\pi_*\cB_{ex}$ we have
$L\subset L_{ex}$ and $\cB_{ex}$ is the smallest with this property. Then $L_{ex}/L\,\iso\,\cO_x/t^m$, where $t\in\cO_x$ is a local parameter. 
 
  Denote by $\pr_{\cW}: \Wald_{\pi}^{x,m}\to\Pic \tilde X$ the map sending the above point to $\cB_{ex}$. 
 
\medskip\noindent
8.2.1 Here is one more description. Denote by $(\Pic \tilde X)^x$ the scheme classifying a line bundle $\cB_{ex}$ on $\tilde X$ together with a trivialization $\cB\otimes\tilde\cO_x\,\iso\,\tilde\cO_x$. The group $\tilde\cO_x^*$ acts on $(\Pic \tilde X)^x$ by changing the trivialization. It is well-known that this action extends to an action of the group ind-scheme $\tilde F_x^*$ on $(\Pic \tilde X)^x$.

 Consider the action of $\tilde F_x^*$ on $(\Pic \tilde X)^x\times \Gr_{\tilde F_x}$ 
which is the product of natural actions on the factors. Then $\Wald_{\pi}^x$ identifies with the stack quotient of $(\Pic \tilde X)^x\times \Gr_{\tilde F_x}$ by $\tilde F_x^*$. Let $f_{\cW}: (\Pic \tilde X)^x\times \Gr_{\tilde F_x}\to \Wald_{\pi}^x$ be the corresponding map.  

\medskip\noindent
8.3  Fix a rank one local system $\tilde E$ on $\tilde X$. The $\tilde\cO_x^*$-orbits on $\Gr_{\tilde F_x}$ are finite-dimensional. So, we have the category of $\tilde\cO_x^*$-equivariant perverse sheaves on $\Gr_{\tilde F_x}$. 
 
\begin{Def} \select{Waldspurger category}  $\,\P^{\tilde E}(\Gr_{\tilde F_x})$ is the category of those $\tilde\cO_x^*$-equivariant perverse sheaves on $\Gr_{\tilde F_x}$ that 
\begin{itemize}
\item  (the nonsplit case) under the action of a uniformazer $\in\tilde F_x^*/\tilde\cO_x^*$ change by $\tilde E_{\tilde x}$, where $\pi(\tilde x)=x$.

\item (the split case) under the action of a uniformizer $t_{\tilde x}\in \tilde F_x^*/\tilde\cO_x^*$
change by $\tilde E_{\tilde x}$ for both $\tilde x\in\pi^{-1}(x)$.
\end{itemize}
\end{Def}

 One should be carefull about the following. Though $\P^{\tilde E}(\Gr_{\tilde F_x})$ is a full subcategory of the category $\P(\Gr_{\tilde F_x})$ of perverse sheaves on $\Gr_{\tilde F_x}$, the Ext groups in these two categories may be different. This is due to the fact that the $\tilde\cO_x^*$-orbits on $\Gr_{\tilde F_x}$ are not contractible. 

 Denote by $A\tilde E$ the automorphic local system on $\Pic\tilde X$ corresponding to $\tilde E$. For $d\ge 0$ its inverse image under $\tilde X^{(d)}\to\Pic ^d \tilde X$ identifies with the symmetric power $\tilde E^{(d)}$ of $\tilde E$. Define the perverse sheaf $\cW_m$ on $\Wald_{\pi}^x$ as the Goresky-MacPherson extension of 
$$
\pr_{\cW}^*A\tilde E\otimes\Qlb[1](\frac{1}{2})^{\otimes \dim \Wald_{\pi}^{x, m}}
$$
under (\ref{open_substack_Wald_x_m}). 
  
  For any $k$-point of $\Gr_{\tilde F_x}$ its stabilizor in $\tilde F_x^*$ is connected. So, the irreducible objects of $\P^{\tilde E}(\Gr_{\tilde F_x})$ are indexed by $m\ge 0$, the irreducibe object $\tilde\cW_m\in \P^{\tilde E}(\Gr_{\tilde F_x})$, defined up to a scalar automorphism, can be described by the following property: for the diagram
$$
\Wald_{\pi}^x\;\getsup{f_{\cW}}\; (\Pic\tilde X)^x\times\Gr_{\tilde F_x}\toup{\pr\times\id} \Pic\tilde X\times\Gr_{\tilde F_x} 
$$
we have $(\pr^*A\tilde E)\boxtimes \tilde\cW_m\,\iso\, f_{\cW}^*\cW_m$. 

 The group scheme $(\Pic\tilde E)^x$ acts on $\Wald_{\pi}^x$ as follows. The action map 
$$
\act: (\Pic\tilde E)^x\times \Wald_{\pi}^x\to \Wald_{\pi}^x
$$ 
sends $(\cB', \, \nu: \cB'\otimes\tilde\cO_x\,\iso\,\tilde\cO_x)\in (\Pic\tilde E)^x$ and $(\cB, L, \pi_*\cB\,\iso\, L\mid_{X-x})\in\Wald_{\pi}^x$ to
$$
(\cB\otimes\cB', \pi_*(\cB\otimes\cB')\,\iso\, L'\mid_{X-x})\in\Wald_{\pi}^x,
$$ 
where the vector bundle $L'$ on $X$ is the gluing of $\pi_*(\cB\otimes\cB')\mid_{X-x}$ and $L\mid_{\Spec\cO_x}$ via the isomorphism
$(\pi_*(\cB\otimes\cB'))\otimes F_x\,\iso\, L\otimes F_x$ induced by $\nu$.  

 Let $\P^{\tilde E}(\Wald_{\pi}^x)$ be the category of perverse sheaves on $\Wald_{\pi}^x$ 
that change by $\pr^*\tilde E$ under the action of $(\Pic\tilde E)^x$, where $\pr: (\Pic\tilde E)^x\to \Pic\tilde E$ is the projection. 

 Here is one more description of this category. Let
$$
\gq_{\Wald}: \pi^{-1}(X-x)\times\Wald_{\pi}^x\to\Wald_{\pi}^x
$$    
be the map sending $(\tilde x, \cB, \pi_*\cB\,\iso\, L\mid_{X-x})$ to $(\cB(\tilde x), \pi_*\cB(\tilde x)\,\iso\, L'\mid_{X-x})$, where the vector bundle $L'$ on $X$ is the gluing of $\pi_*\cB(\tilde x)\mid_{X-x}$ and $L\otimes\cO_x$ via the isomorphism $(\pi_*\cB(\tilde x))\otimes F_x\,\iso\, L\otimes F_x$, which is due to the fact that $\pi(\tilde x)\ne x$.  

 Then $\P^{\tilde E}(\Wald_{\pi}^x)$ is equivalent to the category of pairs: a perverse sheaf $F$  on $\Wald_{\pi}^x$ and an isomorphism $\gq_{\Wald}^*F\,\iso\, \tilde E\boxtimes F$.
  
  The irreducible objects of $\P^{\tilde E}(\Wald_{\pi}^x)$ are exactly $\cW_m$, $m\ge 0$. 
  
\medskip\noindent
8.4 Let $\Sph(\Gr_{\GL_2})$ be the category of $\GL_2(\cO_x)$-equivariant (spherical) perverse sheaves on the affine grassmanian $\Gr_{\GL_2}$. This is a tensor category equivalent to the category of representations of $\GL_2$ over $\Qlb$ (\cite{MV}). It acts on  $\D(\Wald_{\pi}^x)$ by Hecke functors as follows.  

 Let $_x\cH_{\GL_2}$ denote the Hecke stack classifying vector bundles $L,L'$ on $X$ together with an isomorphism $\beta: L\,\iso\,L'\mid_{X-x}$ over $X-x$. Consider the diagram
$$
\Wald_{\pi}^x \; \getsup{\gp_{\cW}}\; \Wald_{\pi}^x\times_{\Bun_2}\, {_x\cH_{\GL_2}} \,\toup{\gq_{\cW}}\, \Wald_{\pi}^x,
$$
where  $\gp_{\cW}$ sends a collection $(L,L',\beta, \cB, \pi_*\cB\,\iso\,L\mid_{X-x})$ to $(L, \cB, \pi_*\cB\,\iso\,L\mid_{X-x})$ and $\gq_{\cW}$ sends this collection to $(L', \cB, \pi_*\cB\,\iso\, L'\mid_{X-x})$.   

 Let $\Bun_2^x$ be the stack classifying $L\in\Bun_2$ together with its trivialization over $\Spec\cO_x$. The projection $\gq_{\GL_2}: {_x\cH_{\GL_2}}\to \Bun_2$ forgetting $L$ can be realized as a fibration  
$$
 \Bun_2^x\times_{\GL_2(\cO_x)} \Gr_{\GL_2}\to\Bun_2,
$$
so for $K\in \D(\Wald_{\pi}^x)$ and $\cA\in \Sph(\Gr_{\GL_2})$ we may form the corresponding twisted exterior product
$K\tboxtimes\cA$. It is normilized so that it is perverse for $K$ perverse and 
$$
\DD(K\tboxtimes\cA)\,\iso\, \DD(K)\tboxtimes\, \DD(\cA)
$$
Let $\H(\cA,\cdot): \D(\Wald_{\pi}^x)\to \D(\Wald_{\pi}^x)$ be the functor given by
$$
 \H(\cA, K)=(\gp_{\cW})_!(K\tboxtimes\cA)
$$
These functors are compatible with the tensor structure on $\Sph(\Gr_{\GL_2})$ in the sense that we have isomorphisms 
\begin{equation}
\label{iso_Hecke_compatibility}
\H(\cA_1, \H(\cA_2, K))\,\iso\,\H(\cA_1\ast\cA_2, K),
\end{equation}
where $\cA_1\ast\cA_2\in\Sph(\Gr_{\GL_2})$ is the convolution (cf. \cite{FGV}, Sect.~5).  One checks that $\P^{\tilde E}(\Wald_{\pi}^x)$ is preserved by Hecke functors.  
 
\begin{Th} 
\label{Th_4}
1) For $d\ge 0$ let $\lambda=(d,0)\in\Lambda^+_{\GL_2}$. We have a canonical isomorphism
$$
\H(\cA_{\lambda}, \cW_0)\,\iso\, \cW_d
$$
2) For $\lambda=(1,1)$ and $d\ge 0$ we have canonically
$$
\H(\cA_{\lambda}, \cW_d)\,\iso\, \left\{
\begin{array}{ll}
\cW_d\otimes \tilde E_{\tilde x}^{\otimes 2},& \mbox{the nonsplit case,} \; \,\pi(\tilde x)=x \\ \\
\cW_d\otimes \tilde E_{\tilde x_1}\otimes \tilde E_{\tilde x_2}, & \mbox{the split case,}\; \,\pi^{-1}(x)=\{x_1,x_2\} 
\end{array}
\right.
$$
\end{Th}

\medskip\noindent
8.5  Set $\Lambda^+_{\GL_2}=\{(a_1\ge a_2)\mid a_i\in\ZZ\}$. We view $\Lambda^+_{\GL_2}$ as the set of dominant coweights for $\GL_2$. For $\lambda=(a_1,a_2)\in\Lambda^+_{\GL_2}$ denote by $\Gr_{\tilde F_x}^{\lambda}\subset\Gr_{\tilde F_x}$ the locally closed subscheme classifying $\cO_x$-sublattices $L\subset t^{a_2}\tilde\cO_x$ such that 
$$
t^{a_2}\tilde\cO_x/L\,\iso\, \cO_x/t^{a_1-a_2}
$$ 
as $\cO_x$-modules. Let $\ov{\Gr}^{\lambda}_{\tilde F_x}$ denote the closure of $\Gr_{\tilde F_x}^{\lambda}$ in $\Gr_{\tilde F_x}$.
 
 Our proof of Theorem~\ref{Th_4} is inspired by (\cite{FGV}, Theorem~4), the following is a key point.  

\begin{Pp} For $m\ge 0$ and a dominant coweight $\lambda=(a_1\ge a_2)$ of $\GL_2$ the intersection $\ov{\Gr}^{\lambda}_{\tilde F_x}\cap \Gr^m_{\tilde F_x}$ is non empty iff $\,0\le m\le a_1-a_2$ and has pure dimension $m$. 
\end{Pp}
\begin{Prf}
1) (the split case). Use the matrix realization of $\Gr_{\tilde F_x}$ as in Sect.~8.1. Using the action of the center of $\GL_2$, we may reduce to the case $\lambda=(a,0)$. 
Stratify $\ov{\Gr}^{\lambda}_{\GL_2}$ by intersecting with $N(F_x)$-orbits on the affine grassmanian, where $N\subset \GL_2$ is the standard maximal unipotent subgroup. 
For all strata the argument is the same, let us explain it for the open stratum 
$$
\left\{\left(
\begin{array}{cc}
t^a & b\\ 
0 & 1
\end{array}
\right), b\in\cO_x\right\}/\left\{\left(
\begin{array}{cc}
1 & c\\
0 & 1
\end{array}
\right),
c\in\cO_x\right\},
$$
which we identify with $\cO_x/t^a$ via the map $\left(\begin{array}{cc}
t^a & b\\ 
0 & 1
\end{array}
\right)\mapsto b$. The point $b\in\cO/t^a$ lies in $\Gr^m_{\tilde F_x}$ iff $b\in t^{a-m}\cO_x^*$. 

\smallskip\noindent
2) (the nonsplit case). Let $t\in\cO_x$ be a local parameter. Multiplying by an appropiate power of $t$ we are reduced to the case $\lambda=(a,0)$. Then $\ov{\Gr}^{\lambda}_{\tilde F_x}$ is the scheme of $\cO_x$-sublattices $L\subset \tilde\cO_x$ such that $\dim(\tilde\cO_x/L)=d$. The intersection $\ov{\Gr}^{\lambda}_{\tilde F_x}\cap \Gr^m_{\tilde F_x}$ is then the scheme of 
sublattices 
$$
L\subset t^{\frac{1}{2}(a-m)}\tilde\cO_x\subset \tilde\cO_x
$$
such that $\dim(t^{\frac{1}{2}(a-m)}\tilde\cO_x)/L=m$ and $L\nsubseteq t^{\frac{a-m+1}{2}}\tilde\cO_x$. Our assertion follows.
\end{Prf}
 
\medskip 
\begin{Rem} In the nonsplit case the schemes $\ov{\Gr}^{\lambda}_{\tilde F_x}\cap \Gr^m_{\tilde F_x}$ are connected, whence in the split case they admit several connected components.
\end{Rem}

 Actually, we need the following a bit different result. Given $\lambda\in\Lambda^+_{\GL_2}$ and $\cO_x$-lattice $L\subset \tilde F_x$, denote by $\ov{\Gr}^{\lambda}_{\tilde F_x}(L)\subset \Gr_{\tilde F_x}$ the closed subscheme of $\cO_x$-lattices $L'\subset \tilde F_x$ such that 
$$
(L', L, \, L\otimes F_x\,\iso\, L'\otimes F_x)\in \ov{\Gr}^{\lambda}_{\GL_2}(L)
$$ 
More precisely, for any isomorphism $L\,\iso\,\cO_x\oplus\cO_x$ of $\cO_x$-modules the corresponding point 
$$
(L', L'\otimes F_x\,\iso\, F_x\oplus F_x)\in\ov{\Gr}^{\lambda}_{\GL_2}
$$ 

\begin{Pp}  
\label{Pp_17}
Let $m\ge 0$ and $L\subset \tilde F_x$ be a $\cO_x$-lattice lying in $\Gr^m_{\tilde F_x}$. For a dominant coweight $\lambda=(d, 0)$ of $\GL_2$ the intersection 
$\ov{\Gr}^{\lambda}_{\tilde F_x}(L)\cap \Gr^0_{\tilde F_x}$ is empty unless $d\ge m$. For $d\ge m$  it is a point (resp., a union of $d-m+1$ points) in the nonsplit case (resp., in the split case).
\end{Pp}
\begin{Prf} 
1) (the nonsplit case). Multiplying by a suitable element of $\tilde F_x^*$, we may assume $L=\cO_x\oplus \cO_x t^{m+\frac{1}{2}}$. The scheme  $\ov{\Gr}^{\lambda}_{\tilde F_x}(L)$ classifies $\cO_x$-sublattices $L'\subset L$ such that $\dim(L/L')=d$. A point $L'$ of this scheme lies in $\Gr^0_{\tilde F_x}$ if and only if $L'=t^{\frac{d+m}{2}}\tilde \cO_x$.
Our assertion follows.

\smallskip\noindent
2) (the split case). Choose a base $\{e_1,e_2\}$ in $\tilde\cO_x$ over $\cO_x$ as in 8.1. Multiplying by an suitable element of $\tilde F_x^*$ we may assume 
$$
L=t^m\cO_xe_1\oplus \cO_x(e_1+e_2)
$$ 
The scheme  $\ov{\Gr}^{\lambda}_{\tilde F_x}(L)$ classifies $\cO_x$-sublattices $L'\subset L$ such that $\dim(L/L')=d$. A point $L'$ of this scheme lies in $\Gr^0_{\tilde F_x}$ if and only if $L'=t^{a_1}\cO_xe_1\oplus t^{a_2}\cO_x e_2$ for some $a_1,a_2\ge 0$ such that $d+m=a_1+a_2$. So, the intersection in question identifies with the set of pairs $\{(a_1,a_2)\mid a_i\ge m,\; d+m=a_1+a_2\}$.    
\end{Prf} 
 
\bigskip\noindent
\begin{Prf}\select{ of Theorem~\ref{Th_4}} \   2) is easy and left to the reader.

\smallskip\noindent
1) We change the notation letting $\lambda=(0,-d)\in\Lambda^+_{\GL_2}$ for given $d\ge 0$. We will establish canonical isomorphisms 
$$
\H(\cA_{\lambda}, \cW_0)\,\iso\, \left\{
\begin{array}{ll}
\cW_d\otimes \tilde E_{\tilde x}^{\otimes -2d},& \mbox{the nonsplit case,} \; \,\pi(\tilde x)=x \\ \\
\cW_d\otimes \tilde E_{\tilde x_1}^{\otimes -d}\otimes \tilde E_{\tilde x_2}^{\otimes -d}, & \mbox{the split case,}\; \,\pi^{-1}(x)=\{x_1,x_2\} 
\end{array}
\right.
$$
Denote by $K_m$ the $*$-restriction of $\H(\cA_{\lambda},\cW_0)$ to $\Wald_{\pi}^{x,m}$. Since $\Wald_{\pi}^{x,0}$ is closed in $\Wald_{\pi}^x$ and $\cW_0$ is self-dual (up to replacing $\tilde E$ by $\tilde E^*$), our assertion is reduced to the following lemma.
\end{Prf}

\begin{Lm} 
\label{Lm_10}
We have $K_m=0$ unless $m\le d$. The complex $K_m$ is placed in non positive (resp., strictly negative) perverse degrees for $m=d$ (resp., for $m<d$). We have canonically
$$
K_d\,\iso\, (\pr_{\cW}^*A\tilde E)\otimes R\otimes\Qlb[1](\frac{1}{2})^{\otimes \dim\Wald_{\pi}^{x,d}}, 
$$
where $\pr_{\cW}:\Wald_{\pi}^{x,d}\to\Pic \tilde X$ is the projection and 
$$
R\,\iso\,  \left\{
\begin{array}{ll}
\tilde E_{\tilde x}^{\otimes -2d},& \mbox{the nonsplit case,} \; \,\pi(\tilde x)=x \\ \\
\tilde E_{\tilde x_1}^{\otimes -d}\otimes \tilde E_{\tilde x_2}^{\otimes -d}, & \mbox{the split case,}\; \,\pi^{-1}(x)=\{x_1,x_2\} 
\end{array}
\right.
$$
\end{Lm}
\begin{Prf}
Consider a point $\eta=(\cB_{ex}, L\subset L_{ex}=\pi_*\cB_{ex})$ of $\Wald_{\pi}^{x,m}$, so $mx=\div(L_{ex}/L)$. Write $_x\ov{\cH}_{\GL_2}^{\lambda}$ for the closed substack of $_x\cH_{\GL_2}$ that under the projection $\gq_{\GL_2}$ identifies with 
$$
\Bun_2^x\times_{\GL_2(\cO_x)}\ov{\Gr}_{\GL_2}^{\lambda}\to \Bun_2
$$
Choose a trivialization of $\cB_{ex}$ over $\Spec\tilde\cO_x$. 
The fibre of 
$\gp_{\cW}: \Wald_{\pi}^x\times_{\Bun_2}\, {_x\ov{\cH}^{\lambda}_{\GL_2}}\to  \Wald_{\pi}^x$ 
over $\eta$ identifies with $\ov{\Gr}^{-w_0(\lambda)}_{\GL_2}(L)$, where we have set $-w_0(\lambda)=(d,0)$. For the diagram
$$
\Wald_{\pi}^x\;\getsup{\gp_{\cW}}\;\Wald_{\pi}^x\times_{\Bun_2} \,  {_x\ov{\cH}^{\lambda}_{\GL_2}}
\,\toup{\gq_{\cW}}\,  \Wald_{\pi}^x
$$
we get $\H(\cA_{\lambda}, \cdot)=(\gp_{\cW})_!\gq_{\cW}^*(\cdot)[d](\frac{d}{2})$. Only the stratum 
$$
\ov{\Gr}^{-w_0(\lambda)}_{\GL_2}(L)\cap \Gr^0_{\tilde F_x}
$$ 
contributes to $K_m$. By Proposition~\ref{Pp_17}, for $m=d$ this is a point whose image under $\gq_{\cW}$ is 
$$
L'=\left\{
\begin{array}{ll}
\pi_*(\cB_{ex}(-2d\tilde x)),& \mbox{the nonsplit case,} \; \pi(\tilde x)=x\\
\pi_*(\cB_{ex}(-d\tilde x_1-d\tilde x_2)), & \mbox{the split case,}\; \pi^{-1}(x)=\{\tilde x_1, \tilde x_2\}
\end{array}
\right.
$$
Since $\dim\Wald_{\pi}^{x,m}=m+\dim\Pic\tilde X$, our assertion follows from the automorphic property of $A\tilde E$. Namely, for the map $m_{\tilde x}:\Pic\tilde X\to\Pic\tilde X$ sending $\cB$ to $\cB(\tilde x)$ we have canonically $m^*_{\tilde x} A\tilde E\,\iso\, A\tilde E\otimes \tilde E_{\tilde x}$. 
\end{Prf}

\medskip
\begin{Rems} i) Our proof of Theorem~\ref{Th_4} also shows the following. The stratum $\Wald_{\pi}^{x,d}$ is dense in $\Wald_{\pi}^{x,\le d}$. Besides, $\cW_d[-\dim\Wald_{\pi}^{x,d}]$ is a constructable sheaf on $\Wald_{\pi}^x$ placed in usual cohomological degree zero. Its fibres over points of $\Wald_{\pi}^{x,m}$ are 1-dimensional (resp., $d-m+1$-dimensional) in the non split (resp., split) case for $m\le d$. 

\smallskip\noindent
ii) The category $\P^{\tilde E}(\Wald_{\pi}^x)$ is not semisimple. Indeed, for $\lambda=(0,-1)$ consider the finite map 
$$
\gq_{\cW}: \Wald_{\pi}^{x,0}\times_{\Bun_2}\, {_x\ov{\cH}^{\lambda}_{\GL_2}}\to \Wald_{\pi}^{x,\le 1}
$$ 
It is an isomorphism over the open substack $\Wald_{\pi}^{x,1}$. Since the open immersion
$\Wald_{\pi}^{x,1}\hook{} \Wald_{\pi}^{x,0}\times_{\Bun_2}\, {_x\ov{\cH}^{\lambda}_{\GL_2}}$ is affine, the open immersion $\Wald_{\pi}^{x,1}\hook{}\Wald_{\pi}^{x,\le 1}$ is also affine. 
Let $\cW_{m,!}$ denote the $!$-extension of $\cW_m\mid_{\Wald_{\pi}^{x,m}}$ under (\ref{open_substack_Wald_x_m}). Then $\cW_{1,!}\in\P^{\tilde E}(\Wald_{\pi}^x)$. So, if this category was semisimple, the exact sequence of perverse sheaves
$$
0\to K\to \cW_{1,!}\to \cW_1\to 0
$$
would split, which contradict the fact that the $*$-restriction $\cW_1\mid_{\Wald_{\pi}^{x,0}}$ is non zero.
\end{Rems}

\medskip\noindent
8.6 {\scshape Casselman-Shalika formula} \   For $\lambda\in\Lambda^+_{\GL_2}$ write $U^{\lambda}$ for the irreducible representation of the Langlands dual group $\check{\GL}_2$ over $\Qlb$. Let $E$ be a $\check{\GL}_2$-local system on $X$ equiped with an isomorphism
$$
U^{(1,1)}_E\,\iso\, \left\{
\begin{array}{ll}
\tilde E_{\tilde x}^{\otimes 2},& \mbox{the nonsplit case,} \; \,\pi(\tilde x)=x \\ \\
\tilde E_{\tilde x_1}\otimes \tilde E_{\tilde x_2}, & \mbox{the split case,}\; \,\pi^{-1}(x)=\{x_1,x_2\} 
\end{array}
\right.
$$
We associate to $E$ the ind-object $K_E$ of $\P^{\tilde E}(\Wald^x_{\pi})$ given by
$$
K_E=\mathop{\oplus}\limits_{d\ge 0} \cW_d\otimes U^{(0,-d)}_E
$$
For a representation $U$ of $\check{\GL}_2$ write $\cA_U$ for the object of $\Sph(\Gr_{{\GL}_2})$ corresponding to $V$ via the Satake equivalence $\Rep(\check{\GL}_2)\,\iso\, \Sph(\Gr_{{\GL}_2})$.  
 
  One formally derives from Theorem~\ref{Th_4} the following.
  
\begin{Cor} For any $U\in\Rep(\check{\GL}_2)$ there is an isomorphism $\alpha_U: \H(\cA_U, K_E)\,\iso\, K_E\otimes U_E$. For $U,U'\in\Rep(\check{\GL}_2)$ the diagram commutes
$$
\begin{array}{ccc}
\H(\cA_{U'}, \H(\cA_U, K_E)) & \toup{\alpha_U} & \H(\cA_{U'}, K_E\otimes U_E)\\ 
\downarrow\lefteqn{\scriptstyle \gamma} && \downarrow\lefteqn{\scriptstyle\alpha_{U'}\otimes\id}\\
\H(\cA_{U\otimes U'}, K_E) & \toup{\alpha_{U\otimes U'}} & K_E\otimes (U\otimes U')_E,
\end{array}
$$
where $\gamma$ is the isomorphism (\ref{iso_Hecke_compatibility}). 
\end{Cor} 
 
\begin{Rem} One may view $\Gr_{\tilde F_x}$ as the ind-scheme classifying a rank 2 vector bundle $L$ on $X$ together with an isomorphism $L\,\iso\, \pi_*\cO_{\tilde X}\mid_{X-x}$. This yields a natural map $\Gr_{\tilde F_x}\to\Wald_{\pi}^x$. 

 The results of Sect.~8 hold also in the case of a finite base field $k=\Fq$. In this case we have the Waldpurger module $W\! A_{\chi}$ introduced in 1.4. For $d\ge 0$ consider the function trace of Frobenius of $\cW_d$ on $\Wald_{\pi}^x(k)$, let $W_d$ be its restriction to $\Gr_{\tilde F_x}$. Then $\{W_d, \, d\ge 0\}$ is a base of the vector space $W\! A_{\chi}$. 
 
  The space $W\! A_{\chi}$ also has the base (indexed by $d\ge 0$) consisting of functions supported over the $\tilde F_x^*$-orbit corresponding to $d$. The Casselman-Shalika formula in this base is given by (\cite{BFF}, Theorem~1.1), it involves some nontrivial denominators. This corresponds to the fact that our ind-object $K_E$ is not locally finite on $\Wald_{\pi}^x$.
\end{Rem} 

\appendix
\bigskip\medskip
\centerline{\scshape Appendix A. Fourier transforms}

\medskip\noindent
For the convenience of the reader, we collect some well-known observations about equivariant categories and Fourier transforms that we need. The proofs are omitted.

\medskip\noindent
A.1 Let $S$ be a scheme of finite type and $\pr: G\to S$ be a groupoid. Assume that $\pr$ is of finite type, with contractible fibres and smooth of relative dimension $k$.  
Assume also that $\act: G\to S$ is smooth of relative dimension $k$. Let $\cL$ be a local system on $G$ whose restriction to the unit section $S\to G$ is trivialized.

 By (\cite{G}, Lemma~4.8), 
we have the Serre subcategory $\P^W(S)\subset \P(S)$ of perverse sheaves $K\in \P(S)$ such that there exists an isomorphism $\act^*K\otimes\cL\,\iso\,\pr^*K$ whose restriction to the unit section is the identity. Let $\D^W(S)\subset \D(S)$ denote the full triangulated subcategory generated by $\P^W(S)$. 

 We write $\D^W_{\cL}(S)$ if we need to express the dependence on $\cL$. For $K\in\D(S)$ we have
 $K\in \D^W_{\cL}(S)$ if and only if $\DD(K)\in \D^W_{\cL^*}(S)$.

 Let $\beta: S'\to S$ be an $S$-scheme of finite type. The groupoid $G$ "lifts" to $S'$ if we have two cartesian squares
$$
\begin{array}{ccc}
G & \toup{\pr}& S\\
\uparrow\lefteqn{\scriptstyle \beta'} && \uparrow\lefteqn{\scriptstyle \beta}\\
G' & \toup{\pr'}& S'
\end{array}
$$
and
$$
\begin{array}{ccc}
G & \toup{\act} & S\\
\uparrow\lefteqn{\scriptstyle \beta'} && \uparrow\lefteqn{\scriptstyle \beta}\\
G' & \toup{\act'} & S'
\end{array}
$$
that make $G'$ a groupoid over $S'$. 

 For the local system $\beta'^*\cL$ we get the category $\D^W(S')$. The functors $\beta_!$ and $\beta_*$ send $\D^W(S')$ to $\D^W(S)$. The functors $\beta^*$ and $\beta^!$ send $\D^W(S)$ to $\D^W(S')$.

\medskip\noindent
A.2  Let $Y\to Z$ be a morphism of schemes of finite type
and $E\to Z$ be a vector bundle over $Z$. Assume that $E$ acts on $Y$ over $Z$, and $\act:E\times_Z Y\to Y$ is smooth of relative dimension $\rk E$.  We have a natural pairing $\chi: E^*\times_Z E\times_Z Y\to \A^1$. For the local system $\cL=\chi^*\cL_{\psi}$ we get the category $\D^W(E^*\times_Z Y)$ as in A.1.

 Let  $F:\D(Y)\to \D^W(E^*\times_Z Y)$ be the functor 
$$
F(K)=\Four(\act^*K)[\rk E](\frac{\rk E}{2})
$$
 Then $F$ is an equivalence of triangulated categories, t-exact and commutes with Verdier duality (up to replacing $\psi$ by $\psi^{-1}$). The quasi-inverse functor is given by $K\mapsto \pr_!(K)$, where $\pr: E^*\times_Z Y\to Y$ is the projection.

Moreover, for any $K\in \D^W(E^*\times_Z Y)$ the natural map
$\pr_!(K)\,\iso\,\pr_*(K)$ is an isomorphism.

\medskip\noindent
A.3 Suppose we are in the situation of A.2. Assume in addition that $p: E'\to E$ is a 
morphism of vector bundles over $Z$. Then $E'$ also acts on $Y$ over $Z$ (via $E$), and we
have the functor $F': \D(Y)\to \D^W(E'^*\times_Z Y)$ defined as in A.2.

 Then we have an isomorphism of functors $F'\;\iso\; (\check{p}\times\id)_! \comp F$,
where $\check{p}\times\id: E^*\times_Z Y\to E'^*\times_Z Y$ is the dual map (cf. \cite{G}, 
5.16).

\bigskip\noindent
{\scshape Acknowlegments.} The author is grateful to G. Laumon for constant support and
thanks D.~Gaitsgory and E. Frenkel for useful discussions.


\begin{thebibliography}{99}
\bibitem{BD} A. Beilinson, V. Drinfeld, Quantization of Hitchin's integrable system and
Hecke eigen-sheaves, preprint
\bibitem{BG} A. Braverman, D. Gaitsgory, Geometric Eisenstein series, Invent. Math.  150 
(2002),  no. 2, 287--384.
 
\bibitem{BFF} D. Bump, S. Friedberg, M. Furusawa, Explicit formulas for the Waldspurger and Bessel models, math.RT/9410202 
 
\bibitem{FGV} E. Frenkel, D. Gaitsgory, K. Vilonen, Whittaker patterns in the geometry of moduli spaces of bundles on curves,  Ann. of Math. (2) 153 (2001), no. 3, 699--748. 
 
\bibitem{G} D. Gaitsgory, On a vanishing conjecture appearing in the geometric Langlands correspondence, math.AG/0204081
\bibitem{Ly} S. Lysenko, On automorphic sheaves on $\Bun_G$, math.RT/0211067

\bibitem{MV} I. Mirkovic, K. Vilonen, Geometric Langlands duality and representations of algebraic groups over commutative rings, math.RT/0401222 

\bibitem{MVW} C. Moeglin, M.-F. Vigneras, J.L. Waldspurger, Correspondence de Howe
sur un corps $p$-adique, Lecture Notes in Math. 1291 (1987)  

\bibitem{PSh} I.I.Piatetski-Shapiro, On the Saito-Kurokawa lifting, Invent. Math. 71 (1983),
309-338

\bibitem{R} Rallis, S. On the Howe duality conjecture.  Compositio Math.  51  (1984), 
no. 3, 333--399.
\bibitem{V} J.-L. Verdier, Des cat\'egories d\'eriv\'ees des cat\'egories ab\'eliennes,
Asterisque 239 (1996), Soci\'et\'e Math. de France

\bibitem{W} J.-L. Waldspurger, Sur les valeurs de certaines fonctions L automorphe en leur centre de sym\'etrie, Compositio Math. 54 (1985), 173-242 

\end{thebibliography}
\end{document}